\documentclass{scrartcl}

\usepackage{amsmath}
 \numberwithin{equation}{section}
\usepackage{amssymb}
\usepackage{amsthm}

\theoremstyle{definition}
 
\theoremstyle{remark}
 
\usepackage[main=english]{babel}
\usepackage[noisbn,noissn,fixlanguage]{babelbib}
 \setbtxfallbacklanguage{english}
\usepackage{cite}%
\usepackage[T1]{fontenc}
\usepackage{hyperref}
\usepackage{lmodern}
\usepackage{mathdots}
\usepackage{mathscinet}
\usepackage{mathtools}
 \mathtoolsset{mathic=true}
\usepackage{microtype}
\usepackage{url}

\usepackage{176arxiv}

\title{Matricial Canonical Moments and Parametrization of Matricial Hausdorff Moment Sequences}
\author{Bernd Fritzsche \and Bernd Kirstein \and Conrad M\"adler}

\begin{document}
\maketitle

\begin{abstract}
 In this paper we study moment sequences of matrix-valued measures on compact intervals.
 A complete parametrization of such sequences is obtained via a symmetric version of matricial canonical moments.
 Furthermore, distinguished extensions of finite moment sequences are characterized in this framework.
 The results are applied to the underlying matrix-valued measures, generalizing some results from the scalar theory of canonical moments.
\end{abstract}

\begin{description}
 \item[Keywords:] Moments of matrix-valued measures, matricial Hausdorff moment problem, coupled quadruples of \tnnH{} block Hankel matrices, matricial canonical moments.
\end{description}

\section{Introduction}
 In this paper, we continue our program of studying the inner structure of moment sequences of matrix-valued measures on subsets of the real line.
 Now we deal with matrix-valued measures on compact intervals.
 Guided by our former investigations on the whole real axis (see~\zitas{MR2570113,MR2805417,MR3014199}) on the one side and of closed semi-bounded infinite intervals (see~\zitas{MR2735313,MR3014201,MR3133464,MR3611479}) on the other side, our approach is mainly of matrix theoretical nature.
 It is divided into two principal steps.
 In the first step the focus lies on the associated extension problem for \tnnH{} (not necessarily invertible) \tbHms{}.
 In the second step we look for an appropriate parametrization of the corresponding matricial moment sequences which reflects the main features of the inner dependencies between the matrices forming the moment sequence under consideration.
 What concerns the first step in the case of the matricial Hausdorff moment sequences a complete treatment could be realized in our paper~\zita{arXiv:1701.04246}.
 The main goal of this paper is to handle the second step.
 More precisely, our main purpose is to construct a parametrization of the set of all matricial Hausdorff moment sequences of a given length, \tie{}, the moment space on \(\ab\) of given order.

 Before explaining our strategy in more detail we give a short sketch about the history of this topic.
 In the scalar case this theme was handled by associating the canonical moments \(p_k\).
 For this subject, we refer to the monograph Dette/Studden~\zita{MR1468473}, who also give there a brief outline of the corresponding history, outgoing from considerations of moment spaces, mentioning Kre\u{\i}n~\zita{MR0044591}, Karlin/Studden~\zitas{MR0204922,MR0196871}, Kre\u{\i}n/Nudel\cprime man~\zita{MR0458081}, and Karlin/Shapley~\zita{MR0059329}, to the first description in Seall/Wetzel~\zita{MR0124654}, and its independent rediscovery in Skibinsky~\zita{MR0228040}.
 For an initial expository introduction to moment space geometry and canonical moment sequences for distributions on the unit interval we refer also to Skibinsky~\zita{MR861911}.
 
 Focusing on characterizing the domain of the measure of orthogonality in the framework of Favard's theorem for matrix-valued measures and orthogonal matrix polynomials, Dette and Studden introduced in~\zita{MR1883272} matrix canonical moments \(U_k\) and generalized several scalar identities involving classical canonical moments to the matrix case.
 The unsymmetric construction used in~\zita{MR1883272} is suitable to obtain right matrix product representations, compatible with the right inner product structure of the right module of matrix polynomials considered there.
 Due to this unbalance, the matrix canonical moments \(U_k\) seem not to be \tH{} in general.
 However, their spectra are located in the interval \([0,1]\).
 To obtain a one-to-one correspondence between matrix-valued measures on \(\ab\), or equivalently their complete sequences of power moments, and a class of sequences of matrices, belonging indeed to a corresponding matricial interval \(\matint{\Oqq}{\Iq}\) with respect to the L\"owner semi-ordering for \tH{} matrices, we choose a slightly different symmetric construction.
 
 This paper is organized as follows.
 In \rsec{S.MP}, we introduce two related types of matricial power moment problems associated with a fixed non-empty Borelian subset of the real axis.
 Since our considerations make intensive use of former results on moment sequences of matrix-valued measures on the real axis \(\R\) and the intervals \(\rhl\) and \(\lhl\), we recall these facts in \rsecss{H.S}{K.S}.
 In \rsec{F.S} basic observations concerning the matricial moment problem on \(\ab\) are stated.
 Furthermore, we recall the answer to the one-step extension problem for truncated matricial Hausdorff sequences (see~\rthm{165.T112}).
 Towards our main result, in \rsec{F.S.fp} we associate to a Hausdorff moment sequence \(\seqska\) a sequence \(\fpseqka\), reducing thereby in a first step the inner dependencies of \(\seqska\).
 To uncover the remaining dependencies between the matrices \(\fpseqka\), we proceed in \rdefn{F.D.cia} to a sequence \(\seqciaka\) of \tnnH{} matrices, where \(\cia{j}\) is located in a matricial interval, the upper bound \(\Pd{j-1}\) of which corresponds to the orthogonal projection onto the column space of the width \(\dia{j-1}\) of the preceding section of the moment space.
 In this way, we have constructed the desired parametrization in our main result \rthm{F.T.Fggcia}.
 Several special configurations in the sequence \(\seqciaka\) are considered in \rprop{F.C.ecase}.
 Some of them turn out to be connected to the notions of centrality and complete degeneracy studied already in our former paper~\zita{arXiv:1701.04246}.
 \rsec{F.BT} is devoted to demonstrate the invariance of the parameter sequence \(\seqciaka\) under affine transformations of the measure having the moment sequence \(\seqska\).
 In the final \rsec{F.CM}, we translate our results on infinite Hausdorff moment sequences \(\seqsinf\) into the language of \tnnH{} measures on \(\ab\).
 This concerns, \teg{}, the characterization of symmetry properties of measures.
 In \rsec{F.s2.cm}, some concrete connections to formulas in Dette/Studden~\zita{MR1883272} are described and using Moore--Penrose inverses some of their results are proved to be valid for arbitrary points in the moment space without regularity assumptions.

\section{On matricial power moment problems on Borel subsets of the real axis}\label{S.MP}
 In this section we are going to formulate two related general classes of matricial power moment problems.
 Before doing this we have to introduce some terminology.
 We denote by \(\Z\) the set of all integers.
 Let \(\N\defeq\setaa{n\in\Z}{n\geq1}\).
 Furthermore, we write \(\R\) for the set of all real numbers and \(\C\) for the set of all complex numbers.
 In the whole paper \(p\) and \(q\) are arbitrarily fixed integers from \(\N\).
 We write \(\Cpq\) for the set of all complex \tpqa{matrices} and \(\Cp\) is short for \(\Coo{p}{1}\).
 When using \(m,n,r,s,\dotsc\) instead of \(p,q\) in this context, we always assume that these are integers from \(\N\).
 We write \(A^\ad\) for the conjugate transpose of a complex \tpqa{matrix} \(A\).
 Denote by \(\Cggq\defeq\setaa{M\in\Cqq}{v^\ad Mv\in[0,\infp)\text{ for all }v\in\Cq}\) the set of \tnnH{} matrices from \(\Cqq\).

 Let \((\mathcal{X},\mathfrak{X})\) be a measurable space.
 Each countably additive mapping whose domain is \(\mathfrak{X}\) and whose values belong to \(\Cggq\) is called a \tnnH{} \tqqa{measure} on \((\mathcal{X}, \mathfrak{X})\).
 For the integration theory with respect to \tnnH{} measures, we refer to Kats~\zita{MR0080280} and Rosenberg~\zita{MR0163346}.
 
 Let \(\BsAR\) (\tresp{}\ \(\BsAC \)) be the \(\sigma\)\nobreakdash-algebra of all Borel subsets of \(\R\) (\tresp{}\ \(\C \)).
 In the whole paper, \(\Omega\) stands for a non-empty set belonging to \(\BsAR\).
 Let \(\BsAO \) be the \(\sigma\)\nobreakdash-algebra of all Borel subsets of \(\Omega\) and let \(\MggqO \) be the set of all \tnnH{} \tqqa{measures} on \((\Omega,\BsAO )\).
 Observe that \(\Mggoa{1}{\Omega}\) coincides with the set of ordinary measures on \((\Omega,\BsAu{\Omega})\) with values in \([0,\infp)\).

 Let \(\NO\defeq\setaa{m\in\Z}{m\geq0}\).
 In the whole paper, \(\kappa\) is either an integer from \(\NO\) or \(\infi\).
 In the latter case, we have \(2\kappa=\infi\) and \(2\kappa+1=\infi\).
 Given \(\upsilon,\omega\in\R\cup\set{-\infty,\infp}\), we denote by \(\mn{\upsilon}{\omega}\defeq\setaa{k\in\Z}{\upsilon\leq k\leq\omega}\) the perhaps empty or unbounded corresponding section of \(\Z\).
 Let \(\Mggoua{q}{\kappa}{\Omega}\) be the set of all \(\mu\in\MggqO \) such that for each \(j\in\mn{0}{\kappa}\) the power function \(x\mapsto x^j\) defined on \(\Omega\) is integrable with respect to \(\mu\) .
 If \(\mu\in\MggquO{\kappa}\), then, for all \(j\in\mn{0}{\kappa}\), the matrix
\beql{I.G.mom}
 \mpm{\mu}{j}
 \defeq\int_\Omega x^j\mu\rk{\dif x}
\eeq
 is called (power) \emph{moment of \(\mu\)} of order \(j\).
 Obviously, we have \(\Mggoa{q}{\Omega}=\Mggoua{q}{0}{\Omega}\subseteq \MggquO{\ell}\subseteq \MggquO{\ell+1}\subseteq\MggquO{\infi}\) for every choice of \(\ell\in\NO \) and \(\mpm{\mu}{0}=\mu\rk{\Omega}\) for all \(\mu\in\Mggoa{q}{\Omega}\).
 Furthermore, if \(\Omega\) is bounded, then one can easily see that \(\MggqO =\MggquO{\infi}\).
 In particular, for \(\ug,\obg\in\R\) with \(\ug<\obg\), we have \(\MggqF=\MggquF{\infi}\).
 We now state the general form of the moment problem lying in the background of our considerations:

\begin{Problem}[\mprob{\Omega}{\kappa}{=}]
 Given a sequence \(\seqska \) of complex \tqqa{matrices}, parametrize the set \(\MggqOsg{\kappa}\) of all \(\sigma\in\Mggoua{q}{\kappa}{\Omega}\) satisfying \(\mpm{\sigma}{j}=\su{j}\) for all \(j\in\mn{0}{\kappa}\).
\end{Problem}

 The set \(\CHq\defeq\setaa{M\in\Cqq}{M^\ad=M}\) of \tH{} matrices from \(\Cqq\) is a partially ordered vector space over the field \(\R\) with positive cone \(\Cggq\).
 For two complex \tqqa{matrices} \(A\) and \(B\), we write \(A\lleq B\) or \(B\lgeq A\) if \(A,B\in\CHq\) and \(B-A\in\Cggq\) are fulfilled.
 Let \(\Opq\) be the zero matrix from \(\Cpq\).
 Sometimes, if the size of the zero matrix is clear from the context, we will omit the indices and write \(\NM\).
 For a complex \tqqa{matrix} \(A\), we have obviously \(A\lgeq\NM\) if and only if \(A\in\Cggq\).
 The above mentioned partial order \(\lleq\) on the set of \tH{} matrices is sometimes called \emph{L\"owner semi-ordering}.
 For \(\kappa<\infi\) the following modification of the above mentioned problem is also of interest:

\begin{Problem}[\mprob{\Omega}{m}{\lleq}]
 Given \(m\in\NO\) and a sequence \(\seqs{m}\) of complex \tqqa{matrices}, parametrize the set \(\MggqOskg{m}\) of all \(\sigma\in\Mggoua{q}{m}{\Omega}\) satisfying \(\mpm{\sigma}{m}\lleq\su{m}\) and \(\mpm{\sigma}{j}=\su{j}\) for all \(j\in\mn{0}{m-1}\).
\end{Problem}

 In this paper, the main interest is directed to the case that \(\Omega\) is a compact interval \(\ab\) of the real axis \(\R\).
 It is obvious that each solution of Problem~\mprob{\ab}{m}{=} generates in a natural way solutions to each of the problems \mprob{\rhl}{m}{=}, \mprob{\lhl}{m}{=}, and \mprob{\R}{m}{=}.
 The last mentioned three matricial moment problems were intensively studied in our former work (see~\zitas{MR2735313,MR2570113,MR2805417,MR3014201,MR3133464,MR3611479,MR3014199}).
 In particular, we analyzed the inner structure of matricial moment sequences associated with each of the sets \(\rhl\), \(\lhl\), and \(\R\).
 These investigations led us to corresponding parametrizations which provide essential insights into the inner structure of the matricial moment sequences under consideration.
 Our techniques were mainly based on the use of Schur complements in appropriately built \tnnH{} \tbHms{}.
 The matricial Hausdorff moment sequences which are the central object of this paper have a much more difficult structure which is caused by the interplay of four coupled matricial Hamburger moment sequences.
 The essential point of our strategy will be to find an appropriate way of organizing this interplay.
 In order to prepare the cornerstones of our approach we have to recall some material which is connected with matricial power moment problems associated with the case that \(\Omega\) is one of the sets \(\R\), \(\rhl\) or \(\lhl\).
 This will be done in the next two sections.

\section{Matricial Hamburger moment sequences and \hhp{s}}\label{H.S}
 We recall classes of sequences of complex \tqqa{matrices} corresponding to solvability criteria for the matricial moment problem on \(\Omega=\R\).
 
 For each \(n\in\NO\), denote by \(\Hggqu{2n}\) the set of all sequences \(\seqs{2n}\) of complex \tqqa{matrices}, for which the corresponding \tbHm{}
\(
  \Hu{n}
  \defeq\matauuo{\su{j+k}}{j,k}{0}{n}
\)
 is \tnnH{}.
 Furthermore, denote by \(\Hggqinf\) the set of all sequences \(\seqsinf\) of complex \tqqa{matrices} satisfying \(\seqs{2n}\in\Hggqu{2n}\) for all \(n\in\NO\).
 The sequences belonging to \(\Hggqu{2n}\) or \(\Hggqinf\) are said to be \emph{\tHnnd}.
 In~\zita{MR2570113} and subsequent papers of the corresponding authors, the sequences belonging to \(\Hggqu{2\kappa}\) were called \emph{Hankel non-negative definite}.
 Our terminology here differs for the sake of consistency.

\bpropnl{\tcf{}~\zitaa{MR1624548}{\cthm{3.2}{213}}}{115.T1615}
 Let \(n\in\NO\) and let \(\seqs{2n}\) be a sequence of complex \tqqa{matrices}.
 Then \(\MggquRskg{2n}\neq\emptyset\) if and only if \(\seqs{2n}\in\Hggqu{2n}\).
\eprop

 The characterization of the solvability of the ``\(=\)''\nobreakdash-version of the moment problem associated with \(\Omega=\R\) requires additional tools which will be prepared now. 
 Let \(n\in\NO\).
 Denote by \(\Hggequ{2n}\) the set of all sequences \(\seqs{2n}\) of complex \tqqa{matrices}, for which there exists a pair \((\su{2n+1},\su{2n+2})\) of complex \tqqa{matrices}, such that the sequence \(\seqs{2n+2}\) belongs to \(\Hggqu{2n+2}\).
 Denote by \(\Hggequ{2n+1}\) the set of all sequences \(\seqs{2n+1}\) of complex \tqqa{matrices}, for which there exists a complex \tqqa{matrix} \(\su{2n+2}\), such that the sequence \(\seqs{2n+2}\) belongs to \(\Hggqu{2n+2}\).
 Furthermore, let \(\Hggeqinf\defeq\Hggqinf\).
 The sequences belonging to \(\Hggequ{2n}\), \(\Hggequ{2n+1}\), or \(\Hggeqinf\) are said to be \emph{\tHnnde}.

\bpropnl{\zitaa{MR2805417}{\cthm{6.6}{486}} (\tcf{}~\zitaa{MR1624548}{\cthm{3.1}{211}})}{H.P.MPsolv}
 Let \(\seqska \) be a sequence of complex \tqqa{matrices}.
 Then \(\MggqkappaRsg\neq\emptyset\) if and only if \(\seqska\in\Hggeqkappa\).
\eprop

 In order to describe the inner dependencies of \tHnnd{} sequences (even in the general case of not necessarily invertible \tnnH{} \tbHms{} \(\Hu{n}\)) in~\zita{MR2570113} the so-called canonical Hankel parametrization was introduced and further discussed in~\zitas{MR2805417,MR3014199}.
 We slightly reformulate this notion in a more convenient form in \rdefn{102.HPN}.
 Therefore, we need several (block) matrices.
 For later use we mention some simple observations concerning the arithmetics of these objects:

\bnotal{N.HKG}
 Let \(\seqska \) be a sequence of complex \tpqa{matrices}.
 Then let the \tbHms{} \(\Hu{n}\), \(\Ku{n}\), and \(\Gu{n}\) be given by
 \(
  \Hu{n}
  \defg\matauuuo{s_{j+k}}{j}{k}{0}{n}
  \)
 for all \(n\in\NO\) with \(2n\leq\kappa\), by
 \(
  \Ku{n}
  \defg\matauuuo{s_{j+k+1}}{j}{k}{0}{n}
 \)
 for all \(n\in\NO\) with \(2n+1\leq\kappa\), and by
 \(
  \Gu{n}
  \defg\matauuuo{s_{j+k+2}}{j}{k}{0}{n}
 \)
 for all \(n\in\NO\) with \(2n+2\leq\kappa\), \tresp{}
\enota

 Let \(\Iq\defeq\matauuo{\Kronu{jk}}{j,k}{1}{q}\) be the identity matrix from \(\Cqq\), where \(\Kronu{jk}\) is the Kronecker delta.
 Sometimes, we will omit the indices and write \(\EM\).
 Given \(n\in\N\) arbitrary rectangular complex matrices \(A_1,A_2,\dotsc,A_n\), we write \(\diag\seq{A_j}{j}{1}{n}=\diag\rk{A_1,A_2,\dotsc,A_n}\defg\matauuuo{\Kronu{jk}A_j}{j}{k}{1}{n}\) for the corresponding block diagonal matrix.
 
\breml{H.R.l*H}
 Let \(\lambda\in\C\) and let \(\seqska\) be a sequence of complex \tpqa{matrices}.
 Let the sequences \(\seq{u_j}{j}{0}{\kappa}\) and \(\seq{v_j}{j}{0}{\kappa}\) be given by \(u_j\defeq\lambda\su{j}\) and \(v_j\defeq \lambda^j\su{j}\), \tresp{}
 Then \(\Huo{n}{u}=\lambda\Huo{n}{s}\) and \(\Huo{n}{v}=\ek{\diag\seq{\lambda^j\Iq}{j}{0}{n}}\Huo{n}{s}\ek{\diag\seq{\lambda^j\Ip}{j}{0}{n}}\) for all \(n\in\NO\) with \(2n\leq\kappa\).
\erem

 We write \(\sdiag{A}{m}\defeq\diag\seq{A}{j}{0}{m}\), given \(A\in\Cpq\) and \(m\in\NO\).
 
\breml{H.R.L*H*R}
 Let \(L\in\Coo{\ell}{p}\), let \(R\in\Coo{q}{r}\), and let \(\seqska\) be a sequence of complex \tpqa{matrices}.
 Let the sequence \(\seq{w_j}{j}{0}{\kappa}\) be given by \(w_j\defeq L\su{j}R\).
 For all \(n\in\NO\) with \(2n\leq\kappa\), then \(\Huo{n}{w}=\sdiag{L}{n+1}\Huo{n}{s}\sdiag{R}{n+1}\).
\erem

\breml{H.R.H+H}
 Let \(\seqska \) and \(\seqt{\kappa}\) be sequences of complex \tpqa{matrices} and let the sequence \(\seq{x_j}{j}{0}{\kappa}\) be given by \(x_j\defeq\su{j}+t_j\).
 For all \(n\in\NO\) with \(2n\leq\kappa\), then \(\Huo{n}{x}=\Huo{n}{s}+\Huo{n}{t}\).
\erem

 Given \(n\in\N\) arbitrary rectangular complex matrices \(A_1,A_2,\dotsc,A_n\), we write \(\col\seq{A_j}{j}{1}{n}=\col\rk{A_1,A_2,\dotsc,A_n}\) (\tresp{}, \(\row\seq{A_j}{j}{1}{n}\defeq\mat{A_1,A_2,\dotsc,A_n}\)) for the block column (\tresp{}, block row) build from the matrices \(A_1,A_2,\dotsc,A_n\) if their numbers of columns (\tresp{}, rows) are all equal.

\bnotal{N.yz}
 Let \(\seqska \) be a sequence of complex \tpqa{matrices}.
 Then let the block row \(\yuu{\ell}{m}\) and the block column \(\zuu{\ell}{m}\) be given by \(\yuu{\ell}{m}\defg\col\seq{s_j}{j}{\ell}{m}\) and \(\zuu{\ell}{m}\defg\row\seq{s_j}{j}{\ell}{m}\)
 for all \(\ell,m\in\NO\) with \(\ell\leq m\leq\kappa\), \tresp{}
\enota

 In addition to \rremss{H.R.l*H}{H.R.L*H*R}, we have:

\breml{H.R.l*yz}
 Let \(\lambda\in\C\) and let \(\seqska\) be a sequence of complex \tpqa{matrices}.
 Let the sequences \(\seq{u_j}{j}{0}{\kappa}\) and \(\seq{v_j}{j}{0}{\kappa}\) be given by \(u_j\defeq\lambda\su{j}\) and \(v_j\defeq \lambda^j\su{j}\) \tresp{}
 Then \(\yuuo{m}{n}{u}=\lambda\yuuo{m}{n}{s}\) and \(\zuuo{m}{n}{u}=\lambda\zuuo{m}{n}{s}\) for all \(m,n\in\mn{0}{\kappa}\).
 Furthermore, \(\yuuo{\ell}{m}{v}=\ek{\diag\seq{\lambda^j\Iq}{j}{\ell}{m}}\yuuo{\ell}{m}{s}\) and \(\zuuo{\ell}{m}{v}=\zuuo{\ell}{m}{s}\ek{\diag\seq{\lambda^j\Ip}{j}{\ell}{m}}\) for all \(\ell,m\in\NO\) with \(\ell\leq m\leq\kappa\).
\erem

\breml{H.R.L*yz*R}
 Let \(L\in\Coo{\ell}{p}\), let \(R\in\Coo{q}{r}\), and let \(\seqska\) be a sequence of complex \tpqa{matrices}.
 Let the sequence \(\seq{w_j}{j}{0}{\kappa}\) be given by \(w_j\defeq L\su{j}R\).
 For all \(m,n\in\mn{0}{\kappa}\) with \(m\leq n\), then \(\yuuo{m}{n}{w}=\sdiag{L}{n-m+1}\yuuo{m}{n}{s}R\) and \(\zuuo{m}{n}{w}=L\zuuo{m}{n}{s}\sdiag{R}{n-m+1}\).
\erem

 The \tbHm{} \(\Hu{n}\) admits the following \tbr{}:

\breml{H.R.Hblock}
 Let \(\seqska\) be a sequence of complex \tpqa{matrices}.
 For all \(n\in\N\) with \(2n\leq\kappa\), then \(\Hu{n}=\smat{\Hu{n-1} & \yuu{n}{2n-1}\\ \zuu{n}{2n-1} & s_{2n}}\).
\erem

 The considerations in this paper heavily rely on the use of the following generalized inverse for complex matrices:
 For each matrix \(A\in\Cpq\), there exists a uniquely determined matrix \(X\in\Cqp\), satisfying the four equations
\begin{align}\label{mpi}
 AXA&=A,&
 XAX&=X,&
 \rk{AX}^\ad&=AX,&
&\text{and}&
 \rk{XA}^\ad&=XA.
\end{align}
 This matrix \(X\) is called the \emph{Moore--Penrose inverse of \(A\)} and is denoted by \(A^\mpi\).
 Concerning a detailed treatment of the machinery of Moore--Penrose inverses we refer to~\zita{MR0338013},~\zitaa{MR1105324}{\cch{1}}, and~\zitaa{MR1987382}{\cch{1}}.
 For our purposes it is convenient to apply~\zitaa{MR1152328}{\csec{1.1}}.
 In \rapp{M.S} we summarize a collection of properties of Moore--Penrose inverses used in this paper.

\bnotal{N.Lambda}
 Let \(\seqska \) be a sequence of complex \tpqa{matrices}.
 Then let \(\Tripu{0}\defg\Opq\), \(\Sigmau{0}\defg\Opq\), \(\Tripu{n}\defg\zuu{n}{2n-1}\Hu{n-1}^\mpi\yuu{n}{2n-1}\), and \(  \Sigmau{n}\defg\zuu{n}{2n-1}\Hu{n-1}^\mpi\Ku{n-1}\Hu{n-1}^\mpi\yuu{n}{2n-1}\) for all \(n\in\N\) with \(2n-1\leq\kappa\).
 Furthermore, let \(\Mu{0}\defg\Opq\), \(\Nu{0}\defg\Opq\), \(\Mu{n}\defg\zuu{n}{2n-1}\Hu{n-1}^\mpi\yuu{n+1}{2n}\), and \(\Nu{n}\defg\zuu{n+1}{2n}\Hu{n-1}^\mpi\yuu{n}{2n-1}\) for all \(n\in\N\) with \(2n\leq\kappa\).
 Moreover, let \(\Lambdau{n}\defg\Mu{n}+\Nu{n}-\Sigmau{n}\) for all \(n\in\NO\) with \(2n\leq\kappa\).
\enota

 In view of \rrem{A.R.l*A} and \eqref{mpi}, we conclude from \rremss{H.R.l*H}{H.R.l*yz}:

\breml{H.R.l*MN}
 Let \(\lambda\in\C\) and let \(\seqska\) be a sequence of complex \tpqa{matrices}.
 Let the sequence \(\seq{u_j}{j}{0}{\kappa}\) be given by \(u_j\defeq\lambda\su{j}\).
 Then \(\Tripuo{n}{u}=\lambda\Tripuo{n}{s}\) and \(\Sigmauo{n}{u}=\lambda\Sigmauo{n}{s}\) for all \(n\in\NO\) with \(2n-1\leq\kappa\).
 Furthermore, \(\Muo{n}{u}=\lambda\Muo{n}{s}\) and \(\Nuo{n}{u}=\lambda\Nuo{n}{s}\) and, consequently, \(\Lambdauo{n}{u}=\lambda\Lambdauo{n}{s}\) for all \(n\in\NO\) with \(2n\leq\kappa\).
\erem

 Using \rrem{A.R.UA+V}, we obtain from \rremss{H.R.L*H*R}{H.R.L*yz*R} moreover:

\breml{H.R.U*MN*V}
 Let \(U\in\Coo{u}{p}\) with \(U^\ad U=\Ip\), let \(V\in\Coo{q}{v}\) with \(VV^\ad=\Iq\), and let \(\seqska\) be a sequence of complex \tpqa{matrices}.
 Let the sequence \(\seq{w_j}{j}{0}{\kappa}\) be given by \(w_j\defeq U\su{j}V\).
 Then \(\Tripuo{n}{w}=U\Tripuo{n}{s}V\) and \(\Sigmauo{n}{w}=U\Sigmauo{n}{s}V\) for all \(n\in\NO\) with \(2n-1\leq\kappa\).
 Furthermore, \(\Muo{n}{w}=U\Muo{n}{s}V\) and \(\Nuo{n}{w}=U\Nuo{n}{s}V\) and, consequently, \(\Lambdauo{n}{w}=U\Lambdauo{n}{s}V\) for all \(n\in\NO\) with \(2n\leq\kappa\).
\erem

 If \(A\) is a square matrix, then denote by \(\det A\) the determinant of \(A\).

\bleml{H.L.z^j*MN}
 Let \(\xi\in\C\) with \(\abs{\xi}=1\) and let \(\seqska \) be a sequence of complex \tpqa{matrices}.
 Let the sequence \(\seq{v_j}{j}{0}{\kappa}\) be given by \(v_j\defeq\xi^j\su{j}\).
 Then \(\Tripuo{n}{v}=\xi^{2n}\Tripuo{n}{s}\) and \(\Sigmauo{n}{v}=\xi^{2n+1}\Sigmauo{n}{s}\) for all \(n\in\NO\) with \(2n-1\leq\kappa\).
 Furthermore, \(\Muo{n}{v}=\xi^{2n+1}\Muo{n}{s}\) and \(\Nuo{n}{v}=\xi^{2n+1}\Nuo{n}{s}\) and, consequently, \(\Lambdauo{n}{v}=\xi^{2n+1}\Lambdauo{n}{s}\) for all \(n\in\NO\) with \(2n\leq\kappa\).
\elem
\bproof
 For all \(r\in\N\), all \(\ell,m\in\NO\) with \(\ell\leq m\), and all \(\omega\in\C\setminus\set{0}\), the matrix \(D_{r,\ell,m}(\omega)\defeq\diag\seq{\omega^j\Iu{r}}{j}{\ell}{m}\) fulfills \(D_{r,\ell,m}(\omega)=\omega^\ell D_{r,0,m-\ell}(\omega)\) and \(\ek{D_{r,\ell,m}(\omega)}^\ad=D_{r,\ell,m}(\ko \omega)\) and, in view of \(\omega\neq0\), furthermore \(\det D_{r,\ell,m}(\omega)\neq0\) and \(\ek{D_{r,\ell,m}(\omega)}^\inv=D_{r,\ell,m}(\omega^\inv)\).
 Since \(\abs{\xi}=1\), for all \(r\in\N\) and all \(n\in\NO\), the matrix \(U_{r,n}\defeq D_{r,0,n}(\xi)=\diag\seq{\xi^j\Iu{r}}{j}{0}{n}\) fulfills then \(\det U_{r,n}\neq0\) and \(U_{r,n}^\inv=U_{r,n}^\ad\).
 Using these notations, we have \(\yuuo{\ell}{m}{v}=\xi^\ell U_{q,m-\ell}\yuuo{\ell}{m}{s}\) and \(\zuuo{\ell}{m}{v}=\xi^\ell\zuuo{\ell}{m}{s}U_{p,m-\ell}\) for all \(\ell,m\in\NO\) with \(\ell\leq m\leq\kappa\), according to \rrem{H.R.l*yz}, and furthermore \(\Huo{m}{v}=U_{q,m}\Huo{m}{s}U_{p,m}\) for all \(m\in\NO\) with \(2m\leq\kappa\) and \(\Kuo{m}{v}=\xi U_{q,m}\Kuo{m}{s}U_{p,m}\) for all \(m\in\NO\) with \(2m+1\leq\kappa\), according to \rrem{H.R.l*H}.
 By virtue of \rrem{A.R.UA+V}, then \(\rk{\Huo{m}{v}}^\mpi=U_{p,m}^\inv\rk{\Huo{m}{s}}^\mpi U_{q,m}^\inv\) for all \(m\in\NO\) with \(2m\leq\kappa\) follows.
 In view of \rnota{N.Lambda}, from these formulas we easily conclude the asserted identities.
\eproof

 Now we introduce the central object of this section.

\bdefnl{102.HPN}
 Let \(\seqska \) be a sequence of complex \tpqa{matrices}.
 Let the sequence \(\seq{\hpu{j}}{j}{0}{\kappa}\) be given by
 \(
  \hpu{2k}
  \defg\su{2k}-\Tripu{k}
 \)
 for all \(k\in\NO\) with \(2k\leq\kappa\) and by
 \(
  \hpu{2k+1}
  \defg\su{2k+1}-\Lambdau{k}
 \)
 for all \(k\in\NO\) with \(2k+1\leq\kappa\).
 Then we call \(\seq{\hpu{j}}{j}{0}{\kappa}\) the \emph{\thpfa{\(\seqska \)}}.
\edefn

\breml{102.MV}
 Let \(\seqska \) be a sequence of complex \tpqa{matrices} with \thpf{} \(\hpseqo{\kappa}\).
 For each \(k\in\mn{0}{\kappa}\), then the matrix \(\hpu{k}\) is built from the matrices \(\su{0},\su{1},\dotsc,\su{k}\).
 In particular, for each \(m\in\mn{0}{\kappa}\), the \thpfa{\(\seqs{m}\)} coincides with \(\hpseqo{m}\).
\erem

 For all \(\iota\in\NOinf\) and each non-empty set \(\mathcal{X}\), let
\beql{seqset}
 \seqset{\iota}{\mathcal{X}}
 \defeq\setaca*{\seq{X_j}{j}{0}{\iota}}{X_j\in\mathcal{X}\text{ for all }j\in\mn{0}{\iota}}.
\eeq
 It is readily checked that a sequence \(\seqska\in\seqset{\kappa}{\Cpq}\) can be recursively reconstructed from its \thpf{} \(\hpseqo{\kappa}\) via \(\su{2k}=\Tripu{k}+\hpu{2k}\) for all \(k\in\NO\) with \(2k\leq\kappa\) and \(\su{2k+1}=\Lambdau{k}+\hpu{2k+1}\) for all \(k\in\NO\) with \(2k+1\leq\kappa\).
 Furthermore, given an arbitrary sequence \(\hpseqo{\kappa}\in\seqset{\kappa}{\Cpq}\), we can build recursively, in the above mentioned way, a sequence \(\seqska\in\seqset{\kappa}{\Cpq}\) such hat \(\hpseqo{\kappa}\) is the \thpfa{\(\seqska\)}.
 Consequently, we have:

\breml{R1429}
 The mapping \(\Phi\colon\seqset{\kappa}{\Cpq}\to\seqset{\kappa}{\Cpq}\) defined by \(\seqska\mapsto\hpseqo{\kappa}\) is bijective.
\erem

 We write \(\ran{A}\defeq\setaa{Ax}{x\in\Cq}\) for the column space of a complex \tpqa{matrix} \(A\).
 As already mentioned, the \thp{s} were introduced in a slightly different form as a pair of two sequences in~\zita{MR2570113} under the name \emph{canonical Hankel parametrization}.
 Several properties of a sequence connected with \tHnnd{ness} can be characterized in terms of its \thp{s}, the main result being the following:
 
\bpropnl{\tcf{}~\zitaa{MR2805417}{\cprop{2.10(b)}{454,} and \cprop{2.15(b)}{457}}}{H.P.HggHP}
 A sequence \(\seqs{2\kappa}\) of complex \tqqa{matrices} belongs to \(\Hggqu{2\kappa}\) if and only if
 \bAeqi{0}
  \item \(\hpu{2k+1}\in\CHq\) and \(\ran{\hpu{2k+1}}\subseteq \ran{\hpu{2k}}\) for all \(k\in\mn{0}{\kappa-1}\)
 \eAeqi
 and furthermore
 \bAeqi{1}
  \item\(\hpu{2k}\in\Cggq\) for all \(k\in\mn{0}{\kappa}\) and \(\ran{\hpu{2k+2}}\subseteq \ran{\hpu{2k}}\) for all \(k\in\mn{0}{\kappa-2}\).
 \eAeqi
\eprop

 The \thp{s} also occur in connection with three term recurrence relations for systems of orthogonal matrix polynomials with respect to \tnnH{} measures on \(\rk{\R,\BsAR}\) (see~\zitaa{MR2570113}{\csec{3}} and~\zita{MR2805417}).
 In the remaining part of this section, we summarize some later used results on the arithmetics of \thp{s}.
 From \rremss{H.R.l*MN}{H.R.U*MN*V} we obtain:

\breml{H.R.l*h}
 Let \(\lambda\in\C\) and let \(\seqska\) be a sequence of complex \tpqa{matrices}.
 Then \(\seq{\lambda\hpu{j}}{j}{0}{\kappa}\) coincides with the \thpfa{\(\seq{\lambda s_j}{j}{0}{\kappa}\)}.
\erem

\breml{H.R.U.h.V}
 Let \(U\in\Coo{u}{p}\) with \(U^\ad U=\Ip\), let \(V\in\Coo{q}{v}\) with \(VV^\ad=\Iq\), and let \(\seqska\) be a sequence of complex \tpqa{matrices}.
 Then \(\seq{U\hpu{j}V}{j}{0}{\kappa}\) coincides with the \thpfa{\(\seq{Us_jV}{j}{0}{\kappa}\)}.
\erem

 In view of \rexam{H.E.d^j*s} below, the following observation, which can be easily seen from \rlem{H.L.z^j*MN}, is of interest:
 
\breml{H.R.d^j*h}
 Let \(\xi\in\C\) with \(\abs{\xi}=1\) and let \(\seqska\) be a sequence of complex \tpqa{matrices}.
 Then \(\seq{\xi^j\hpu{j}}{j}{0}{\kappa}\) coincides with the \thpfa{\(\seq{\xi^j s_j}{j}{0}{\kappa}\)}.
\erem

 Now we turn our attention to the Schur complement associated with block matrices.
 If \(M=\tmat{A & B\\ C & D}\)
 is the \tbr{} of a complex \taaa{(p+q)}{(r+s)}{matrix} \(M\) with \taaa{p}{r}{block} \(A\), then the matrix
\[%
 M\schca A
 \defg D-CA^\mpi B
\]
 is called the \emph{Schur complement of \(A\) in \(M\)}.
 Concerning this concept, we refer to~\zitas{MR2160825}.
 In view of \eqref{mpi}, we get from \rrem{A.R.l*A}:

\breml{ab.R1111}
 If \(\eta\in\C\), then
 \(
  \eta M
  =
  \tmat{
   \eta A&\eta B\\
   \eta C&\eta D
  }
 \)
 and \((\eta M)\schca (\eta A)=\eta(M\schca A)\).
\erem

\blemnl{\tcf{}~\zitaa{MR1152328}{\clemss{1.1.9}{18}{1.1.10}{19}}}{L.AEP}
 Let \(M=\tmat{A&B\\C&D}\) be the \tbr{} of a complex \taaa{(p+q)}{(p+q)}{matrix} \(M\) with \tppa{block} \(A\).
 Then the following statements are equivalent:
 \baeqi{0}
  \item \(M\) is \tnnH{}.
  \item \(A\) and \(M\schca A\defg D-CA^\mpi B\) are both \tnnH{} and furthermore \(\ran{B}\subseteq\ran{A}\) and \(C=B^\ad\).
  \item \(D\) and \(M\schca D\defg A-BD^\mpi C\) are both \tnnH{} and furthermore \(\ran{C}\subseteq\ran{D}\) and \(B=C^\ad\).
 \eaeqi
\elem

\bpropnl{\tcf{}~\zitaa{MR2160825}{\cthm{5.9(ii)}{148}}}{A.R.sc+sc}
 Let \(M_1=\tmat{A_1& B_1\\ C_1 & D_1}\) and \(M_2=\tmat{A_2& B_2\\ C_2 & D_2}\)
 be the \tbr{s} of two \tnnH{} \taaa{(p+q)}{(p+q)}{matrices} \(M_1\) and \(M_2\) with \tppa{block} \(A_1\) and \(A_2\), \tresp{}
 Then \(M\defg M_1+M_2\) admits the \tbr{} \(M=\tmat{A_1+A_2&B_1+B_2\\C_1+C_2 & D_1+D_2}\)
 and \(M\schca (A_1+A_2)\lgeq\rk{M_1\schca A_1}+\rk{M_2\schca A_2}\).
\eprop

 By virtue of \rrem{H.R.Hblock}, we use the following notation:

\bnotal{H.N.L}
 Let \(\seqska \) be a sequence of complex \tpqa{matrices}.
 Then let \(\Lu{0}\defg\Hu{0}\) and let \(\Lu{n}\defg\Hu{n}\schca \Hu{n-1}\) for all \(n\in\N\) with \(2n\leq\kappa\).
\enota

\breml{H.R.h2L}
 Let \(\seqska\) be a sequence of complex \tpqa{matrices}.
 For all \(n\in\NO\) with \(2n\leq\kappa\), then \(\hpu{2n}=\Lu{n}\).
 In particular, if \(n\geq1\), then \(\hpu{2n}\) is the Schur complement \(\Hu{n}\schca \Hu{n-1}\) of \(\Hu{n-1}\) in \(\Hu{n}\).
\erem

\bleml{H.L.Hgg-sr}
 Let \(n\in\NO\) and let \(\seqs{2n}\in\Hggqu{2n}\).
 Then:
\benui
 \il{H.L.Hgg-sr.a} \(\su{j}\in\CHq\) for all \(j\in\mn{0}{2n}\) and \(\su{2k}\in\Cggq\) for all \(k\in\mn{0}{n}\).
 \il{H.L.Hgg-sr.b} \(\ran{\su{2k+1}}\subseteq\ran{\su{2\ell}}\) for all \(k\in\mn{0}{n-1}\) and all \(\ell\in\mn{0}{n}\).
 \il{H.L.Hgg-sr.c} If \(n\geq2\), then
 \beql{H.L.Hgg-sr.A}
  \ran{\su{0}}
  \supseteq\ran{\su{2}}
  =\ran{\su{4}}
  =\dotsb
  =\ran{\su{2n-2}}
  \subseteq\ran{\su{2n}}.
 \eeq
\eenui
\elem
\bproof
 \rPart{H.L.Hgg-sr.a} was already proved in~\zitaa{MR3014199}{\clem{3.2}{126}}.
 First observe that each principal submatrix of the \tnnH{} matrix \(\Hu{n}\) is \tnnH{} itself.
 Assume \(n\geq1\).
 For each \(k\in\mn{0}{n-1}\), we consider the principal submatrix \(\tmat{\su{2k}&\su{2k+1}\\\su{2k+1}&\su{2k+2}}\) of \(\Hu{n}\).
 \rlem{L.AEP} yields then
 \begin{align}\label{H.L.Hgg-sr.1}
  \ran{\su{2k+1}}&\subseteq\ran{\su{2k}}&
  &\text{and}&
  \ran{\su{2k+1}}&\subseteq\ran{\su{2k+2}}.
 \end{align}
 Now, assume \(n\geq2\).
 The application of \rlem{L.AEP} to the principal submatrices \(\tmat{\su{0}&\su{2}\\\su{2}&\su{4}}\) and \(\tmat{\su{2n-4}&\su{2n-2}\\\su{2n-2}&\su{2n}}\) of \(\Hu{n}\) yields \(\ran{\su{2}}\subseteq\ran{\su{0}}\) and \(\ran{\su{2n-2}}\subseteq\ran{\su{2n}}\).
 Now, assume \(n\geq3\).
 For each \(k\in\mn{1}{n-2}\), we consider the principal submatrices \(\tmat{\su{2k-2}&\su{2k}\\\su{2k}&\su{2k+2}}\) and \(\tmat{\su{2k}&\su{2k+2}\\\su{2k+2}&\su{2k+4}}\) of \(\Hu{n}\).
 From \rlem{L.AEP} we conclude then \(\ran{\su{2k}}\subseteq\ran{\su{2k+2}}\) and \(\ran{\su{2k+2}}\subseteq\ran{\su{2k}}\).
 Consequently, we have
 \(
  \ran{\su{2}}
  =\ran{\su{4}}
  =\dotsb
  =\ran{\su{2n-2}}
 \)
 also in the case \(n\geq2\).
 Hence, \eqref{H.L.Hgg-sr.A} is proved.
 From \eqref{H.L.Hgg-sr.1} and \eqref{H.L.Hgg-sr.A} we now conclude \(\ran{\su{2k+1}}\subseteq\ran{\su{2\ell}}\) for all \(k\in\mn{0}{n-1}\) and all \(\ell\in\mn{0}{n}\).
\eproof

\breml{H.R.Hgg-s}
 Let \(\kappa\in\NOinf\) and let \(\seqs{2\kappa} \in \Hggqu{2\kappa}\).
 In view of \rlem{H.L.Hgg-sr}, then \(\su{j}\in\CHq\) for all \(j\in\mn{0}{2\kappa}\) and, furthermore, \(\su{2k}\in\Cggq\) for all \(k\in\mn{0}{\kappa}\).
\erem

\bpropl{H.R.Hgg+Hgg}
 Let \(\seqs{2\kappa}\) and \(\seqt{2\kappa}\) be sequences belonging to \(\Hggqu{2\kappa}\) with \thpf{s} \(\seq{\frx_j}{j}{0}{2\kappa}\) and \(\seq{\fry_j}{j}{0}{2\kappa}\), \tresp{}
 Denote by \(\seq{\frz_j}{j}{0}{2\kappa}\) the \thpfa{\(\seq{s_j+t_j}{j}{0}{2\kappa}\)}.
 Then \(\seq{s_j+t_j}{j}{0}{2\kappa}\in\Hggqu{2\kappa}\) and, for all \(k\in\mn{0}{\kappa}\), furthermore \(\frz_{2k}\lgeq\frx_{2k}+\fry_{2k}\).
\eprop
\bproof
 Combine \rremss{H.R.H+H}{H.R.h2L} and \rprop{A.R.sc+sc}.
\eproof

\bpropnl{\tcf{}~\zitaa{MR3014199}{\cprop{3.11}{128}}}{H.R.A.Hgg.A^*}
 Let \(A\in\Cpq\), let \(\seqs{2\kappa}\in\Hggqu{2\kappa}\) with \thpf{} \(\hpseqo{2\kappa}\), and let \(\seq{\frt_j}{j}{0}{2\kappa}\) be the \thpfa{\(\seq{A\su{j}A^\ad}{j}{0}{2\kappa}\)}.
 Then \(\seq{A\su{j}A^\ad}{j}{0}{2\kappa}\in\Hgguu{p}{2\kappa}\) and \(\frt_{2k}\lgeq A\hpu{2k}A^\ad\) for all \(k\in\mn{0}{\kappa}\).
\eprop

\section{Matricial $\rhl$-Stieltjes moment sequences and \hkp{s}}\label{K.S}
 Let \(\ug,\obg\in\R\).
 In this section we summarize some results on \(\rhl\)\nobreakdash-Stieltjes moment sequences and \(\lhl\)\nobreakdash-Stieltjes moment sequences which are important for our subsequent considerations.
 These results are mostly taken from~\zita{MR3014201,MR3133464,MR3611479}.
 We concentrate mainly on the case of an interval \(\rhl\).
The dual case of an interval \(\lhl\) can be reduced to the former case applied to the interval \([-\obg,\infp)\) by a reflection procedure.

 In our following considerations let \(\ug\) be a fixed real number.
\bnotal{K.N.sa}
 Let \(\seqska\) be a sequence of complex \tpqa{matrices} and assume \(\kappa\geq1\).
 Then let the sequence \(\seqsa{\kappa-1}\) be given by
 \(
  \sa{j}
  \defg-\ug\su{j}+\su{j+1}
 \).
 For each matrix \(X_k=X_k^{\langle s\rangle}\) built from the sequence \(\seqska\), denote (if possible) by \(X_{\ug,k}\defeq X_k^{\langle a\rangle}\) the corresponding matrix built from the sequence \(\seqsa{\kappa-1}\) instead of \(\seqska\).
\enota

 In view of \rnota{N.HKG}, we get in particular \(\Hau{n}=-\ug\Hu{n}+\Ku{n}\) for all \(n\in\NO\) with \(2n+1\leq\kappa\).
 In the classical case \(\ug=0\), we have furthermore \(\sa{j}=\su{j+1}\) for all \(j\in\mn{0}{\kappa-1}\).
 
 We start with some classes of sequences of complex matrices corresponding to solvability criteria for matricial Stieltjes moment problems on the half-line \(\rhl\).
 Let \(\Kggqu{0}\defg\Hggqu{0}\).
 For each \(n\in\N\), denote by \(\Kggqu{2n}\) the set of all sequences \(\seqs{2n}\) of complex \tqqa{matrices} for which the \tbHms{} \(\Hu{n}\) and \(\Hau{n-1}\) are both \tnnH{}.
 For each \(n\in\NO\), denote by \(\Kggqu{2n+1}\) the set of all sequences \(\seqs{2n+1}\) of complex \tqqa{matrices} for which the \tbHms{} \(\Hu{n}\) and \(\Hau{n}\) are both \tnnH{}.
 Furthermore, denote by \(\Kggqinf\) the set of all sequences \(\seqsinf\) of complex \tqqa{matrices} satisfying \(\seqs{m}\in\Kggqu{m}\) for all \(m\in\NO\).
 The sequences belonging to \(\Kggqu{0}\), \(\Kggqu{2n}\), \(\Kggqu{2n+1}\), or \(\Kggqinf\) are said to be \emph{\tKnnd}.
 In~\zita{MR3014201} and subsequent papers of the corresponding authors, the sequences belonging to \(\Kggqu{\kappa}\) were called \emph{\(\ug\)\nobreakdash-Stieltjes (right-sided) non-negative definite}.
 Our terminology here differs for the sake of consistency.

\bpropnl{\zitaa{MR2735313}{\cthm{1.4}{909}}}{dm2.T11}
 Let \(m\in\NO\) and let \(\seqs{m}\) be a sequence of complex \tqqa{matrices}.
 Then \(\MggqKskg{m}\neq\emptyset\) if and only if \(\seqs{m}\in\Kggquu{m}{\ug}\).
\eprop

 The characterization of the solvability of the ``\(=\)''\nobreakdash-version of the moment problem associated with \(\Omega=\rhl\) requires additional tools which will be prepared now.
 For each \(m\in\NO\), denote by \(\Kggequ{m}\) the set of all sequences \(\seqs{m}\) of complex \tqqa{matrices} for which there exists a complex \tqqa{matrix} \(\su{m+1}\) such that the sequence \(\seqs{m+1}\) belongs to \(\Kggqu{m+1}\).
 Furthermore, let \(\Kggeqinf\defeq\Kggqinf\).
 The sequences belonging to \(\Kggequ{m}\) or \(\Kggeqinf\) are said to be \emph{\tKnnde}.

\bpropnl{\zitaa{MR3014201}{\cthm{1.6}{214}}}{dm2.T10}
 Let \(\seqska \) be a sequence of complex \tqqa{matrices}.
 Then \(\MggqKsg{\kappa}\neq\emptyset\) if and only if \(\seqska \in\Kggequu{\kappa}{\ug}\).
\eprop

 For the convenience of the reader we mention that we use in the following definition the \(\Theta\)'s introduced in \rnota{N.Lambda}.

\bdefnl{K.D.kpf}
 Let \(\seqska \) be a sequence of complex \tpqa{matrices}.
 Let the sequence \(\seq{\kpu{j}}{j}{0}{\kappa}\) be given by
 \(
  \kpu{2k}
  \defg\su{2k}-\Tripu{k}
 \)
 for all \(k\in\NO\) with \(2k\leq\kappa\) and by
 \(
  \kpu{2k+1}
  \defg\sau{2k}-\Tripa{k}
 \)
 for all \(k\in\NO\) with \(2k+1\leq\kappa\).
 Then we call \(\seq{\kpu{j}}{j}{0}{\kappa}\) the \emph{\tkpfa{\(\seqska \)}}.
\edefn

\breml{K.R.kpftr}
 Let \(\seqska \) be a sequence of complex \tpqa{matrices} with \tkpf{} \(\kpseqka \).
 For each \(k\in\mn{0}{\kappa}\), the matrix \(\kpu{k}\) is built from the matrices \(\su{0},\su{1},\dotsc,\su{k}\).
 In particular, for each \(m\in\mn{0}{\kappa}\), the \tkpfa{\(\seqs{m}\)} coincides with \(\kpseq{m}\).
\erem

 It is readily checked that a sequence \(\seqska\) belonging to the set \(\seqset{\kappa}{\Cpq}\) given in \eqref{seqset} can be recursively reconstructed from its \tkpf{} \(\kpseqka\) via \(\su{2k}=\Tripu{k}+\kpu{2k}\) for all \(k\in\NO\) with \(2k\leq\kappa\) and \(\su{2k+1}=\ug\su{2k}+\Tripa{k}+\kpu{2k+1}\) for all \(k\in\NO\) with \(2k+1\leq\kappa\).
 Furthermore, given an arbitrary sequence \(\kpseqka\in\seqset{\kappa}{\Cpq}\), we can build recursively, in the above mentioned way, a sequence \(\seqska\in\seqset{\kappa}{\Cpq}\) such hat \(\kpseqka\) is the \tkpfa{\(\seqska\)}.
 Consequently, we have:

\breml{K.R.kpfCqq}
 The mapping \(\Phi_\ug\colon\seqset{\kappa}{\Cpq}\to\seqset{\kappa}{\Cpq}\) defined by \(\seqska\mapsto\kpseqka\) is bijective.
\erem

 The \tkp{s} were introduced in~\zita{MR3014201} under the name \emph{right-sided \(\ug\)\nobreakdash-Stieltjes parametrization}.
 Analogous to the \thp{s}, several properties connected with \tKnnd{ness} can be characterized in terms of the \tkp{s}, the main result being the following:
 
\bpropnl{\tcf{}~\zitaa{MR3014201}{\cthm{4.12(b)}{225}}}{K.P.KggKP}
 A sequence \(\seqska \) of complex \tqqa{matrices} belongs to \(\Kggqka\) if and only if \(\kpu{j}\in\Cggq\) for all \(j\in\mn{0}{\kappa}\) and \(\Bilda{\kpu{j+1}}\subseteq \Bilda{\kpu{j}}\) for all \(j\in\mn{0}{\kappa-2}\).
\eprop

 We have the following connection to the \thpf{} given in \rdefn{102.HPN}:

\breml{K.R.kphp}
 Let \(\seqska \) be a sequence of complex \tpqa{matrices} with \tkpf{} \(\kpseqka \) and \thpf{} \(\hpseqo{\kappa}\).
 In the case \(\kappa\geq1\), denote by \(\hpaseqo{\kappa-1}\) the \thpfa{the} sequence \(\seqsa{\kappa-1}\) given in \rnota{K.N.sa}.
 Then \(\kpu{2k}=\hpu{2k}\) for all \(k\in\NO\) with \(2k\leq\kappa\) and \(\kpu{2k+1}=\hpau{2k}\) for all \(k\in\NO\) with \(2k+1\leq\kappa\).
\erem

 Using \rrem{K.R.kphp}, we obtain from \rremss{H.R.h2L}{102.MV} the following two remarks:

\breml{K.R.k2L}
 Let \(\seqska\) be a sequence of complex \tpqa{matrices}.
 Then \(\kpu{2n}=\Lu{n}\) for all \(n\in\NO\) with \(2n\leq\kappa\) and \(\kpu{2n+1}=\Lau{n}\) for all \(n\in\NO\) with \(2n+1\leq\kappa\).
 In particular, if \(n\geq1\), then \(\kpu{2n}\) is the Schur complement \(\Hu{n}\schca \Hu{n-1}\) of \(\Hu{n-1}\) in \(\Hu{n}\) and \(\kpu{2n+1}\) is the Schur complement \(\Hau{n}\schca \Hau{n-1}\) of \(\Hau{n-1}\) in \(\Hau{n}\).
\erem

 Using \rrem{K.R.kphp}, we obtain from \rrem{H.R.l*h} immediately:

\breml{K.R.l*k}
 Let \(\lambda\in\C\) and let \(\seqska\) be a sequence of complex \tpqa{matrices}.
 Then \(\seq{\lambda\kpu{j}}{j}{0}{\kappa}\) coincides with the \tkpfa{\(\seq{\lambda s_j}{j}{0}{\kappa}\)}.
\erem

 Since we consider the \tkpf{} for real numbers \(\ug\) only, we can obtain analogous results to \rrem{H.R.d^j*h} in two different ways:

\bleml{K.L.-1^j*k}
 Let \(\seqska \) be a sequence of complex \tpqa{matrices}.
 Let the sequence \(\seqr{\kappa}\) be given by \(r_j\defg(-1)^j\su{j}\).
 Then \(\seq{(-1)^j\kpu{j}}{j}{0}{\kappa}\) coincides with the \tkupfa{-\ug}{\(\seqr{\kappa}\)}.
\elem
\bproof
 We consider the nontrivial case \(\kappa\geq1\) only.
 Denote by \(\seq{\frr_j}{j}{0}{\kappa}\) the \tkupfa{-\ug}{\(\seqr{\kappa}\)}.
 Because of \rrem{K.R.kphp}, then we conclude \(\frr_{2k}=(-1)^{2k}\kpu{2k}\) for all \(k\in\NO\) with \(2k\leq\kappa\) from \rrem{H.R.d^j*h}.
 Observe that \(-(-\ug)r_j+r_{j+1}=-(-1)^j\sau{j}\) for all \(j\in\mn{0}{\kappa-1}\) holds true.
 Hence, by virtue of \rrem{K.R.kphp}, we get \(\frr_{2k+1}=-(-1)^{2k}\kpu{2k+1}\) for all \(k\in\NO\) with \(2k+1\leq\kappa\), using \rremss{H.R.l*h}{H.R.d^j*h}.
 Hence, we have \(\frr_j=(-1)^j\kpu{j}\) for all \(j\in\mn{0}{\kappa}\).
\eproof

\bleml{K.L.d^j*k}
 Let \(\xi\in\C\) with \(\abs{\xi}=1\) and let \(\seqska\) be a sequence of complex \tpqa{matrices} with \tkupf{0} \(\kpseqka\).
 Then \(\seq{\xi^j\kpu{j}}{j}{0}{\kappa}\) coincides with the \tkupfa{0}{\(\seq{\xi^j s_j}{j}{0}{\kappa}\)}.
\elem
\bproof
 We consider the non-trivial case \(\kappa\geq1\) only.
 Let the sequence \(\seqr{\kappa}\) be given by \(r_j\defg\xi^j\su{j}\) and denote by \(\seq{\frr_j}{j}{0}{\kappa}\) the \tkupfa{-0}{\(\seqr{\kappa}\)}.
 Observe, that \(-(-0)r_j+r_{j+1}=\xi\xi^j(-0\cdot\su{j}+\su{j+1})\) holds true for all \(j\in\mn{0}{\kappa-1}\).
 For all \(j\in\mn{0}{\kappa}\), we thus obtain \(\frr_j=\xi^j\kpu{j}\) by the same reasoning as in the proof of \rlem{K.L.-1^j*k}.
\eproof

 Now we consider for an arbitrarily fixed number \(\obg\in\R\) a class of sequences of complex matrices, which turns out to be dual to the class \(\Kggqka\) and which takes the corresponding role for the moment problem on the left half-line \(\lhl\).
 (This class was already studied in~\zita{MR3014201}.):
 Let \(\Lggqu{0}\defg\Hggqu{0}\).
 For each \(n\in\N\), denote by \(\Lggqu{2n}\) the set of all sequences \(\seqs{2n}\) of complex \tqqa{matrices} for which the \tbHms{} \(\Hu{n}\) and \(\obg\Hu{n-1}-\Ku{n-1}\) are both \tnnH{}.
 For each \(n\in\NO\), denote by \(\Lggqu{2n+1}\) the set of all sequences \(\seqs{2n+1}\) of complex \tqqa{matrices} for which the \tbHms{} \(\Hu{n}\) and \(\obg\Hu{n}-\Ku{n}\) are both \tnnH{}.
 Furthermore, denote by \(\Lggqinf\) the set of all sequences \(\seqsinf\) of complex \tqqa{matrices} satisfying \(\seqs{m}\in\Lggqu{m}\) for all \(m\in\NO\).
 The sequences belonging to \(\Lggqu{0}\), \(\Lggqu{2n}\), \(\Lggqu{2n+1}\), or \(\Lggqinf\) are said to be \emph{\tLnnd}.
 (Note that in~\zita{MR3014201} the sequences belonging to \(\Lggqu{\kappa}\) were called \emph{\(\obg\)\nobreakdash-Stieltjes left-sided non-negative definite}.)
 For each \(m\in\NO\), denote by \(\Lggequ{m}\) the set of all sequences \(\seqs{m}\) of complex \tqqa{matrices} for which there exists a complex \tqqa{matrix} \(\su{m+1}\) such that the sequence \(\seqs{m+1}\) belongs to \(\Lggqu{m+1}\).
 Furthermore, let \(\Lggeqinf\defeq\Lggqinf\).
 The sequences belonging to \(\Lggequ{m}\) or \(\Lggeqinf\) are said to be \emph{\tLnnde}.

\bpropnl{\zitaa{MR3014201}{\cthm{1.8}{215}}}{121.T1609}
 Let \(\obg\in\R\) and let \(\seqska \) be a sequence of complex \tqqa{matrices}.
 Then \(\MggqLsg{\kappa}\neq\emptyset\) if and only if \(\seqska\in\Lggequu{\kappa}{\obg}\).
\eprop

 The proof of \rprop{121.T1609} is based on a combination of \rprop{dm2.T10} applied to the interval \([-\obg,\infp)\) with a reflection procedure as in \rlem{K.L.-1^j*k}.

\section{Matricial $\ab$-Hausdorff moment sequences}\label{F.S}
 At the beginning of this section we recall a first collection of results on the matricial Hausdorff moment problem which are mostly taken from~\zita{arXiv:1701.04246}.
 In the remaining part of this paper, let \(\ug\) and \(\obg\) be two arbitrarily given real numbers satisfying \(\ug<\obg\).
 We now start considering moment sequences of \tnnH{} \tqqa{measures} on the interval \(\ab\) which form the main object of this paper.
 To state a solvability criterion for the corresponding matricial \(\ab\)\nobreakdash-Hausdorff moment problem, we first extend \rnota{K.N.sa}:

\bnotal{F.N.sa}
 Let \(\seqska\) be a sequence of complex \tpqa{matrices} and assume \(\kappa\geq1\).
 Let the sequences \(\seqsa{\kappa-1}\) and \(\seqsb{\kappa-1}\) be given by \(\sa{j}\defg-\ug\su{j}+\su{j+1}\) and \(\sub{j}\defg\obg\su{j}-\su{j+1}\), \tresp{}
 Furthermore, if \(\kappa\geq2\), then let the sequence \(\seqsab{\kappa-2}\) be given by
 \(
  \sab{j}
  \defg-\ug\obg\su{j}+(\ug+\obg)\su{j+1}-\su{j+2}
 \).
 
 For each matrix \(X_k=X_k^{\langle s\rangle}\) built from the sequence \(\seqska\), we denote (if possible) by \(X_{\ug,k}\defeq X_k^{\langle a\rangle}\), by \(X_{k,\obg}\defeq X_k^{\langle b\rangle}\), and by \(X_{\ug,k,\obg}\defeq X_k^{\langle c\rangle}\) the corresponding matrix built from the sequence \(\seqsa{\kappa-1}\), from the sequence \(\seqsb{\kappa-1}\), and from the sequence \(\seqsab{\kappa-2}\) instead of \(\seqska\), \tresp{}
\enota

 In view of \rnota{N.HKG}, we get in particular \(\Hau{n}=-\ug\Hu{n}+\Ku{n}\) and \(\Hub{n}=\obg\Hu{n}-\Ku{n}\) for all \(n\in\NO\) with \(2n+1\leq\kappa\) and
\(
 \Hab{n}
 =-\ug\obg\Hu{n}+(\ug+\obg)\Ku{n}-\Gu{n}
\)
 for all \(n\in\NO\) with \(2n+2\leq\kappa\).
 In the classical case \(\ug=0\) and \(\obg=1\), we have furthermore \(\sa{j}=\su{j+1}\) and \(\sub{j}=\su{j}-\su{j+1}\) for all \(j\in\mn{0}{\kappa-1}\) and \(\sab{j}=\su{j+1}-\su{j+2}\) for all \(j\in\mn{0}{\kappa-2}\).
 The following observation is immediate from the definitions:

\breml{F.R.c=ab}
 We have \(\sab{j}=-\ug\sub{j}+\sub{j+1}\) and \(\sab{j}=\obg\sau{j}-\sau{j+1}\) for all \(j\in\mn{0}{\kappa-2}\).
\erem

 We denote by \(\Cgq\defeq\setaa{M\in\Cqq}{v^\ad Mv\in(0,\infp)\text{ for all }v\in\Cq\setminus\set{\Ouu{q}{1}}}\)
 the set of \tpH{} matrices from \(\Cqq\).
 Now we introduce those classes of (finite or infinite) sequences from \(\Cqq\) which are the central objects in this paper.
 Let \(\Fggqu{0}\) (\tresp{}, \(\Fgqu{0}\)) be the set of all sequences \((s_j)_{j=0}^0\) with \(\su{0}\in\Cggq\) (\tresp{}, \(\su{0}\in\Cgq\)).
 For each \(n\in\N\), denote by \(\Fggqu{2n}\) (\tresp{}, \(\Fgqu{2n}\)) the set of all sequences \(\seqs{2n}\) of complex \tqqa{matrices}, for which the \tbHms{} \(\Hu{n}\) and \(\Hab{n-1}\) are both non-negative (\tresp{}, positive) \tH{}.
 For each \(n\in\NO\), denote by \(\Fggqu{2n+1}\) (\tresp{}, \(\Fgqu{2n+1}\)) the set of all sequences \(\seqs{2n+1}\) of complex \tqqa{matrices} for which the \tbHms{} \(\Hau{n}\) and \(\Hub{n}\) are both non-negative (\tresp{}, positive) \tH{}.
 Furthermore, denote by \(\Fggqinf\) (\tresp{}, \(\Fgqinf\)) the set of all sequences \(\seqsinf \) of complex \tqqa{matrices} such that \(\seqs{m}\) belongs to \(\Fggqu{m}\) (\tresp{}, \(\Fgqu{m}\)) for all \(m\in\NO\).
 The sequences belonging to \(\Fggqu{0}\), \(\Fggqu{2n}\), \(\Fggqu{2n+1}\), or \(\Fggqinf\) (\tresp{}, \(\Fgqu{0}\), \(\Fgqu{2n}\), \(\Fgqu{2n+1}\), or \(\Fgqinf\)) are said to be \emph{\(\ab\)\nobreakdash-non-negative} (\tresp{}, \emph{positive}) \emph{definite}.
 (Note that in~\zita{arXiv:1701.04246} the sequences belonging to \(\Fggqu{\kappa}\) (\tresp{}, \(\Fgqu{\kappa}\)) were called \emph{\(\ab\)\nobreakdash-Hausdorff non-negative} (\tresp{}, \emph{positive}) \emph{definite}.)
 The following result indicates the importance of the just introduced objects.

\bthmnl{\tcf{}~\zitaa{MR2222521}{\cthm{1.3}{127}} and~\zitaa{MR2342899}{\cthm{1.3}{106}}}{I.P.ab}
 Let \(\seqska\) be a sequence of complex \tqqa{matrices}.
 Then \(\MggqFksg\neq\emptyset\) if and only if \(\seqska\in\Fggqka\).
\ethm

 Since \(\Omega=\ab\) is bounded, one can easily see that \(\MggqF=\MggquF{\infi}\), \tie{}\ the power moment \eqref{I.G.mom} exists for all \(j\in\NO\).

 Let \(\sigma\in\MggqF\).
 Then we call \(\seqmpm{\sigma}\) given by \(\mpm{\sigma}{j}\defeq\int_\ab x^j\sigma\rk{\dif x}\) the \emph{\tfpmfa{\(\sigma\)}}.

 Given the complete sequence of prescribed power moments \(\seqsinf\), the moment problem on the compact interval \(\Omega=\ab\) differs from the moment problems on the unbounded sets \(\Omega=\rhl\) and \(\Omega=\R\) in having necessarily a unique solution, assumed that a solution exists:

\bpropl{I.P.ab8}
 If \(\seqsinf\in\Fggqinf\), then the set \(\MggqFsg{\infi}\) consists of exactly one element.
\eprop

 We can summarize \rprop{I.P.ab8} and \rthm{I.P.ab} for \(\kappa=\infi\):

\bpropl{I.P.ab8Fgg}
 The mapping \(\Xi_\ab\colon\MggqF\to\Fggqinf\) given by \(\sigma\mapsto\seqmpm{\sigma}\) is well defined and bijective.
\eprop

 Observe that, for each \(\ell\in\NO\), the \emph{\(\ell\)\nobreakdash-th moment space}
\beql{msp}
 \MomqF{\ell}
 \defeq\setaca*{\col\rk{\mpm{\sigma}{0},\mpm{\sigma}{1},\dotsc,\mpm{\sigma}{\ell}}}{\sigma\in\MggqF}
\eeq
 \emph{on \(\ab\)} corresponds to the class \(\Fggqu{\ell}\) via
\[
 \MomqF{\ell}
 =\setaca*{\yuu{0}{\ell}}{\seqs{\ell}\in\Fggqu{\ell}},
\]
 whereas its interior \(\intr\rk{\MomqF{\ell}}\) corresponds to \(\Fgqu{\ell}\) via
\[
 \intr\rk*{\MomqF{\ell}}
 =\setaca*{\yuu{0}{\ell}}{\seqs{\ell}\in\Fgqu{\ell}}
\]
 (\tcf{}~\zitaa{MR1883272}{\cthm{2.6}{177}}).
 Here \(\yuu{0}{\ell}\) is given in \rnota{N.yz} .

 Obviously, we have
\beql{Fgg2n}
 \Fggqu{2n}
 =\setaa*{\seq{\su{j}}{j}{0}{2n}\in\Hggqu{2n}}{\seq{\sab{j}}{j}{0}{2(n-1)}\in\Hggqu{2(n-1)}}
\eeq
 for all \(n\in\N\) and
\beql{Fgg2n+1}
 \Fggqu{2n+1}
 =\setaa*{\seq{\su{j}}{j}{0}{2n+1}\in\seqset{2n+1}{\Cqq}}{\set*{\seqsa{2n},\seqsb{2n}}\subseteq\Hggqu{2n}}
\eeq
 for all \(n\in\NO\).
 The following result can also be easily obtained from \rthm{I.P.ab}:

\bpropnl{\zitaa{arXiv:1701.04246}{\cprop{9.1}{27}}}{F.L.sabF}
 Let \(\seqska\in\Fggqka\).
 If \(\kappa\geq1\), then \(\set{\seqsa{\kappa-1},\seqsb{\kappa-1}}\subseteq\Fggqu{\kappa-1}\).
 If \(\kappa\geq2\), furthermore \(\seqsab{\kappa-2}\in\Fggqu{\kappa-2}\).
\eprop

\bleml{F.R.Fgg-s}
 Let \(\seqska  \in \Fggqka\).
 Then \(\su{j}\in\CHq\) for all \(j\in\mn{0}{\kappa}\) and \(\su{2k}\in\Cggq\) for all \(k\in\NO\) with \(2k\leq\kappa\).
 Furthermore, \(\ug\su{2k}\lleq\su{2k+1}\lleq\obg\su{2k}\) for all \(k\in\NO\) with \(2k+1\leq\kappa\).
\elem
\bproof
 In view of~\zitaa{arXiv:1701.04246}{\cprop{7.7(a)}{18}} and \eqref{Fgg2n}, \eqref{Fgg2n+1}, and \(\ug,\obg\in\R\), this can be seen from \rrem{H.R.Hgg-s}.
\eproof

 Particular relations between the moment problems on \(\ab\), \(\rhl\), and \(\lhl\) are expressed in the following two results:

\bpropnl{\zitaa{arXiv:1701.04246}{\cprop{9.2}{28}}}{ab.R1532}
 Let \(m\in\NO\).
 Then \(\Fggqu{m}\subseteq\Kggequ{m}\cap\Lggequ{m}\).
\eprop

\bpropnl{\tcf{}~\zitaa{arXiv:1701.04246}{\cprop{8.24}{26}}}{ab.R1535}
 \(\Fggqu{2\kappa+1}=\Kggqu{2\kappa+1}\cap\Lggqu{2\kappa+1}\).
\eprop

 Observe that in general the sets \(\Fggqu{2n}\) and \(\Kggequ{2n}\cap\Lggequ{2n}\) do not necessarily coincide (\tcf{}~\zitaa{arXiv:1701.04246}{\cexa{1}{28}}).
 In comparison to \tKnnd{} sequences, the inner structure of \tFnnd{} sequences is more complicated.
 On the other hand, finite sequences from \(\Fggqu{m}\) can always be extended to sequences from \(\Fggqu{\ell}\) for all \(\ell\in\minf{m+1}\), which is due to the fact that a \tnnH{} measure on the bounded set \(\ab\) possesses power moments of all non-negative orders.
 We first introduce several matrices and recall their role in the corresponding extension problem for \tFnnd{} sequences, studied in~\zita{arXiv:1701.04246}.
 
 For two \tH{} \tqqa{matrices} \(A\) and \(B\) with \(A\lleq B\), the closed matricial interval \(\matint{A}{B}\defeq\setaca{X\in\CHq}{A\lleq X\lleq B}\) is non-empty.
 We write \(A\lls B\) or \(B\lgs A\) if \(A,B\in\CHq\) and \(B-A\in\Cgq\) are fulfilled, in which case the open matricial interval \(\matinto{A}{B}\defeq\setaca{X\in\CHq}{A\lls X\lls B}\) is non-empty.
 For each matrix \(A\in\Cggq\), there exists a uniquely determined matrix \(Q\in\Cggq\) with \(Q^2=A\) called the \emph{\tnnH{} square root} \(Q=A^\varsqrt \) of \(A\).
 
\blemnl{\zitaa{arXiv:1701.04246}{\clem{10.1}{29}}}{165.L101}
 Let \(A,B\in\CHq\) and let \(D\defeq B-A\).
 Then:
\benui
 \il{165.L101.a} \(\matint{A}{B}\neq\emptyset\) if and only if \(D\) is \tnnH{}.
 In this case,
\(
 \matint{A}{B}
 =\setaca{A+ D^\varsqrt K D^\varsqrt }{K\in\matint{\Oqq}{\Iq}}
\).
 \il{165.L101.b} \(\matinto{A}{B}\neq\emptyset\) if and only if \(D\) is \tpH{}.
 In this case,
\(
 \matinto{A}{B}
 =\setaca{A+ D^\varsqrt K D^\varsqrt }{K\in\matinto{\Oqq}{\Iq}}
\).
\eenui
\elem

 From \rlem{165.L101} we immediately obtain the following observations:
 
\breml{F.R.int}
 Let \(A,B\in\CHq\) and let \(D\defeq B-A\).
 Then:
\benui
 \il{F.R.int.a} Let \(D\in\Cggq\) and let \(\lambda\in[0,1]\), then \(A+\lambda D\in\matint{A}{B}\).
 \il{F.R.int.b} Let \(D\in\Cgq\) and let \(\lambda\in(0,1)\), then \(A+\lambda D\in\matinto{A}{B}\).
\eenui
\erem

 In the scalar case \(q=1\) the intervals \(\matint{A}{B}\) and \(\matinto{A}{B}\) are obviously exhausted by the set of points described in \rpartss{F.R.int.a}{F.R.int.b} of \rrem{F.R.int}, \tresp{}
 It is easily checked that in the case \(q>1\) the intervals \(\matint{A}{B}\) and \(\matinto{A}{B}\) obtain also other points.
 
 One of the main results in~\zita{arXiv:1701.04246} states that the possible one-step extensions \(\su{m+1}\in\Cqq\) of a sequence \(\seqs{m}\) belonging to \(\Fggqu{m}\) (\tresp{}, \(\Fgqu{m}\)) to an \(\ab\)\nobreakdash-non-negative (\tresp{}, positive) definite sequence \(\seqs{m+1}\) fill out a matricial intervals.
 In order to give an exact description we need a little preparation.
 
\bdefnl{F.N.umom}
 Let \(\seqska\) be a sequence of complex \tpqa{matrices}.
 Then the sequences \(\seq{\umg{j}}{j}{0}{\kappa}\) and \(\seq{\omg{j}}{j}{0}{\kappa}\) given by \(  \umg{2k}\defg\ug\su{2k}+\Tripa{k}\) and \(\omg{2k}\defg\obg\su{2k}-\Tripb{k}\) for all \(k\in\NO\) with \(2k\leq\kappa\) and by \(  \umg{2k+1}\defg\Tripu{k+1}\) and \(\omg{2k+1}\defg-\ug\obg\su{2k}+(\ug+\obg)\su{2k+1}-\Tripab{k}\) for all \(k\in\NO\) with \(2k+1\leq\kappa\) are called the \emph{\tflep{\(\seqska\)}} and the \emph{\tfrep{\(\seqska\)}}.
\edefn

 By virtue of \rnota{N.Lambda}, we have in particular
\begin{align}
 \umg{0}&=\ug\su{0},&
 \omg{0}&=\obg\su{0}\label{F.G.uo0}
 \shortintertext{and}
 \umg{1}&=\su{1}\su{0}^\mpi\su{1},&
 \omg{1}&=-\ug\obg\su{0}+(\ug+\obg)\su{1}.\label{F.G.uo1}
\end{align}

 Using \rlem{F.R.Fgg-s} and \rrem{A.R.A++*}, we easily obtain:

\breml{ab.L0911}
 Let \(\seqska  \in \Fggqka\).
 Then \(\set{\umg{j},\omg{j}}\subseteq\CHq\) for all \(j\in\mn{0}{\kappa}\).
\erem

 The following result is of central importance for our subsequent considerations.

\bthmnl{\tcf{}~\zitaa{arXiv:1701.04246}{\cthm{11.2}{44}}}{165.T112}
 Let \(m\in\NO\) and let \(\seqs{m}\) belong to \(\Fggqu{m}\) (\tresp{}, \(\Fgqu{m}\)).
 Then the matricial interval \(\matint{\umg{m}}{\omg{m}}\) (\tresp{}, \(\matinto{\umg{m}}{\omg{m}}\)) is non-empty and coincides with the set of all complex \tqqa{matrices} \(\su{m+1}\) for which \(\seqs{m+1}\) belongs to \(\Fggqu{m+1}\) (\tresp{}, \(\Fgqu{m+1}\)).
\ethm

 For later use, we recall several results and notations from~\zita{arXiv:1701.04246} supplemented by some additional observations:

\bdefnl{F.N.AB}
 Let \(\seqska\) be a sequence of complex \tpqa{matrices}.
 Then the sequence \(\seq{\usc{j}}{j}{0}{\kappa}\) given by \(\usc{0}\defeq\su{0}\) and by \(\usc{j}\defg\su{j}-\umg{j-1}\) is called the \emph{\tflsc{\(\seqska\)}}.
 Furthermore, if \(\kappa\geq1\), then the sequence \(\seq{\osc{j}}{j}{1}{\kappa}\) given by \(\osc{j}\defg\omg{j-1}-\su{j}\) is called the \emph{\tfusc{\(\seqska\)}}.
\edefn

 Because of \eqref{F.G.uo0} and  \eqref{F.G.uo1}, we have in particular
\begin{align}\label{F.G.AB1B2}
 \usc{1}&=\sau{0},&
 \osc{1}&=\sub{0},&
 &\text{and}&
 \osc{2}&=\sab{0}.
\end{align}
 In view of \rrem{H.R.h2L}, we obtain:

\breml{F.R.ABL}
 We have \(\usc{2n}=\Lu{n}\) (\tresp{}, \(\usc{2n+1}=\Lau{n}\)) for all \(n\in\NO\) with \(2n\leq\kappa\) (\tresp{}, \(2n+1\leq\kappa\)).
 In particular, if \(n\geq1\), then \(\usc{2n}\) is the Schur complement of \(\Hu{n-1}\) in \(\Hu{n}\) and \(\usc{2n+1}\) is the Schur complement of \(\Hau{n-1}\) in \(\Hau{n}\). 
 Furthermore, we have \(\osc{2n+1}=\Lub{n}\) (\tresp{}, \(\osc{2n+2}=\Lab{n}\)) for all \(n\in\NO\) with \(2n+1\leq\kappa\) (\tresp{}, \(2n+2\leq\kappa\)).
 In particular, if \(n\geq1\), then \(\osc{2n+1}\) is the Schur complement of \(\Hub{n-1}\) in \(\Hub{n}\) and \(\osc{2n+2}\) is the Schur complement of \(\Hab{n-1}\) in \(\Hab{n}\).
\erem

\bdefnl{F.D.dia}
 Let \(\seqska \) be a sequence of complex \tpqa{matrices}.
 Then we call \(\seq{\dia{j}}{j}{0}{\kappa}\) be given by
\(
  \dia{j}
  \defeq\omg{j}-\umg{j}
\)
 the \emph{\tfdfa{\(\seqska \)}}.
\edefn

 By virtue of \eqref{F.G.uo0} and \eqref{F.G.uo1}, we have in particular
\begin{align}\label{F.G.d01}
 \dia{0}&=\ba\su{0}&
 &\text{and}&
 \dia{1}&=-\ug\obg\su{0}+(\ug+\obg)\su{1}-\su{1}\su{0}^\mpi\su{1}.
\end{align}
 Observe that, in view of \rthm{165.T112}, the interval lengths \(\dia{j}\) can be understood as the widths of the sections of the moment space \eqref{msp}.
 From \rdefn{F.N.umom} we easily obtain:

\breml{F.R.diatr}
 Let \(\seqska\) be a sequence of complex \tpqa{matrices} with \tfdf{} \(\seqdiaka \).
 For each \(k\in\mn{0}{\kappa}\), the matrix \(\dia{k}\) is built  from the matrices \(\su{0},\su{1},\dotsc,\su{k}\).
 In particular, for each \(m\in\mn{0}{\kappa}\), the \tfdfa{\(\seqs{m}\)} coincides with \(\seqdia{m}\).
\erem

\breml{ab.R1420}
 Let \(\seqska\) be a sequence of complex \tpqa{matrices}.
 For all \(j\in\mn{0}{\kappa-1}\), then \(\dia{j}=\usc{j+1}+\osc{j+1}\).
\erem

 We now consider the behavior of the \tfdf{} under two elementary transformations of the underlying sequence \(\seqska \) of complex \tpqa{matrices}.
 In view of \(\ug,\obg\in\R\), the construction of the sequences \(\seqsa{\kappa-1}\), \(\seqsb{\kappa-1}\), and \(\seqsab{\kappa-2}\) from \rnota{F.N.sa} is compatible with the transformations in question.
 Thus, extending, in the case \(\kappa<\infi\), the sequence \(\seqska\) to a sequence \(\seqsinf\) of complex \tpqa{matrices} in an arbitrary way and taking into account \rremssss{H.R.h2L}{F.R.ABL}{F.R.diatr}{ab.R1420}, we conclude from \rremss{H.R.l*h}{H.R.U.h.V} immediately:

\breml{F.R.l*d}
 If \(\lambda\in\C\) and if \(\seqska\) is a sequence of complex \tpqa{matrices}, then \(\seq{\lambda\dia{j}}{j}{0}{\kappa}\) is exactly the \tfdfa{\(\seq{\lambda s_j}{j}{0}{\kappa}\)}.
\erem

\breml{F.R.U.d.V}
 Let \(U\in\Coo{u}{p}\) with \(U^\ad U=\Ip\), let \(V\in\Coo{q}{v}\) with \(VV^\ad=\Iq\), and let \(\seqska\) be a sequence of complex \tpqa{matrices}.
 Then \(\seq{U\dia{j}V}{j}{0}{\kappa}\) coincides with the \tfdfa{\(\seq{Us_jV}{j}{0}{\kappa}\)}.
\erem

 Looking back to~\zita{arXiv:1701.04246} we see that the parallel sum of matrices is an important tool for the description of several statements connected with the extension problem for \tFnnd{} sequences.
 For the convenience of the reader we recall the parallel sum and one of its main properties.
 If \(A\) and \(B\)  are two complex \tpqa{matrices}, then the matrix
\beql{ps}
  A\ps B
  \defg A(A+B)^\mpi B
\eeq
 is called the \emph{parallel sum of \(A\) and \(B\)}.

\bthmnl{\tcf{}~\zitaa{MR0325642}{\cthm{2.2(c)}{93}}}{ab.T0909}
 If \(A,B\in\Cggq\), then \(A\ps B\in\Cggq\).
\ethm

 The following results indicate the key role of the parallel sum in the framework of our extension problem.
 
\bthmnl{\zitaa{arXiv:1701.04246}{\cthm{10.14}{32}}}{ab.P1422}
 Let \(\seqska\in\Fggqka\).
 Then \(\dia{0}=\ba\usc{0}\) and, for all \(k\in\mn{1}{\kappa}\), furthermore \(\dia{k}=\ba(\usc{k}\ps\osc{k})\) and \(\dia{k}=\ba(\osc{k}\ps\usc{k})\).
\ethm

\bpropnl{\zitaa{arXiv:1701.04246}{\cprop{10.15}{37}}}{ab.C0929}
 If \(\seqska\) belongs to \(\Fggqka\) (\tresp{} \(\Fgqu{\kappa}\)), then \(\dia{j}\) is non-negative (\tresp{}, positive) \tH{} for all \(j\in\mn{0}{\kappa}\).
\eprop

 Now we state a first consequence of \rprop{ab.C0929}.

\bcorl{ab.R1011}
 Let \(m\in\NO\), let \(\seqs{m}\in\Fggqu{m}\), let \(\lambda\in[0,1]\), and let \(\su{m+1}\defeq\umg{m}+\lambda\dia{m}\).
 Then:
\benui
 \il{ab.R1011.a} The sequence \(\seqs{m+1}\) belongs to \(\Fggqu{m+1}\).
 \il{ab.R1011.b} The identities \(\usc{m+1}=\lambda\dia{m}\), \(\osc{m+1}=(1-\lambda)\dia{m}\), and \(\dia{m+1}=\ba\lambda(1-\lambda)\dia{m}\) hold true.
\eenui
\ecor
\bproof
 \eqref{ab.R1011.a} Obviously, \(\usc{m+1}=\su{m+1}-\umg{m}=\lambda\dia{m}\) and \(\osc{m+1}=\omg{m}-\su{m+1}=\omg{m}-\umg{m}-\lambda\dia{m}=(1-\lambda)\dia{m}\).
 In view of \rprop{ab.C0929}, we have \(\dia{m}\in\Cggq\).
 Thus, \rremp{F.R.int}{F.R.int.a} implies \(\su{m+1}\in\matint{\umg{m}}{\omg{m}}\).
 Because of \rthm{165.T112}, hence \(\seqs{m+1}\in\Fggqu{m+1}\).
 
 \eqref{ab.R1011.b} By virtue of \rthm{ab.P1422}, \rpart{ab.R1011.a}, \eqref{ps}, and \eqref{mpi}, we get with \(\delta\defeq\bam\) then
\[\begin{split}
 \dia{m+1}
 =\delta(\usc{m+1}\ps\osc{m+1})
 &=\delta\rk{\lambda\dia{m}}\ek*{\lambda\dia{m}+(1-\lambda)\dia{m}}^\mpi\ek*{(1-\lambda)\dia{m}}\\
 &=\delta\lambda(1-\lambda)\dia{m}\dia{m}^\mpi\dia{m}
 =\delta\lambda(1-\lambda)\dia{m}.\qedhere
\end{split}\]
\eproof

 Our next considerations can be summarized as follows:
 We consider a sequence \(\seqska\in\Fggqka\) and the matrices introduced in \rdefnss{F.N.AB}{F.D.dia}.
 Then we get a lot of interesting relations between ranges and null spaces of consecutive elements belonging to these matrix sequences.
 We write \(\nul{A}\defeq\setaa{x\in\Cq}{Ax=\Ouu{p}{1}}\) for the null space of a complex \tpqa{matrix} \(A\).

\bpropnl{\zitaa{arXiv:1701.04246}{\ccor{10.20}{39}}}{ab.C1101}
 Let \(\seqska\in\Fggqka\).
 For all \(j\in\mn{1}{\kappa}\), then
 \(\ran{\dia{j}}\subseteq\ran{\dia{j-1}}\), \(\nul{\dia{j-1}}\subseteq\nul{\dia{j}}\) and \(\ran{\usc{j}}\subseteq\ran{\usc{j-1}}\), \(\nul{\usc{j-1}}\subseteq\nul{\usc{j}}\).
 For all \(j\in\mn{2}{\kappa}\), furthermore \(\ran{\osc{j}}\subseteq\ran{\osc{j-1}}\) and \(\nul{\osc{j-1}}\subseteq\nul{\osc{j}}\).
\eprop

\bpropnl{\zitaa{arXiv:1701.04246}{\ccor{10.21}{39}}}{ab.C1343}
 If \(\seqska\in\Fggqka\), then \(\dia{j}=\ba\usc{j}\dia{j-1}^\mpi\osc{j}\) and \(\dia{j}=\ba\osc{j}\dia{j-1}^\mpi\usc{j}\) for all \(j\in\mn{1}{\kappa}\).
\eprop

 The geometry of a matricial interval suggests that in addition to the end points \(\umg{m}\) and \(\omg{m}\) of the extension interval in \rthm{165.T112}, the center of the matricial interval \(\matint{\umg{m}}{\omg{m}}\) is of particular interest:
 
\bdefnl{F.D.mi}
 Let \(\seqska \) be a sequence of complex \tpqa{matrices}.
 Then \(\seq{\mi{j}}{j}{0}{\kappa}\) given by \(\mi{j}\defeq\frac{1}{2}\rk{\umg{j}+\omg{j}}\) is called the \emph{\tfmfa{\(\seqska \)}}.
\edefn
 
 In view of \eqref{F.G.uo0} and \eqref{F.G.uo1}, we have in particular \(\mi{0}=\frac{\ug+\obg}{2}\su{0}\) and \(\mi{1}=\frac{1}{2}\ek{\su{1}\su{0}^\mpi\su{1}-\ug\obg\su{0}+(\ug+\obg)\su{1}}\).
 It should be mentioned that the choice \(\su{j+1}=\mi{j}\) corresponds to the maximization of the width \(\dia{j+1}\) of the corresponding section of the moment space \eqref{msp}:

\bpropnl{\zitaa{arXiv:1701.04246}{\cprop{10.23}{39}}}{ab.P1057}
 Let \(\seqska\in\Fggqka\).
 For all \(j\in\mn{0}{\kappa-1}\), then \(\dia{j+1}\lleq\frac{\bam}{4}\dia{j}\) with equality if and only if \(\su{j+1}=\mi{j}\).
\eprop

\section{The \hfpf{}}\label{F.S.fp}
 This section occupies a central position in this paper.
 It contains the realization of our main goal.
 For this we still need some preparations.
 By subsuming the \tkp{s} of the sequence \(\seqska\) and of the sequence \(\seqsb{\kappa-1}\) associated to it via \rnota{F.N.sa} to a new parameter sequence \(\fpseqka\) we are now led to one of the central objects of this paper:

\bdefnl{F.D.fpf}
 Let \(\seqska \) be a sequence of complex \tpqa{matrices}.
 Let the sequence \(\fpseqka\) be given by \(\fpu{0}\defg\usc{0}\), by \(\fpu{4k+1}\defeq\usc{2k+1}\) and \(\fpu{4k+2}\defeq\osc{2k+1}\) for all \(k\in\NO\) with \(2k+1\leq\kappa\), and by \(\fpu{4k+3}\defeq\osc{2k+2}\) and \(\fpu{4k+4}\defeq\usc{2k+2}\) for all \(k\in\NO\) with \(2k+2\leq\kappa\). 
 Then we call \(\fpseqka\) the \emph{\tfpfa{\(\seqska \)}}.
\edefn

 In view of \rdefn{F.N.AB}, \eqref{F.G.AB1B2}, \eqref{F.G.uo0}, and \eqref{F.G.uo1}, we have in particular
\begin{align}
 \fpu{0}&=\su{0},&
 \fpu{1}&=\sau{0}=\su{1}-\ug\su{0},&
 \fpu{2}&=\sub{0}=\obg\su{0}-\su{1}\label{F.G.f012}
\shortintertext{and}
 &&
 \fpu{3}&=\sab{0}=-\ug\obg\su{0}+(\ug+\obg)\su{1}-\su{2},&
 \fpu{4}&=\su{2}-\su{1}\su{0}^\mpi\su{1}.\label{F.G.f34}
\end{align}

 The following observation is immediate from the definitions:

\breml{F.R.fp12}
 For all \(m\in\mn{1}{\kappa}\), we have \(\set{\fpu{2m-1},\fpu{2m}}=\set{\usc{m},\osc{m}}\).
\erem

 By virtue of \rremss{F.R.ABL}{F.R.c=ab}, the connection to the \tkpf{} given in \rdefn{K.D.kpf} follows from \rrem{K.R.k2L}:

\breml{F.R.fpkap}
 Let \(\seqska\) be a sequence of complex \tpqa{matrices} with \tfpf{} \(\fpseqka\) and \tkpf{} \(\kpseq{\kappa}\).
 In the case \(\kappa\geq1\), denote by \(\lkpseq{\kappa-1}\) the \tkpfa{\(\seqsb{\kappa-1}\)}.
 Then \(\kpu{2k}=\fpu{4k}\) (\tresp{}\ \(\kpu{2k+1}=\fpu{4k+1}\)) for all \(k\in\NO\) with \(2k\leq\kappa\) (\tresp{}\ \(2k+1\leq\kappa\)) and \(\lkpu{2k}=\fpu{4k+2}\) (\tresp{}\ \(\lkpu{2k+1}=\fpu{4k+3}\)) for all \(k\in\NO\) with \(2k\leq\kappa-1\) (\tresp{}\ \(2k+1\leq\kappa-1\)).
\erem

 Now we state some observations on the arithmetics of the \tfp{s}.

\bleml{F.L.fpkbp}
 Let \(\seqska \) be a sequence of complex \tpqa{matrices} with \tfpf{} \(\fpseqka\) and \tkupf{\obg} \(\pkpseq{\kappa}\).
 If \(\kappa\geq1\), denote by \(\qkpseq{\kappa-1}\) the \tkupfa{\obg}{\(\seqsa{\kappa-1}\)}.
 Then \(\pkpu{j}=(-1)^j\fpu{2j}\) for all \(j\in\mn{0}{\kappa}\) and \(\qkpu{j}=(-1)^j\fpu{2j+1}\) for all \(j\in\mn{0}{\kappa-1}\).
\elem
\bproof
 Using \rremss{K.R.k2L}{F.R.ABL}, we obtain \(\pkpu{2k}=\Lu{k}=\usc{2k}=\fpu{4k}\) for all \(k\in\NO\) with \(2k\leq\kappa\).
 Now assume \(\kappa\geq1\).
 Let the sequence \(\seq{d_j}{j}{0}{\kappa-1}\) be given by \(d_j\defeq-\obg\su{j}+\su{j+1}\).
 Then \(d_j=-\sub{j}\) for all \(j\in\mn{0}{\kappa-1}\).
 Hence, in view of \rrem{ab.R1111}, we get from \rrem{K.R.k2L} furthermore \(\pkpu{2k+1}=\Luo{k}{d}=-\Lub{k}=-\osc{2k+1}=-\fpu{4k+2}\) for all \(k\in\NO\) with \(2k+1\leq\kappa\).
 Thus, \(\pkpu{j}=(-1)^j\fpu{2j}\) holds true for all \(j\in\mn{0}{\kappa}\).
 Using \rremss{K.R.k2L}{F.R.ABL}, we obtain \(\qkpu{2k}=\Lau{k}=\usc{2k+1}=\fpu{4k+1}\) for all \(k\in\NO\) with \(2k\leq\kappa-1\).
 Now assume \(\kappa\geq2\).
 Let the sequence \(\seq{e_j}{j}{0}{\kappa-2}\) be given by \(e_j\defeq-\obg\sau{j}+\sau{j+1}\).
 Because of \rrem{F.R.c=ab}, then \(e_j=-\sab{j}\) for all \(j\in\mn{0}{\kappa-2}\).
 Hence, in view of \rrem{ab.R1111}, we get from \rrem{K.R.k2L} furthermore \(\qkpu{2k+1}=\Luo{k}{e}=-\Lab{k}=-\osc{2k+2}=-\fpu{4k+3}\) for all \(k\in\NO\) with \(2k+1\leq\kappa-1\).
 Thus, \(\qkpu{j}=(-1)^j\fpu{2j+1}\) for all \(j\in\mn{0}{\kappa-1}\).
\eproof

\breml{F.R.fpftr}
 Let \(\seqska \) be a sequence of complex \tpqa{matrices} with \tfpf{} \(\fpseqka\).
 Then \(\fpu{0}=\su{0}\).
 Furthermore, in view of \rlem{F.L.fpkbp} and \rrem{K.R.kpftr}, for each \(k\in\mn{1}{\kappa}\), the matrices \(\fpu{2k-1}\) and \(\fpu{2k}\) are built from the matrices \(\su{0},\su{1},\dotsc,\su{k}\).
 In particular, for each \(m\in\mn{0}{\kappa}\), the \tfpfa{\(\seqs{m}\)} coincides with \(\seq{\fpu{j}}{j}{0}{2m}\).
\erem

\breml{F.R.fpfCqq}
 In view of \rremss{F.R.fpkap}{K.R.kpfCqq}, the mapping \(\Phi_{\ug,\obg}\colon\seqset{\kappa}{\Cpq}\to\seqset{2\kappa}{\Cpq}\) defined by \(\seqska\mapsto\fpseqka\) is injective.
\erem

 We now consider the behavior of the \tfpf{} under two elementary transformations of the underlying sequence \(\seqska \) of complex \tpqa{matrices}.
 In view of \(\ug,\obg\in\R\), the construction of the sequences \(\seqsa{\kappa-1}\), \(\seqsb{\kappa-1}\), and \(\seqsab{\kappa-2}\) from \rnota{F.N.sa} is compatible with the transformations in question.
 Thus, taking into account \rremss{H.R.h2L}{F.R.ABL} and \rdefn{F.D.fpf}, we conclude from \rremss{H.R.l*h}{H.R.U.h.V} immediately:

\breml{F.R.l*f}
 Let \(\lambda\in\C\) and let \(\seqska \) be a sequence of complex \tpqa{matrices}.
 Then \(\seq{\lambda\fpu{j}}{j}{0}{2\kappa}\) coincides with the \tfpfa{\(\seq{\lambda s_j}{j}{0}{\kappa}\)}.
\erem

\breml{F.R.U.f.V}
 Let \(U\in\Coo{u}{p}\) with \(U^\ad U=\Ip\), let \(V\in\Coo{q}{v}\) with \(VV^\ad=\Iq\), and let \(\seqska \) be a sequence of complex \tpqa{matrices}.
 Then \(\seq{U\fpu{j}V}{j}{0}{2\kappa}\) coincides with the \tfpfa{\(\seq{Us_jV}{j}{0}{\kappa}\)}.
\erem

 Using \rrem{F.R.fpkap} and \rlem{F.L.fpkbp}, we can obtain a result analogous to \rlem{K.L.-1^j*k} for the \tfp{} as well:
 
\bleml{F.L.-1^j*f}
 Let \(\seqska \) be a sequence of complex \tpqa{matrices} with \tfpf{} \(\fpseqka\).
 Let the sequence \(\seqr{\kappa}\) be given by \(r_j\defg(-1)^j\su{j}\).
 Denote by \(\rfpseq{2\kappa}\) the \tfupfa{-\obg}{-\ug}{\(\seqr{\kappa}\)}.
 Then \(\rfpu{4k}=\fpu{4k}\) for all \(k\in\NO\) with \(2k\leq\kappa\) and \(\rfpu{4k+1}=\fpu{4k+2}\) for all \(k\in\NO\) with \(2k+1\leq\kappa\).
 Moreover, \(\rfpu{4k+2}=\fpu{4k+1}\) for all \(k\in\NO\) with \(2k+1\leq\kappa\) and \(\rfpu{4k+3}=\fpu{4k+3}\) for all \(k\in\NO\) with \(2k+2\leq\kappa\).
\elem
\bproof
 We consider the nontrivial case \(\kappa\geq1\) only.
 Denote by \(\skpseqka\) the \tkupfa{-\obg}{\(\seqr{\kappa}\)}.
 Because of \rrem{F.R.fpkap}, then \(\rfpu{4k}=\skpu{2k}\) for all \(k\in\NO\) with \(2k\leq\kappa\) and \(\rfpu{4k+1}=\skpu{2k+1}\) for all \(k\in\NO\) with \(2k+1\leq\kappa\).
 Denote by \(\pkpseqka\) the \tkupfa{\obg}{\(\seqska\)}.
 For all \(j\in\mn{0}{\kappa}\), then \((-1)^j\pkpu{j}=\skpu{j}\), according to \rlem{K.L.-1^j*k}, and furthermore \(\pkpu{j}=(-1)^j\fpu{2j}\), by virtue of \rlem{F.L.fpkbp}.
 Consequently, we have \(\rfpu{4k}=\fpu{4k}\) for all \(k\in\NO\) with \(2k\leq\kappa\) and \(\rfpu{4k+1}=\fpu{4k+2}\) for all \(k\in\NO\) with \(2k+1\leq\kappa\). 
 Let the sequence \(\seqt{\kappa-1}\) be given by \(t_j\defeq(-\ug)r_j-r_{j+1}\).
 Denote by \(\tkpseq{\kappa-1}\) the \tkupfa{-\obg}{\(\seqt{\kappa-1}\)}.
 Because of \rrem{F.R.fpkap}, then \(\rfpu{4k+2}=\tkpu{2k}\) for all \(k\in\NO\) with \(2k\leq\kappa-1\) and \(\rfpu{4k+3}=\tkpu{2k+1}\) for all \(k\in\NO\) with \(2k+1\leq\kappa-1\).
 Denote by \(\qkpseq{\kappa-1}\) the \tkupfa{\obg}{\(\seqsa{\kappa-1}\)}.
 Observe that \((-1)^j\sau{j}=t_j\) holds true for all \(j\in\mn{0}{\kappa-1}\).
 According to \rlem{K.L.-1^j*k}, then \((-1)^j\qkpu{j}=\tkpu{j}\) for all \(j\in\mn{0}{\kappa-1}\).
 From \rlem{F.L.fpkbp}, we obtain furthermore \(\qkpu{j}=(-1)^j\fpu{2j+1}\) for all \(j\in\mn{0}{\kappa-1}\).
 Consequently, we have \(\rfpu{4k+2}=\fpu{4k+1}\) for all \(k\in\NO\) with \(2k+1\leq\kappa\) and \(\rfpu{4k+3}=\fpu{4k+3}\) for all \(k\in\NO\) with \(2k+2\leq\kappa\).
\eproof

\bpropl{F.L.A.Fgg.A^*}
 Let \(A\in\Cpq\) and let \(\seqska\in\Fggqka\) with \tfpf{} \(\fpseqka\).
 Let the sequence \(\seqt{\kappa}\) be given by \(\tu{j}\defeq A\su{j}A^\ad\) and denote by \(\seq{\frt_j}{j}{0}{2\kappa}\) the \tfpfa{\(\seqt{\kappa}\)}.
 Then \(\seqt{\kappa}\in\Fgguuuu{p}{\kappa}{\ug}{\obg}\) and, for all \(j\in\mn{0}{2\kappa}\), furthermore \(\frt_j\lgeq A\fpu{j}A^\ad\).
\eprop
\bproof
 In view of \rnota{F.N.sa}, we have obviously \(-\ug\tu{j}+\tu{j+1}=A\sau{j}A^\ad\) and \(\obg\tu{j}-\tu{j+1}=A\sub{j}A^\ad\) for all \(j\in\mn{0}{\kappa-1}\), and furthermore \(-\ug\obg\tu{j}+\rk{\ug+\obg}\tu{j+1}-\tu{j+2}=A\sab{j}A^\ad\) for all \(j\in\mn{0}{\kappa-2}\).
 Because of \eqref{Fgg2n} and \eqref{Fgg2n+1}, hence \(\seqt{\kappa}\in\Fgguuuu{p}{\kappa}{\ug}{\obg}\) follows by virtue of \rprop{H.R.A.Hgg.A^*}.
 Moreover, we can conclude \(\frt_j\lgeq A\fpu{j}A^\ad\) for all \(j\in\mn{0}{2\kappa}\) also from \rprop{H.R.A.Hgg.A^*}, taking into account \rremss{H.R.h2L}{F.R.ABL} and \rdefn{F.D.fpf}.
\eproof

 Using \rprop{H.R.Hgg+Hgg} instead of \rprop{H.R.A.Hgg.A^*}, we get in a similar way:

\bpropl{F.R.Fgg+Fgg}
 Let \(\seqska \) and \(\seqt{\kappa}\) be sequences belonging to \(\Fggqka\) with \tfpf{s} \(\fpseqka\) and \(\seq{\mathfrak{g}_j}{j}{0}{2\kappa}\), \tresp{}
 Let the sequence \(\seq{z_j}{j}{0}{\kappa}\) be given by \(z_{j}\defeq\su{j}+\tu{j}\) and denote by \(\seq{\frz_j}{j}{0}{2\kappa}\) the \tfpfa{\(\seq{z_j}{j}{0}{\kappa}\)}.
 Then \(\seq{z_j}{j}{0}{\kappa}\in\Fggqka\) and \(\frz_j\lgeq\fpu{j}+\mathfrak{g}_j\) for all \(j\in\mn{0}{2\kappa}\).
\eprop

\breml{F.R.l*Fgg}
 Let \(\lambda\in[0,\infp)\) and let \(\seqska\in\Fggqka\).
 Then \(\seq{\lambda s_j}{j}{0}{\kappa}\in\Fggqka\), by virtue of \rprop{F.L.A.Fgg.A^*}.
\erem

 The following result indicates that the \tfpfa{a} sequence \(\seqska\in\Fggqka\) is a convenient tool to handle simultaneously all the four coupled \tnnH{} \tbHms{} associated with the sequence \(\seqska\).
 We write \(\rank A\) for the rank of a complex \tpqa{matrix} \(A\).
 
\bpropl{F.L.FPdet}
 Let \(\seqska\in\Fggqka\).
 Then
\begin{align}
 \rank \Hu{n}&= \sum_{k=0}^{n} \rank \fpu{4k}&&\text{and}&
 \det \Hu{n}&= \prod_{k=0}^{n} \det \fpu{4k}\label{F.L.FPdet.B0}
\end{align}
 for all \(n\in\NO\) with \(2n\leq\kappa\),
\begin{align}
 \rank \Hau{n}&= \sum_{k=0}^{n} \rank \fpu{4k+1},&&&
 \det \Hau{n}&= \prod_{k=0}^{n} \det \fpu{4k+1},\label{F.L.FPdet.B1}\\
 \rank \Hub{n}&= \sum_{k=0}^{n} \rank \fpu{4k+2},&
&\text{and}&
 \det \Hub{n}&= \prod_{k=0}^{n} \det \fpu{4k+2}\label{F.L.FPdet.B2}
\end{align}
 for all \(n\in\NO\) with \(2n+1\leq\kappa\), and
\begin{align}
 \rank \Hab{n}&= \sum_{k=0}^{n} \rank \fpu{4k+3}&
&\text{and}&
 \det \Hab{n}&= \prod_{k=0}^{n} \det \fpu{4k+3}\label{F.L.FPdet.B3}
\end{align}
 for all \(n\in\NO\) with \(2n+2\leq\kappa\).
\eprop
\bproof
 From \rpropss{ab.R1532}{ab.R1535} we obtain \(\seqska\in\Kggqka\).
 In view of \rrem{F.R.fpkap}, the application of~\zitaa{MR3014201}{\clem{4.11}{225}} to the sequence \(\seqska\) yields \eqref{F.L.FPdet.B0} and \eqref{F.L.FPdet.B1}.
 Now assume \(\kappa\geq1\).
 According to \rprop{F.L.sabF}, then \(\seqsb{\kappa-1}\in\Fggqu{\kappa-1}\).
 In particular, \(\seqsb{\kappa-1}\in\Kggqu{\kappa-1}\), by the same reasoning as above.
 Taking into account \rremss{F.R.fpkap}{F.R.c=ab}, the application of~\zitaa{MR3014201}{\clem{4.11}{225}} to the sequence \(\seqsb{\kappa-1}\) yields \eqref{F.L.FPdet.B2} and \eqref{F.L.FPdet.B3}.
\eproof

 Maximization of the determinant in \eqref{F.L.FPdet.B0} was treated by Dette and Studden in~\zita{MR2118663}, using matrix canonical moments introduced in~\zita{MR1883272}.
 We recall their construction and mention the corresponding connections to the present paper briefly in \rsec{F.s2.cm}.

 Now we characterize \tFnnd{ness} analogously to \rpropss{H.P.HggHP}{K.P.KggKP} in terms of the \tfp{s} introduced in \rdefn{F.D.fpf}:
 
\bpropl{F.P.FggFP}
 Let \(\seqska \) be a sequence of complex \tqqa{matrices}.
 Then \(\seqska \in\Fggqka\) if and only if \(\fpu{j}\in\Cggq\) for all \(j\in\mn{0}{2\kappa}\).
\eprop
\bproof
 In view of~\zitaa{arXiv:1701.04246}{\clem{10.10}{31}}, \(\seqska \in\Fggqka\) implies \(\fpu{j}\in\Cggq\) for all \(j\in\mn{0}{2\kappa}\).
 Conversely, assume \(\fpu{j}\in\Cggq\) for all \(j\in\mn{0}{2\kappa}\).
 We are going to show \(\seqska \in\Fggqka\) by mathematical induction.
 Because of \(\fpu{0}=\su{0}=\Hu{0}\), we have \(\seqs{0}\in\Fggqu{0}\).
 Now assume \(\kappa\geq1\) and suppose that there exists an integer \(m\in\mn{0}{\kappa-1}\) such that \(\seqs{m}\in\Fggqu{m}\) holds true.
 According to \rrem{ab.L0911}, then the matrices \(\umg{m}\) and \(\omg{m}\) are both \tH{}.
 By virtue of \rrem{F.R.fp12}, we have furthermore \(\set{\fpu{2m+1},\fpu{2m+2}}=\set{\usc{m+1},\osc{m+1}}=\set{\su{m+1}-\umg{m},\omg{m}-\su{m+1}}\).
 Thus, we conclude \(\su{m+1}\in\CHq\) and \(\umg{m}\lleq\su{m+1}\lleq\omg{m}\), \tie{}\ \(\su{m+1}\in\matint{\umg{m}}{\omg{m}}\).
 \rthm{165.T112} yields then \(\seqs{m+1}\in\Fggqu{m+1}\).
\eproof

\bpropl{F.P.FgFP}
 Let \(\seqska \) be a sequence of complex \tqqa{matrices}.
 Then \(\seqska \in\Fgqu{\kappa}\) if and only if \(\fpu{j}\in\Cgq\) for all \(j\in\mn{0}{2\kappa}\).
\eprop
\bproof
 Combine \rpropss{F.P.FggFP}{F.L.FPdet}.
\eproof

 By virtue of \rlem{F.L.fpkbp} and \rrem{K.R.kpfCqq}, a sequence \(\seqska\) of complex \tpqa{matrices} can be recovered from its even indexed \tfp{s} \(\seq{\fpu{2k}}{k}{0}{\kappa}\), in case that \(\obg\) is known.
 In particular, there are inner dependencies in the complete sequence \(\fpseqka\) of \tfp{s}, which easily can be seen from \rremss{ab.R1420}{F.R.fp12}:

\breml{F.R.f2n-1}
 Let \(\seqska\) be a sequence of complex \tpqa{matrices}.
 For all \(k\in\mn{1}{\kappa}\), then \(\fpu{2k-1}=\dia{k-1}-\fpu{2k}\).
\erem

 Furthermore, in view of \rrem{F.R.fp12}, we get from \rthm{ab.P1422} the following connection between \tfp{s} and \tfd{s} associated to a sequence belonging to \(\Fggqka\):

\breml{F.R.d+=f}
 Let \(\seqska \in\Fggqka\).
 Then \(\dia{0}=\ba\fpu{0}\) and, for all \(k\in\mn{1}{\kappa}\), furthermore \(\dia{k}=\ba\rk{\fpu{2k-1}\ps\fpu{2k}}\) and \(\dia{k}=\ba\rk{\fpu{2k}\ps\fpu{2k-1}}\).
\erem

\breml{F.R.rndf}
 Let \(\seqska\in\Fggqka\).
 By virtue of~\zitaa{arXiv:1701.04246}{\cprop{10.18}{38}} and \rrem{F.R.fp12}, then \(\ran{\dia{0}}=\ran{\fpu{0}}\) and \(\nul{\dia{0}}=\nul{\fpu{0}}\).
 Furthermore, \(\ran{\dia{j}}=\ran{\fpu{2j-1}}\cap\ran{\fpu{2j}}\) and \(\nul{\dia{j}}=\nul{\fpu{2j-1}}+\nul{\fpu{2j}}\) for all \(j\in\mn{1}{\kappa}\) and \(\ran{\dia{j}}=\ran{\fpu{2j+1}}+\ran{\fpu{2j+2}}\) and \(\nul{\dia{j}}=\nul{\fpu{2j+1}}\cap\nul{\fpu{2j+2}}\) for all \(j\in\mn{0}{\kappa-1}\).
\erem

 To single out all sequences \(\fpseqka\) of complex \tqqa{matrices} which indeed occur as \tfp{s} of sequences \(\seqska\in\Fggqka\), we introduce, motivated by \rremss{F.R.f2n-1}{F.R.d+=f}, the following class, using the parallel sum of two matrices, given in \eqref{ps}:

\bnotal{F.N.cs}
 For each \(\eta\in[0,\infp)\), denote by \(\cs{q}{\kappa}{\eta}\) the set of all sequences \(\seq{f_j}{j}{0}{2\kappa}\) of \tnnH{} \tqqa{matrices} satisfying, in the case \(\kappa\geq1\), the equations \(\eta f_{0}=f_{1}+f_{2}\) and \(\eta\rk{f_{2k-1}\ps f_{2k}}=f_{2k+1}+f_{2k+2}\) for all \(k\in\mn{1}{\kappa-1}\).
\enota

 For positive real numbers \(f_j\), the condition \(\eta\rk{f_{2k-1}\ps f_{2k}}=f_{2k+1}+f_{2k+2}\) is fulfilled if and only if the harmonic mean \(h_k\defeq2/\rk{1/f_{2k-1}+1/f_{2k}}\) and the arithmetic mean \(a_{k+1}\defeq\rk{f_{2k+1}+f_{2k+2}}/2\) of consecutive pairs \((f_{2k-1},f_{2k})\) and \((f_{2k+1},f_{2k+2})\) have constant ratio \(a_{k+1}/h_k=\eta/4\).

 The following result provides first information about the importance of the \tfp{ization} of \tFnnd{} sequences.

\bthml{F.T.FggFP}
 The mapping \(\Gamma_{\ug,\obg}\colon\Fggqka\to\csqkad\) given by \(\seqska\mapsto\fpseqka\) is well defined and bijective.
\ethm
\bproof
 Let \(\delta\defeq\bam\).
 First consider \(\seqska \in\Fggqka\) arbitrarily with \tfpf{} \(\fpseqka\).
 In view of \rprop{F.P.FggFP}, we have \(\fpu{j}\in\Cggq\) for all \(j\in\mn{0}{2\kappa}\).
 \rrem{F.R.f2n-1} yields \(\dia{k}=\fpu{2k+1}+\fpu{2k+2}\) for all \(k\in\mn{0}{\kappa-1}\).
 According to \rrem{F.R.d+=f}, we have furthermore \(\dia{0}=\delta\fpu{0}\) and \(\dia{k}=\delta\rk{\fpu{2k-1}\ps\fpu{2k}}\) for all \(k\in\mn{1}{\kappa}\).
 In the case \(\kappa\geq1\), we obtain in particular \(\delta\fpu{0}=\fpu{1}+\fpu{2}\) and \(\delta\rk{\fpu{2k-1}\ps\fpu{2k}}=\fpu{2k+1}+\fpu{2k+2}\) for all \(k\in\mn{1}{\kappa-1}\).
 Hence, \(\fpseqka\) belongs to \(\cs{q}{\kappa}{\delta}\).
 Thus, \(\Gamma_{\ug,\obg}\) is well defined.
 Because of \rrem{F.R.fpfCqq}, the mapping \(\Gamma_{\ug,\obg}\) is furthermore injective.
 To verify the surjectivity of \(\Gamma_{\ug,\obg}\), we consider now an arbitrary sequence \(\seq{f_j}{j}{0}{2\kappa}\) belonging to \(\cs{q}{\kappa}{\delta}\).
 By virtue of \rrem{K.R.kpfCqq}, there exists a sequence \(\seqska\) of complex \tqqa{matrices} with \tkupf{\obg} \(\seq{(-1)^jf_{2j}}{j}{0}{\kappa}\).
 In view of \rlem{F.L.fpkbp}, the \tfpf{} \(\fpseqka\) of \(\seqska\) then fulfills \(\fpu{2k}=f_{2k}\) for all \(k\in\mn{0}{\kappa}\).

 We show \(f_j=\fpu{j}\) for all \(j\in\mn{0}{2\kappa}\) by mathematical induction:
 Obviously, \(f_0=\fpu{0}\).
 Now let \(\kappa\geq1\).
 By construction, we have \(f_0=\fpu{0}\) and \(f_2=\fpu{2}\) and, in view of \(\seq{f_j}{j}{0}{2\kappa}\in\cs{q}{\kappa}{\delta}\) and \eqref{F.G.f012}, \eqref{F.G.d01}, and \rrem{F.R.f2n-1}, then
 \[
  f_1
  =\delta f_0-f_2
  =\delta\fpu{0}-\fpu{2}
  =\delta\su{0}-\fpu{2}
  =\dia{0}-\fpu{2}
  =\fpu{1}.
 \]
 Now assume \(\kappa\geq2\) and suppose that there is an integer \(m\in\mn{1}{\kappa-1}\), such that \(f_j=\fpu{j}\) holds true for all \(j\in\mn{0}{2m}\).
 By virtue of \rrem{F.R.fpftr} and \(\seq{f_j}{j}{0}{2\kappa}\in\cs{q}{\kappa}{\delta}\), the \tfpfa{\(\seqs{m}\)} consists of \tnnH{} matrices.
 Hence, \(\seqs{m}\in\Fggqu{m}\), according to \rprop{F.P.FggFP}.
 \rrem{F.R.d+=f} yields \(\dia{m}=\delta(\fpu{2m-1}\ps\fpu{2m})\).
 We have \(f_{2m+2}=\fpu{2m+2}\) by construction.
 From \(\seq{f_j}{j}{0}{2\kappa}\in\cs{q}{\kappa}{\delta}\) and \rrem{F.R.f2n-1} we get then
 \[
  f_{2m+1}
  =\delta\rk{f_{2m-1}\ps f_{2m}}-f_{2m+2}
  =\delta\rk{\fpu{2m-1}\ps\fpu{2m}}-\fpu{2m+2}
  =\fpu{2m+1}.
 \]
 Thus, \(f_j=\fpu{j}\) is proved for all \(j\in\mn{0}{2\kappa}\) by mathematical induction.
 
 Since \(\seq{f_j}{j}{0}{2\kappa}\) belongs to \(\cs{q}{\kappa}{\delta}\), the \tfpfa{\(\seqska\)} consists of \tnnH{} matrices.
 Hence, \rprop{F.P.FggFP} yields \(\seqska\in\Fggqka\).
 Moreover, \(\seq{f_j}{j}{0}{2\kappa}\) is the image of \(\seqska\) under \(\Gamma_{\ug,\obg}\).
\eproof

 To uncover the above mentioned inner dependencies between the \tfp{s} of a sequence belonging to \(\Fggqka\), we introduce new parameters.
 Observe that the following construction is well-defined due to \rprop{ab.C0929}:

\bdefnl{F.D.cia}
 Let \(\seqska\in\Fggqka\) with \tfpf{} \(\fpseqka\) and \tfdf{} \(\seqdiaka\).
 Then we call \(\seqciaka\) given by \(\cia{0}\defeq\fpu{0}\) and by \(\cia{j}\defeq\rk{\dia{j-1}^\varsqrt}^\mpi\fpu{2j}\rk{\dia{j-1}^\varsqrt}^\mpi\) for each \(j\in\mn{1}{\kappa}\) the \emph{\tfcfa{\(\seqska\)}}.
\edefn

 For the special case that \(\ab=[0,1]\) and \(\su{0}=\Iq\), we obtain, by virtue of \eqref{F.G.f012}, \eqref{F.G.f34}, and \eqref{F.G.d01}, then (\tcf{}~\zitaa{MR1468473}{\cexa{1.3.1}{11}})
 \begin{align*}
 \cia{0}&=\Iq,&
 \cia{1}&=\Iq-\su{1},&
&\text{and}&
 \cia{2}
 &=\ek*{\rk{\su{1}-\su{1}^2}^\varsqrt}^\mpi\rk{\su{2}-\su{1}^2}\ek*{\rk{\su{1}-\su{1}^2}^\varsqrt}^\mpi.
\end{align*}
 Observe that the \tfp{s} \(\cia{j}\) of an \tFnnd{} sequence are connected in the case \(q=1\), \(\ab=[0,1]\), and \(\su{0}=1\) to the classical canonical moments \(p_k\), given in our notation by (\tcf{}~\zitaa{MR1468473}{\ceq{1.3.2}{11}})
\[
 p_k
 =\frac{\su{k}-\umg{k-1}}{\omg{k-1}-\umg{k-1}}
 =\frac{\usc{k}}{\dia{k-1}}, 
\]
 of the corresponding point in the moment space of probability measures on \([0,1]\) via \(p_1=1-\cia{1}\), \(p_2=\cia{2}\), \(p_3=1-\cia{3}\), \(p_4=\cia{4}\), \(p_5=1-\cia{5}\), \ldots , \tie{}, \(\cia{2\ell-1}=q_{2\ell-1}\) and \(\cia{2\ell}=p_{2\ell}\) with the quantities \(q_k=1-p_k\) additionally occurring in the classical framework.

 Now we state some observations concerning the arithmetics of the object introduced in \rdefn{F.D.cia}.
 From~\zitaa{arXiv:1701.04246}{\cprop{7.7(a)}{18}} and \rremss{F.R.diatr}{F.R.fpftr}, we easily get:

\breml{F.R.ciatr}
 Let \(\seqska\in\Fggqka\) with \tfcf{} \(\seqciaka\).
 For each \(k\in\mn{0}{\kappa}\), the matrix \(\cia{k}\) is built from the matrices \(\su{0},\su{1},\dotsc,\su{k}\).
 In particular, for each \(m\in\mn{0}{\kappa}\) the sequence \(\seqs{m}\) belongs to \(\Fggqu{m}\) and has \tfcf{} \(\seqcia{m}\).
\erem

\bleml{F.R.UeU^*}
 Let \(V\in\Cqp\) with \(VV^\ad=\Iq\) and let \(\seqska\in\Fggqka\).
 Then \(\seq{V^\ad\su{j}V}{j}{0}{\kappa}\) belongs to \(\Fgguuuu{p}{\kappa}{\ug}{\obg}\) and has \tfcf{} \(\seq{V^\ad\cia{j}V}{j}{0}{\kappa}\).
\elem
\bproof
 \rprop{F.L.A.Fgg.A^*} yields \(\seq{V^\ad\su{j}V}{j}{0}{\kappa}\in\Fgguuuu{p}{\kappa}{\ug}{\obg}\).
 Denote by \(\seq{\mathfrak{p}_j}{j}{0}{\kappa}\) the \tfcfa{\(\seq{V^\ad\su{j}V}{j}{0}{\kappa}\)}, by \(\seq{\mathfrak{l}_j}{j}{0}{\kappa}\) the \tfdfa{\(\seq{V^\ad\su{j}V}{j}{0}{\kappa}\)}, and by \(\seq{\mathfrak{t}_j}{j}{0}{\kappa}\) the \tfpfa{\(\seq{V^\ad\su{j}V}{j}{0}{\kappa}\)}.
 According to \rrem{F.R.U.f.V}, then \(\mathfrak{t}_j=V^\ad\fpu{j}V\) for all \(j\in\mn{0}{2\kappa}\). In view of \rdefn{F.D.cia}, in particular \(\mathfrak{p}_0=\mathfrak{t}_0=V^\ad\fpu{0}V=V^\ad\cia{0}V\).
 Now assume \(\kappa\geq1\).
 Using \rrem{F.R.U.d.V}, we obtain \(\mathfrak{l}_j=V^\ad\dia{j}V\) for all \(j\in\mn{0}{\kappa}\).
 Consequently, we have \(\rk{\mathfrak{l}_j^\varsqrt}^\mpi=V^\ad\rk{\dia{j}^\varsqrt}^\mpi V\) for all \(j\in\mn{0}{\kappa}\), by virtue of \rremss{A.R.Vsqrt}{A.R.UA+V}.
 In view of \rdefn{F.D.cia} and \(VV^\ad=\Iq\), hence
\begin{multline*}
 \mathfrak{p}_j
 =\rk{\mathfrak{l}_{j-1}^\varsqrt}^\mpi\mathfrak{t}_{2j}\rk{\mathfrak{l}_{j-1}^\varsqrt}^\mpi
 =V^\ad\rk{\dia{j-1}^\varsqrt}^\mpi VV^\ad\fpu{2j}VV^\ad\rk{\dia{j-1}^\varsqrt}^\mpi V\\
 =V^\ad\rk{\dia{j-1}^\varsqrt}^\mpi\fpu{2j}\rk{\dia{j-1}^\varsqrt}^\mpi V
 =V^\ad\cia{j}V
\end{multline*}
 follows for all \(j\in\mn{1}{\kappa}\).
\eproof

 Using \rrem{A.R.l*A}, we obtain from \rremsss{F.R.l*Fgg}{F.R.l*d}{F.R.l*f} in a similar way:

\breml{F.R.l*e}
 Let \(\lambda\in(0,\infp)\) and let \(\seqska\in\Fggqka\) with \tfcf{} \(\seqciaka\).
 Denote by \(\seq{\mathfrak{p}_j}{j}{0}{\kappa}\) the \tfcfa{\(\seq{\lambda\su{j}}{j}{0}{\kappa}\)}.
 Then \(\seq{\lambda s_j}{j}{0}{\kappa}\) belongs to \(\Fggqka\) and furthermore \(\mathfrak{p}_0=\lambda\cia{0}\) and \(\mathfrak{p}_j=\cia{j}\) for all \(j\in\mn{1}{\kappa}\).
\erem

 Extending in the case \(\kappa<\infi\) the sequence \(\seqska\in\Fggqka\) to a sequence \(\seqs{\kappa+1}\in\Fggqu{\kappa+1}\), using, \teg{}, \rcor{ab.R1011}, and applying \rremss{F.R.diatr}{F.R.f2n-1}, we conclude from \rprop{F.L.A.Fgg.A^*} easily:

\bpropl{F.R.A.d.A^*}
 Let \(A\in\Cpq\) and let \(\seqska\in\Fggqka\) with \tfdf{} \(\seqdiaka\).
 Denote by \(\seq{\mathfrak{l}_j}{j}{0}{\kappa}\) the \tfdfa{\(\seq{A\su{j}A^\ad}{j}{0}{\kappa}\)}.
 For all \(j\in\mn{0}{\kappa}\), then \(\mathfrak{l}_j\lgeq A\dia{j}A^\ad\).
\eprop

 Using \rprop{F.R.Fgg+Fgg} instead of \rprop{F.L.A.Fgg.A^*}, we get in a similar way:

\bpropl{F.R.d+d}
 Let \(\seqska \) and \(\seqt{\kappa}\) be sequences belonging to \(\Fggqka\) with \tfdfs{} \(\seqdiaka\) and \(\seq{\mathfrak{v}_j}{j}{0}{2\kappa}\), \tresp{}
 Denote by \(\seq{\mathfrak{l}_j}{j}{0}{\kappa}\) the \tfdfa{\(\seq{\su{j}+\tu{j}}{j}{0}{\kappa}\)}.
 For all \(j\in\mn{0}{\kappa}\), then \(\mathfrak{l}_j\lgeq\dia{j}+\mathfrak{v}_j\).
\eprop

 Our next aim is to single out all sequences \(\seqciaka\) of complex \tqqa{matrices}, which indeed occur as \tfc{s} of sequences \(\seqska\in\Fggqka\).
 With the Euclidean scalar product \(\ipE{\cdot}{\cdot}\colon\x{\Cq}\to\C\) given by \(\ipE{x}{y}\defeq y^\ad x\), which is \(\C\)\nobreakdash-linear in its first argument, the vector space \(\Cq\) over the field \(\C\) becomes a unitary space.
 Let \(\mathcal{U}\) be an arbitrary non-empty subset of \(\Cq\).
 The orthogonal complement \(\mathcal{U}^\orth\defeq\setaa{v\in\Cq}{\ipE{v}{u}=0\text{ for all }u\in\mathcal{U}}\) of \(\mathcal{U}\) is a linear subspace of the unitary space \(\Cq\).
 If \(\mathcal{U}\) is a linear subspace itself, the unitary space \(\Cq\) is the orthogonal sum of \(\mathcal{U}\) and \(\mathcal{U}^\orth\).
 In this case, we write \(\OPu{\mathcal{U}}\) for the transformation matrix corresponding to the orthogonal projection onto \(\mathcal{U}\) with respect to the standard basis of \(\Cq\), \tie{}, \(\OPu{\mathcal{U}}\) is the uniquely determined matrix \(P\in\Cqq\) satisfying \(P^2=P=P^\ad\) and \(\ran{P}=\mathcal{U}\).
 
\bpropl{F.L.cia-fp}
 Let \(\seqska\in\Fggqka\).
 For all \(j\in\mn{1}{\kappa}\), then
\begin{align*}
 \fpu{2j-1}&=\dia{j-1}^\varsqrt \rk{\OPu{\ran{\dia{j-1}}}-\cia{j}}\dia{j-1}^\varsqrt =\dia{j-1}^\varsqrt \rk{\Iq-\cia{j}}\dia{j-1}^\varsqrt ,&
 \fpu{2j}&=\dia{j-1}^\varsqrt \cia{j}\dia{j-1}^\varsqrt ,
\end{align*}
 and furthermore \(\OPu{\ran{\dia{j-1}}}-\cia{j}=\rk{\dia{j-1}^\varsqrt}^\mpi\fpu{2j-1}\rk{\dia{j-1}^\varsqrt}^\mpi\).
\eprop
\bproof
 Let \(j\in\mn{1}{\kappa}\).
 In view of \rprop{ab.C0929}, the matrix \(D\defeq\dia{j-1}\) is \tnnH{}.
 With \(Q\defeq D^\varsqrt \) we have \(\cia{j}=Q^\mpi\fpu{2j}Q^\mpi\) by \rdefn{F.D.cia}.
 From \rremss{F.R.rndf}{A.R.r-sqrt}, we conclude \(\ran{\fpu{2j}}\subseteq\ran{D}=\ran{Q}\) and \(\nul{Q}=\nul{D}\subseteq\nul{\fpu{2j}}\).
 By virtue of \rrem{R.AA+B}, we obtain then
\(
 Q\cia{j}Q
 =QQ^\mpi\fpu{2j}Q^\mpi Q
 =\fpu{2j}
\).
 According to \rrem{A.R.r-sqrt}, we have furthermore
\(
 Q\OPu{\ran{D}}Q
 =Q\OPu{\ran{Q}}Q
 =QQ
 =D
\), whereas \rremss{ab.R1052}{A.R.r-sqrt} yield
\[
 Q^\mpi D Q^\mpi
 =Q^\mpi QQQ^\mpi
 =\OPu{\ran{Q^\ad}}\OPu{\ran{Q}}
 =\OPu{\ran{Q}}\OPu{\ran{Q}}
 =\OPu{\ran{Q}}
 =\OPu{\ran{D}}.
\]
 Using \rrem{F.R.f2n-1}, we get hence
\(
 Q\rk{\OPu{\ran{D}}-\cia{j}}Q
 =Q\rk{\Iq-\cia{j}}Q
 =D-\fpu{2j}
 =\fpu{2j-1}
\)
 and
\(
 Q^\mpi\fpu{2j-1}Q^\mpi
 =Q^\mpi\rk{D-\fpu{2j}}Q^\mpi
 =Q^\mpi DQ^\mpi-Q^\mpi\fpu{2j}Q^\mpi
 =\OPu{\ran{D}}-\cia{j}
\).
\eproof

 Using \rrem{A.R.XAX<=XBX} we see that the following construction is well defined.

\bnotal{F.N.01seq}
 For each \(\eta\in[0,\infp)\), let \(\es{q}{\kappa}{\eta}\) be the set of all sequences \(\seq{e_k}{k}{0}{\kappa}\) from \(\Cggq\) which, in the case \(\kappa\geq1\), fulfill \(e_k\lleq\OPu{\ran{d_{k-1}}}\) for all \(k\in\mn{1}{\kappa}\), where the sequence \(\seq{d_k}{k}{0}{\kappa}\) is recursively given by \(d_0\defeq\eta e_0\) and
\[
 d_k
 \defeq\eta d_{k-1}^\varsqrt  e_k^\varsqrt (\OPu{\ran{d_{k-1}}}-e_k) e_k^\varsqrt  d_{k-1}^\varsqrt .
\]
\enota

 Consider the scalar case \(q=1\) and an \(\eta\in(0,\infp)\):
 The set \(\es{1}{\infi}{\eta}\) consists of all sequences \(\seq{e_k}{k}{0}{\infi}\) of real numbers, having one of the following two forms:
\begin{itemize}
 \item \(0=e_0=e_1=e_2=\dotsb\)
 \item \(e_0\in(0,\infp)\) and there exists an integer \(k\in\N\) such that \(e_1,e_2,\dotsc,e_{k-1}\in(0,1)\) if \(k\geq2\), and, furthermore, \(e_k\in\set{0,1}\) and \(0=e_{k+1}=e_{k+2}=\dotsb\)
\end{itemize}

\bpropl{F.L.cs-es}
 Let \(\eta\in[0,\infp)\) and let \(\Delta_\eta\colon\cs{q}{\kappa}{\eta}\to\es{q}{\kappa}{\eta}\) be defined by \(\seq{f_j}{j}{0}{2\kappa}\mapsto\seq{e_k}{k}{0}{\kappa}\), where the sequence \(\seq{e_k}{k}{0}{\kappa}\) is given via
 \[
  e_k
  \defeq
  \begin{cases}
   f_0\tincase{k=0}\\
   \rk{h_{k-1}^\varsqrt}^\mpi f_{2k}\rk{h_{k-1}^\varsqrt}^\mpi\tincase{k\geq1}
  \end{cases}
 \]
 where the sequence \(\seq{h_k}{k}{0}{\kappa}\) is given by
 \beql{F.L.cs-es.A}
  h_k
  \defeq
  \begin{cases}
   \eta f_0\tincase{k=0}\\
   \eta(f_{2k-1}\ps f_{2k})\tincase{k\geq1}
  \end{cases}.
 \eeq
 Then \(\Delta_\eta\) is well defined and bijective.
 Furthermore, \(\seq{h_k}{k}{0}{\kappa}\) coincides with the sequence \(\seq{d_k}{k}{0}{\kappa}\) associated to \(\seq{e_k}{k}{0}{\kappa}\) in \rnota{F.N.01seq}.
\eprop
\bproof
 We first derive some facts, which will be used several times in the proof:
 Consider a matrix \(D\in\Cggq\).
 Let \(P\defeq\OPu{\ran{D}}\) and \(Q\defeq D^\varsqrt \).
 Because of \rrem{A.R.r-sqrt}, we have \(\OPu{\ran{Q}}=P\).
 Hence, \rrem{ab.R1052} yields
 \begin{align}\label{F.L.cs-es.2}
  QQ^\mpi&=P=DD^\mpi&
 &\text{and}&
  Q^\mpi Q &=P.
 \end{align}
 In view of \rrem{A.R.A+>}, we conclude then
 \begin{align}\label{F.L.cs-es.34}
  P&=P^2=Q\rk{Q^\mpi}^2Q=QD^\mpi Q&
&\text{and}&
  P&=P^2=Q^\mpi Q^2Q^\mpi=Q^\mpi DQ^\mpi.
 \end{align}
 Additionally, we consider now two matrices \(E,F\in\Cqq\).
 First suppose that
 \beql{F.L.cs-es.1}
  E
  =Q^\mpi FQ^\mpi
 \eeq
 and
 \beql{F.L.cs-es.9}
  \Oqq
  \lleq F
  \lleq D.
 \eeq
 hold true.
 According to \rrem{A.R.A++*}, we have \(\rk{Q^\mpi}^\ad=Q^\mpi\).
 In combination with \eqref{F.L.cs-es.1}, \eqref{F.L.cs-es.9}, \rrem{A.R.XAX<=XBX}, and \eqref{F.L.cs-es.34}, we get then
 \beql{F.L.cs-es.10}
  E 
  =Q^\mpi FQ^\mpi
  \lleq Q^\mpi DQ^\mpi
  =P.
 \eeq
 Using \rlem{A.R.rA<rB} and \rremp{R.AA+B}{R.AA+B.a}, we conclude from \eqref{F.L.cs-es.9} easily
 \begin{align}\label{F.L.cs-es.6}
  PF&=F&
&\text{and}&
  FP&=\rk{PF}^\ad=F.
 \end{align}
 By virtue of \eqref{F.L.cs-es.1}, \eqref{F.L.cs-es.2} and \eqref{F.L.cs-es.6}, we obtain
 \beql{F.L.cs-es.7}
  QEQ 
  =QQ ^\mpi FQ^\mpi Q 
  =PFP 
  =F.
 \eeq
 In view of \eqref{ps}, we get from \eqref{F.L.cs-es.2}, \eqref{F.L.cs-es.6}, \eqref{F.L.cs-es.7}, and \eqref{F.L.cs-es.34} then
\[\begin{split}
  \rk{D-F}\ps F 
  &=\rk{D-F}D^\mpi F 
  =DD^\mpi F-FD ^\mpi F
  =PF-FD ^\mpi F\\
  &=F-FD ^\mpi F
  =QEQ-QEQD^\mpi QEQ 
  =Q\rk{E-EPE}Q.
\end{split}\]
 Because of \eqref{F.L.cs-es.1} and \eqref{F.L.cs-es.2}, we have \(PE=E\).
 In combination with \rremsps{R.AA+B}{R.AA+B.a}{A.R.r-sqrt}, immediately \(P E^\varsqrt = E^\varsqrt \) follows.
 Thus, we have
 \[
  E-EPE 
  = E^\varsqrt P E^\varsqrt -E^2
  = E^\varsqrt (P-E) E^\varsqrt 
 \]
 and, consequently,
 \beql{F.L.cs-es.11}
  \rk{D-F}\ps F 
  =Q E^\varsqrt (P-E) E^\varsqrt Q.
 \eeq
 Conversely, suppose that
 \beql{F.L.cs-es.1b}
  F
  =QEQ
 \eeq
 and
 \beql{F.L.cs-es.9b}
  \Oqq
  \lleq E
  \lleq P
 \eeq
 are fulfilled.
 Using \rlem{A.R.rA<rB} and \rremp{R.AA+B}{R.AA+B.a}, from \eqref{F.L.cs-es.9b} we conclude
 \begin{align}\label{F.L.cs-es.6b}
  PE&=E&
&\text{and}&
  EP&=\rk{PE}^\ad=E.
 \end{align}
 By virtue of \eqref{F.L.cs-es.1b}, \eqref{F.L.cs-es.2}, and \eqref{F.L.cs-es.6b}, we infer
 \beql{F.L.cs-es.7b}
  Q^\mpi FQ^\mpi
  =Q^\mpi QEQQ^\mpi
  =PEP 
  =E.
 \eeq
 In view of \eqref{ps}, \eqref{F.L.cs-es.2}, \eqref{F.L.cs-es.1b}, and \eqref{F.L.cs-es.34}, we get
\[\begin{split}
  \rk{D-F}\ps F 
  &=\rk{D-F}D^\mpi F 
  =DD^\mpi F-FD ^\mpi F\\
  &=QQ^\mpi QEQ-QEQD^\mpi QEQ 
  =Q\rk{E-EPE}Q.
\end{split}\]
 From the first equation in \eqref{F.L.cs-es.6b} the identity \(E-EPE= E^\varsqrt (P-E) E^\varsqrt \) follows in the same way as above.
 Hence, we obtain again \eqref{F.L.cs-es.11}. 
 We now start proving the assertion of the lemma, only considering the non-trivial case \(\kappa\geq1\):
 First let \(\seq{f_j}{j}{0}{2\kappa}\in\cs{q}{\kappa}{\eta}\) and let \(\seq{e_k}{k}{0}{\kappa}\) be its image under \(\Delta_\eta\).
 Then \(f_j\in\Cggq\) for all \(j\in\mn{0}{2\kappa}\).
 In view of \rthm{ab.T0909}, hence \(h_k\in\Cggq\) for all \(k\in\mn{0}{\kappa}\).
 Using \rremss{A.R.A++*}{A.R.XAX<=XBX}, we then obtain \(e_k\in\Cggq\) for all \(k\in\mn{0}{\kappa}\).
 Now consider an arbitrary integer \(k\in\mn{1}{\kappa}\).
 Observe that the matrices \(E\defeq e_k\) and \(F\defeq f_{2k}\) fulfill \eqref{F.L.cs-es.1} with \(D\defeq h_{k-1}\).
 Because of \(\seq{f_j}{j}{0}{2\kappa}\in\cs{q}{\kappa}{\eta}\) and \eqref{F.L.cs-es.A}, we have
 \beql{F.L.cs-es.8}
  f_{2k-1}+f_{2k}
  =h_{k-1}.
 \eeq
 In particular, \eqref{F.L.cs-es.9} holds true.
 Thus, from \eqref{F.L.cs-es.10} (\tresp{}, \eqref{F.L.cs-es.7}) we get
 \beql{F.L.cs-es.10a}
  e_k
  \lleq\OPu{\ran{h_{k-1}}}
 \eeq
 and
 \beql{F.L.cs-es.7a}
   h_{k-1}^\varsqrt e_k h_{k-1}^\varsqrt 
  =f_{2k}.
 \eeq
 In view of \eqref{F.L.cs-es.8} and \eqref{F.L.cs-es.11}, we conclude
\[\begin{split}
  h_{k}
  &=\eta\rk{f_{2k-1}\ps f_{2k}}
  =\eta\ek*{\rk{D-F}\ps F}\\
  &=\eta Q E^\varsqrt (P-E) E^\varsqrt Q
  =\eta  h_{k-1}^\varsqrt  e_{k}^\varsqrt (\OPu{\ran{h_{k-1}}}-e_{k}) e_{k}^\varsqrt  h_{k-1}^\varsqrt .
\end{split}\]
 Since \(h_0=\eta e_0\) also holds true, hence \(\seq{h_k}{k}{0}{\kappa}\) coincides with the sequence \(\seq{d_k}{k}{0}{\kappa}\) associated to \(\seq{e_k}{k}{0}{\kappa}\) in \rnota{F.N.01seq}.
 By virtue of \eqref{F.L.cs-es.10a}, thus \(\seq{e_k}{k}{0}{\kappa}\) belongs to \(\es{q}{\kappa}{\eta}\).
 Consequently, \(\Delta_\eta\) is well defined.
 Observe that, for each \(k\in\mn{2}{\kappa}\), the matrix \(h_{k-1}\) is built from \(f_{2k-3}\) and \(f_{2k-2}\).
 Hence, using \eqref{F.L.cs-es.7a} and \eqref{F.L.cs-es.8}, the whole sequence \(\seq{f_j}{j}{0}{2\kappa}\) can be reconstructed from the sequence \(\seq{e_k}{k}{0}{\kappa}\).
 This shows that \(\Delta_\eta\) is injective.
 To prove the surjectivity of \(\Delta_\eta\), we consider an arbitrary sequence \(\seq{e_k}{k}{0}{\kappa}\in\es{q}{\kappa}{\eta}\).
 Denote by \(\seq{d_k}{k}{0}{\kappa}\) the sequence associated to \(\seq{e_k}{k}{0}{\kappa}\) in \rnota{F.N.01seq}.
 Let \(\seq{f_j}{j}{0}{2\kappa}\) be given by \(f_0\defeq e_0\) and by \(f_{2k}\defeq d_{k-1}^\varsqrt e_k d_{k-1}^\varsqrt \) and \(f_{2k-1}\defeq d_{k-1}^\varsqrt \rk{\OPu{\ran{d_{k-1}}}-e_k} d_{k-1}^\varsqrt \) for \(k\in\mn{1}{\kappa}\).
 Observe that \(e_k\in\Cggq\) for all \(k\in\mn{0}{\kappa}\) and \(e_k\lleq\OPu{\ran{d_{k-1}}}\) for all \(k\in\mn{1}{\kappa}\).
 Hence, because of \rrem{A.R.XAX<=XBX}, we have then \(f_{j}\in\Cggq\) for all \(j\in\mn{0}{2\kappa}\).
 Now consider again an arbitrary \(k\in\mn{1}{\kappa}\).
 The matrix \(D\defeq d_{k-1}\) belongs to \(\Cggq\).
 Using \eqref{F.L.cs-es.2}, thus
 \beql{F.L.cs-es.8b}
  f_{2k-1}+f_{2k}
  = d_{k-1}^\varsqrt \OPu{\ran{d_{k-1}}} d_{k-1}^\varsqrt 
  =QPQ
  =QQQ^\mpi Q
  =Q^2
  =D
 \eeq
 follows.
 The matrices \(E\defeq e_k\) and \(F\defeq f_{2k}\) fulfill \eqref{F.L.cs-es.1b} and \eqref{F.L.cs-es.9b}.
 In view of \eqref{F.L.cs-es.8b} and \eqref{F.L.cs-es.11}, hence
\beql{F.L.cs-es.12b}\begin{split}
 \eta\rk{f_{2k-1}\ps f_{2k}}
 &=\eta\ek*{\rk{D-F}\ps F}
 =\eta Q E^\varsqrt (P-E) E^\varsqrt Q\\
 &=\eta  d_{k-1}^\varsqrt  e_{k}^\varsqrt (\OPu{\ran{d_{k-1}}}-e_{k}) e_{k}^\varsqrt  d_{k-1}^\varsqrt 
 =d_k.
\end{split}\eeq
 Observe that, in the case \(k\leq\kappa-1\), we get
 \(
  f_{2k+1}+f_{2k+2}
  = d_k^\varsqrt \OPu{\ran{d_k}} d_k^\varsqrt 
  =d_k
 \)
 in the same way as above.
 Thus, we have \(\eta\rk{f_{2k-1}\ps f_{2k}}=f_{2k+1}+f_{2k+2}\) for all \(k\in\mn{1}{\kappa-1}\).
 Since \(d_0=\eta f_0\) holds true, additionally \(\eta f_0=f_1+f_2\) follows from \eqref{F.L.cs-es.8b} if \(k=1\).
 Consequently, \(\seq{f_j}{j}{0}{2\kappa}\) belongs to \(\cs{q}{\kappa}{\eta}\).
 In view of \(d_0=\eta f_0\) and \eqref{F.L.cs-es.12b}, the sequence \(\seq{d_k}{k}{0}{\kappa}\) coincides with the sequence \(\seq{h_k}{k}{0}{\kappa}\) given via \eqref{F.L.cs-es.A}.
 By virtue of \eqref{F.L.cs-es.7b}, we obtain then
 \(
  \rk{h_{k-1}^\varsqrt}^\mpi f_{2k}\rk{h_{k-1}^\varsqrt}^\mpi
  =Q^\mpi FQ^\mpi
  =E
  =e_k
 \).
 Hence, the sequence \(\seq{e_k}{k}{0}{\kappa}\) is the image of \(\seq{f_j}{j}{0}{2\kappa}\) under \(\Delta_\eta\).
 This shows that \(\Delta_\eta\) is surjective.
\eproof

\bthml{F.T.Fggcia}
 The mapping \(\Sigma_{\ug,\obg}\colon\Fggqka\to\esqkad\) given by \(\seqska\mapsto\seqciaka\) is well defined and bijective.
\ethm
\bproof
 Let \(\Gamma_{\ug,\obg}\colon\Fggqka\to\csqkad\) be defined by \(\seqska\mapsto\fpseqka\), where \(\fpseqka\) is the \tfpfa{\(\seqska\)}.
 According to \rthm{F.T.FggFP}, then \(\Gamma_{\ug,\obg}\) is well defined and bijective.
 By virtue of \rprop{F.L.cs-es}, the mapping \(\Delta_\delta\colon\csqkad\to\esqkad\) given there is well defined and bijective as well.
 From \rdefn{F.D.cia} and \rrem{F.R.d+=f} we see that \(\Sigma_{\ug,\obg}=\Delta_\delta\circ\Gamma_{\ug,\obg}\).
\eproof

 Now we are going to describe some connections between \tfc{s} and \tfd{s}:

\bnotal{F.N.Pd}
 Let \(\seqska\) be a sequence of complex \tpqa{matrices} with \tfdf{} \(\seqdiaka\).
 Then for all \(j\in\mn{0}{\kappa}\), let \(\Pd{j}\defeq\OPu{\ran{\dia{j}}}\) be the matrix corresponding to the orthogonal projection onto \(\ran{\dia{j}}\).
\enota

\bpropl{F.P.ed}
 Let \(\seqska\in\Fggqka\).
 Then \(\cia{0}\lgeq\Oqq\) and \(\cia{j}\in\matint{\Oqq}{\Pd{j-1}}\) for all \(j\in\mn{1}{\kappa}\).
 Furthermore, \(\dia{0}=\ba\cia{0}\) and \(\dia{j}=\ba\dia{j-1}^\varsqrt  \cia{j}^\varsqrt (\Pd{j-1}-\cia{j}) \cia{j}^\varsqrt \dia{j-1}^\varsqrt \) for all \(j\in\mn{1}{\kappa}\).
\eprop
\bproof
 By virtue of \rprop{F.L.cs-es}, the mapping \(\Delta_\delta\colon\csqkad\to\esqkad\) given there is well defined and bijective.
 Denote by \(\fpseqka\) the \tfpfa{\(\seqska\)}.
 According to \rthm{F.T.FggFP}, then \(\fpseqka\in\csqkad\).
 From \rdefn{F.D.cia} and \rrem{F.R.d+=f} we conclude in view of \rprop{F.L.cs-es} that \(\seqciaka\) is the image of \(\fpseqka\) under \(\Delta_\delta\) and that the sequence \(\seqdiaka\) is of the asserted form given in \rnota{F.N.01seq}.
 Consequently, the sequence \(\seqciaka\) belongs to \(\esqkad\).
 In particular, \(\cia{0}\lgeq\Oqq\) and \(\Oqq\lleq\cia{j}\lleq\Pd{j-1}\) for all \(j\in\mn{1}{\kappa}\). 
\eproof

\bpropl{F.L.1-e}
 Let \(\seqska\in\Fggqka\).
 For all \(j\in\mn{1}{\kappa}\), then
\beql{F.L.1-e.A}\begin{split}
 \dia{j}
 &=\ba\dia{j-1}^\varsqrt \cia{j}\rk{\Iq-\cia{j}}\dia{j-1}^\varsqrt \\
 &=\frac{\bam}{4}\dia{j-1}-\ba\dia{j-1}^\varsqrt \rk*{\frac{1}{2}\Pd{j-1}-\cia{j}}^\ad\rk*{\frac{1}{2}\Pd{j-1}-\cia{j}}\dia{j-1}^\varsqrt 
\end{split}\eeq
 and, in particular, \(\det\rk{\dia{j}}=\ba^q\det\rk{\dia{j-1}}\det\rk{\cia{j}}\det\rk{\Iq-\cia{j}}\).
\eprop
\bproof
 Assume \(\kappa\geq1\) and let \(j\in\mn{1}{\kappa}\).
 In view of \rremss{A.R.r-sqrt}{ab.R1842} and \rdefn{F.D.cia}, we have \(\ran{ \cia{j}^\varsqrt }=\ran{\cia{j}}\subseteq\ran{\rk{\dia{j-1}^\varsqrt}^\mpi}=\ran{\dia{j-1}^\varsqrt }=\ran{\dia{j-1}}\), implying \(\Pd{j-1} \cia{j}^\varsqrt = \cia{j}^\varsqrt \) and \(\Pd{j-1}\dia{j-1}^\varsqrt=\dia{j-1}^\varsqrt\).
 In particular, \(\rk{\Pd{j-1}-\cia{j}}\cia{j}^\varsqrt=\rk{\Iq-\cia{j}}\cia{j}^\varsqrt\) and \(\rk{\frac{1}{2}\Pd{j-1}-\cia{j}}\dia{j-1}^\varsqrt=\rk{\frac{1}{2}\Iq-\cia{j}}\dia{j-1}^\varsqrt\).
 According to \rprop{F.P.ed}, we obtain with \(\delta\defeq\bam\) then
\[
 \dia{j}
 =\delta\dia{j-1}^\varsqrt  \cia{j}^\varsqrt (\Iq-\cia{j}) \cia{j}^\varsqrt \dia{j-1}^\varsqrt
 =\delta\dia{j-1}^\varsqrt \cia{j}\rk{\Iq-\cia{j}}\dia{j-1}^\varsqrt.
\]
 and conclude from this furthermore
\begin{multline*}
 \delta\dia{j-1}^\varsqrt \rk*{\frac{1}{2}\Pd{j-1}-\cia{j}}^\ad\rk*{\frac{1}{2}\Pd{j-1}-\cia{j}}\dia{j-1}^\varsqrt
 =\delta\dia{j-1}^\varsqrt\rk*{\frac{1}{2}\Iq-\cia{j}}^2\dia{j-1}^\varsqrt\\
 =\delta\dia{j-1}^\varsqrt\rk*{\frac{1}{4}\Iq-\cia{j}+\cia{j}^2}\dia{j-1}^\varsqrt
 =\frac{\delta}{4}\dia{j-1}-\dia{j}.\qedhere
\end{multline*}
\eproof

 Observe that the representation \eqref{F.L.1-e.A} yields the inequality in \rprop{ab.P1057}.
 The distinguished points \(\umg{m}\), \(\omg{m}\), and \(\mi{m}\), given in \rdefnss{F.N.umom}{F.D.mi}, of the extension interval \(\matint{\umg{m}}{\omg{m}}\) in \rthm{165.T112} stand in a direct correspondence to the from the geometric view analogous points in the interval \(\matint{\Oqq}{\Pd{m}}\) in \rprop{F.P.ed}:

\bthml{F.C.ecase}
 Let \(\seqska\in\Fggqka\), assume \(\kappa\geq1\), and let \(k\in\mn{1}{\kappa}\).
 Then:
\benui
 \il{F.C.ecase.a} If \(k\) is even, then \(\cia{k}=\Oqq\) (\tresp{}, \(\cia{k}=\Pd{k-1}\)) if and only if \(\su{k}=\umg{k-1}\) (\tresp{}, \(\su{k}=\omg{k-1}\)).
 In this case, \(\dia{k}=\Oqq\).
 \il{F.C.ecase.b} If \(k\) is odd, then \(\cia{k}=\Oqq\) (\tresp{}, \(\cia{k}=\Pd{k-1}\)) if and only if \(\su{k}=\omg{k-1}\) (\tresp{}, \(\su{k}=\umg{k-1}\)).
 In this case, \(\dia{k}=\Oqq\).
 \il{F.C.ecase.c} \(\cia{k}=\frac{1}{2}\Pd{k-1}\) if and only if \(\su{k}=\mi{k-1}\).
 In this case, \(\dia{k}=\frac{\bam}{4}\dia{k-1}\).
 \il{F.C.ecase.d} \(\dia{k}=\Oqq\) if and only if \(\cia{k}^2=\cia{k}\).
 In this case, \(\dia{j}=\Oqq\), \(\cia{j}=\Oqq\), and \(\su{j}=\umg{j-1}=\mi{j-1}=\omg{j-1}\) for all \(j\in\mn{k+1}{\kappa}\).
 \il{F.C.ecase.e} \(\det\dia{k}\neq0\) if and only if \(\det\cia{j}\neq0\) for all \(j\in\mn{0}{k}\) and \(\det\rk{\Iq-\cia{j}}\neq0\) for all \(j\in\mn{1}{k}\).
 In this case, \(\det\dia{j}\neq0\) for all \(j\in\mn{0}{k}\) and furthermore \(\cia{0}\lgs\Oqq\) and \(\Oqq\lls\cia{j}\lls\Iq\) for all \(j\in\mn{1}{k}\).
\eenui
\ethm
\bproof
 \eqref{F.C.ecase.a}--\eqref{F.C.ecase.b} In view of \rdefnss{F.D.fpf}{F.N.AB}, the equivalences are direct consequences of \rprop{F.L.cia-fp} and \rdefn{F.D.cia}.
 From \rprop{F.P.ed}, the conclusion for \(\dia{k}\) follows. 

 \eqref{F.C.ecase.c} First assume \(\cia{k}=\frac{1}{2}\Pd{k-1}\).
 From \rprop{F.L.cia-fp} we conclude then \(\fpu{2k-1}=\fpu{2k}\).
 Hence, \rrem{F.R.fp12} yields \(\usc{k}=\osc{k}\).
 In view of \rdefn{F.N.AB}, this means \(\su{k}-\umg{k-1}=\omg{k-1}-\su{k}\), implying \(\su{k}=\frac{1}{2}\rk{\umg{k-1}+\omg{k-1}}=\mi{k-1}\) by virtue of \rdefn{F.D.mi}.
 According to \rprop{ab.P1057}, then \(\dia{k}=\frac{\bam}{4}\dia{k-1}\).
 Conversely, assume \(\su{k}=\mi{k-1}\).
 By the same reasoning as above, we conclude then \(\fpu{2k-1}=\fpu{2k}\).
 Consequently, from \rdefn{F.D.cia} and \rprop{F.L.cia-fp}, we obtain \(\cia{k}=\Pd{k-1}-\cia{k}\), implying \(\cia{k}=\frac{1}{2}\Pd{k-1}\).

 \eqref{F.C.ecase.d} In view of \rdefn{F.D.cia} and \rrem{ab.R1842}, we have
\(
 \ran{\cia{k}}
 \subseteq\ran{\rk{\dia{k-1}^\varsqrt}^\mpi}
 =\ran{\dia{k-1}^\varsqrt }
\)
 and hence \(\OPu{\ran{\dia{k-1}^\varsqrt }}\cia{k}=\cia{k}\).
 Since, by virtue of \rprop{F.P.ed}, the matrix \(\cia{k}\) is \tH{}, we can thus conclude \(\cia{k}\dia{k-1}^\varsqrt \rk{\dia{k-1}^\varsqrt}^\mpi=\cia{k}\OPu{\ran{\dia{k-1}^\varsqrt }}=\rk{\OPu{\ran{\dia{k-1}^\varsqrt }}\cia{k}}^\ad=\cia{k}\), using \rrem{ab.R1052}.
 Taking additionally into account \rrem{A.R.A+>}, furthermore \(\rk{\dia{k-1}^\varsqrt}^\mpi\dia{k-1}^\varsqrt \cia{k}=\rk{\cia{k}\dia{k-1}^\varsqrt \rk{\dia{k-1}^\varsqrt}^\mpi}^\ad=\cia{k}\) follows.
 According to \rprop{F.L.1-e}, we have \(\dia{k}=\ba\dia{k-1}^\varsqrt \rk{\cia{k}-\cia{k}^2}\dia{k-1}^\varsqrt \) and obtain then \(\rk{\dia{k-1}^\varsqrt}^\mpi\dia{k}\rk{\dia{k-1}^\varsqrt}^\mpi=\ba\rk{\cia{k}-\cia{k}^2}\).
 In view of \(\bam>0\), thus \(\cia{k}^2=\cia{k}\) is equivalent to \(\dia{k}=\Oqq\).
 Now assume \(\dia{k}=\Oqq\).
 By virtue of \rprop{ab.C1101}, then \(\dia{\ell}=\Oqq\) and therefore \(\Pd{\ell}=\Oqq\) for all \(\ell\in\mn{k}{\kappa}\).
 In view of \rdefn{F.D.cia}, thus \(\cia{j}=\Oqq\) follows for all \(j\in\mn{k+1}{\kappa}\).
 For all \(j\in\mn{k+1}{\kappa}\), in particular, \(\cia{j}=\Oqq=\Pd{j-1}=\frac{1}{2}\Pd{j-1}\), which because of \rpartstos{F.C.ecase.a}{F.C.ecase.c} implies \(\su{j}=\umg{j-1}=\omg{j-1}=\mi{j-1}\).

 \eqref{F.C.ecase.e} First assume \(\det\dia{k}\neq0\).
 Using \rprop{ab.C1101}, we then can conclude \(\det\dia{\ell}\neq0\) and hence \(\Pd{\ell}=\Iq\) for all \(j\in\mn{0}{k}\).
 This implies \(\det\cia{0}\neq0\) and \(\Oqq\lleq\cia{j}\lleq\Iq\) for all \(j\in\mn{1}{k}\), by virtue of \rprop{F.P.ed}, and furthermore \(\det\cia{j}\neq0\) and \(\det\rk{\Iq-\cia{j}}\neq0\) for all \(j\in\mn{1}{k}\), according to \rprop{F.L.1-e}.
 Since we additionally have \(\cia{0}\lgeq\Oqq\) from \rprop{F.P.ed}, consequently \(\cia{0}\lgs\Oqq\) and \(\Oqq\lls\cia{j}\lls\Iq\) for all \(j\in\mn{1}{k}\) follow.
 Conversely, suppose that \(\det\cia{j}\neq0\) for all \(j\in\mn{0}{k}\) and \(\det\rk{\Iq-\cia{j}}\neq0\) for all \(j\in\mn{1}{k}\) hold true.
 In particular, \(\det\dia{0}=\ba^q\det\cia{0}\neq0\), in view of \rprop{F.P.ed}.
 Using \rprop{F.L.1-e}, we then get \(\det\dia{1}\neq0, \dotsc,\det\dia{k}\neq0\) in a successive way.
\eproof

 In view of \rrem{A.R.P<P}, we obtain  from \rprop{ab.C1101} immediately:
 
\breml{F.R.P<P}
 Let \(\seqska\in\Fggqka\).
 For all \(j\in\mn{1}{\kappa}\), then \(\Pd{j}\lleq\Pd{j-1}\).
\erem

 Hence, the matricial intervals \(\matint{\Oqq}{\Pd{j-1}}\) in \rprop{F.P.ed} form a nested decreasing sequence.
 By virtue of \rprop{ab.C0929}, the situation simplifies in the \tFpd{} case to \(\Pd{j-1}=\Iq\):

\bpropl{F.P.Fg-cia}
 Let \(\seqska \) be a sequence of complex \tqqa{matrices}.
 Then \(\seqska \in\Fgqka\) if and only if \(\cia{0}\in\Cgq\) and \(\cia{j}\in\matinto{\Oqq}{\Iq}\) for all \(j\in\mn{1}{\kappa}\).
\eprop
\bproof
 In view of \(\cia{0}=\fpu{0}=\usc{0}=\su{0}\), the case \(\kappa=0\) is readily checked.
 Now assume \(\kappa\geq1\).
 
 First suppose that \(\seqska\) belongs to \(\Fgqka\).
 According to \rprop{ab.C0929}, then \(\dia{j}\in\Cgq\) for all \(j\in\mn{0}{\kappa}\).
 Consequently, \rthmp{F.C.ecase}{F.C.ecase.e} yields \(\cia{0}\lgs\Oqq\) and \(\Oqq\lls\cia{j}\lls\Iq\) for all \(j\in\mn{1}{\kappa}\).
 
 Conversely, assume \(\cia{0}\in\Cgq\) and \(\cia{j}\in\matinto{\Oqq}{\Iq}\) for all \(j\in\mn{1}{\kappa}\).
 Let \(\eta\defeq\ba\) and let the sequence \(\seq{e_k}{k}{0}{\kappa}\) be given by \(e_k\defeq\cia{k}\).
 Then it is readily checked that the sequence \(\seq{d_k}{k}{0}{\kappa}\) associated to \(\seq{e_k}{k}{0}{\kappa}\) via \rnota{F.N.01seq} is well defined, fulfilling \(d_k\in\Cgq\) and \(\OPu{\ran{d_k}}=\Iq\) for all \(k\in\mn{0}{\kappa}\).
 Consequently, the sequence \(\seq{e_k}{k}{0}{\kappa}\) belongs to \(\es{q}{\kappa}{\eta}\), \tie{}, \(\seqciaka\in\esqkad\).
 According to \rthm{F.T.Fggcia}, thus \(\seqska\in\Fggqka\) follows.
 Hence, we can apply \rthmp{F.C.ecase}{F.C.ecase.e} to obtain \(\det\dia{j}\neq0\) for all \(j\in\mn{0}{\kappa}\).
 From \rrem{F.R.rndf} we conclude then \(\det\fpu{j}\neq0\) for all \(j\in\mn{0}{2\kappa}\).
 Since \rprop{F.P.FggFP} yields \(\fpu{j}\in\Cggq\) for all \(j\in\mn{0}{2\kappa}\), we get \(\fpu{j}\in\Cgq\) for all \(j\in\mn{0}{2\kappa}\).
 Using \rprop{F.P.FgFP}, we infer \(\seqska \in\Fgqka\).
\eproof
  
 In the remaining part of this section, we characterize two special properties for sequences \(\seqska\in\Fggqka\) already considered in~\zita{arXiv:1701.04246} in terms of their \tfc{s} \(\seqciaka\):
 
\bdefnnl{\tcf{}~\zitaa{arXiv:1701.04246}{\cdefnssp{10.24}{10.25}{40}}}{F.D.abCD}
 Let \(\seqska\) be a sequence of complex \tpqa{matrices} with \tfdf{} \(\seqdiaka\) and let \(k\in\mn{0}{\kappa}\).
 Then \(\seqska\) is said to be \emph{\tabCDo{k}} if \(\dia{k}=\Oqq\).
\edefn

\bpropl{F.P.abCDe}
 Let \(\seqska\in\Fggqka\), assume \(\kappa\geq1\), and let \(k\in\mn{1}{\kappa}\).
 Then \(\seqska\) is \tabCDo{k} if and only if \(\cia{k}^2=\cia{k}\).
 In this case, \(\dia{j}=\Oqq\), \(\cia{j}=\Oqq\), and \(\su{j}=\umg{j-1}=\mi{j-1}=\omg{j-1}\) for all \(j\in\mn{k+1}{\kappa}\).
\eprop
\bproof
 This is an immediate consequence of \rthmp{F.C.ecase}{F.C.ecase.d}.
\eproof

 Observe that in the situation of \rprop{F.P.abCDe} due to \(\cia{k}\lgeq\Oqq\) we have \(\cia{k}^\ad=\cia{k}\) and thus the condition \(\cia{k}^2=\cia{k}\) is equivalent to \(\cia{k}\) being a transformation matrix corresponding to an orthogonal projection, \tie{}, \(\cia{k}=\OPu{\ran{\cia{k}}}\).

\bdefnnl{\tcf{}~\zitaa{arXiv:1701.04246}{\cdefn{10.33}{42}}}{F.D.abZ}
 Let \(\seqska\) be a sequence of complex \tpqa{matrices} with \tfmf{} \(\seqmika\).
 Assume \(\kappa\geq1\) and let \(k\in\mn{1}{\kappa}\).
 Then \(\seqska\) is said to be \emph{\tabZo{k}} if \(\su{j}=\mi{j-1}\) for all \(j\in\mn{k}{\kappa}\).
\edefn

\bpropl{F.P.abZe}
 Let \(\seqska\in\Fggqka\), assume \(\kappa\geq1\), and let \(k\in\mn{1}{\kappa}\).
 Then \(\seqska\) is \tabZo{k} if and only if \(\cia{j}=\frac{1}{2}\Pd{j-1}\) for all \(j\in\mn{k}{\kappa}\).
 In this case, \(\dia{j}=\ek{\ba/4}^{j-k+1}\dia{k-1}\) and \(\Pd{j}=\Pd{k-1}\) for all \(j\in\mn{k}{\kappa}\).
\eprop
\bproof
 This is an immediate consequence of \rthmp{F.C.ecase}{F.C.ecase.c}.
\eproof

\section{The \tBTion{} in the class $\Fggqka$}\label{F.BT}
 To motivate the considerations of this and the next section, we give a matrix version of the \emph{change-of-variables formula}.
 Therefore, we consider a measurable function \(T\colon\dbX\to\dbY\) between two measurable spaces \((\dbX,\sigA)\) and \((\dbY,\sigB)\) and a \tnnH{} \tqqa{measure} \(\mu\) on \((\dbX,\sigA)\).
 For all \(Y\in\dbY\), let
\(\ek{T(\mu)}(Y)
 \defeq\mu\rk{T^\inv(Y)}\),
 where
\(T^\inv(Y)\defeq\setaa{x\in\dbX}{T(x)\in Y}\)
 is the preimage of \(Y\) under \(T\).
 Then \(T(\mu)\colon\sigB\to\Cggq\) is a \tnnH{} \tqqa{measure} on \((\dbY,\sigB)\), called the \emph{image measure of \(\mu\) under \(T\)}.

\bpropnl{\tcf{}~\zitaa{MR3133464}{\cprop{B.1}{3926}}}{142.P1100}
 Let \(g\colon\dbY\to\C\) be a measurable function.
 Then \(g\in\LnnH{T(\mu)}\) if and only if \(g\circ T\in\LnnH{\mu}\).
 In this case, \(\int_Yg\dif[T(\mu)]=\int_{T^\inv(Y)}(g\circ T)\dif\mu\) for all \(Y\in\sigB\).
\eprop

 In particular, consider an affine transformation of \(\R\) and let \(T\) be its restriction onto a non-empty Borel subset \(\Omega\) of \(\R\).
 Then the connection between the power moments \(\mpm{\sigma}{j}\defeq\int_\Omega x^j\sigma\rk{\dif x}\) of a \tnnH{} measure \(\sigma\) on \((\Omega,\BsAu{\Omega})\) and the power moments \(t_j\defeq\mpm{T\rk{\sigma}}{j}=\int_{T(\Omega)}y^j\ek{T(\sigma)}\rk{\dif y}\) of the image measure \(T(\sigma)\) of \(\sigma\) under \(T\) is given via the following transformation for sequences of complex matrices:

\bdefnl{H.D.BT}
 Let \(\phi,\psi\in\C\) and let \(\seqska \) be a sequence of complex \tpqa{matrices}.
 Then the sequence \(\seq{w_j}{j}{0}{\kappa}\) given by
 \(
  w_j
  \defg\sum_{\ell=0}^j\binom{j}{\ell}\psi^\ell\phi^{j-\ell}s_\ell
 \)
 is called the \emph{\tBTaaa{\phi}{\psi}{\(\seqska\)}}.
\edefn

\bcorl{B.L.BTmom}
 Let \(\phi,\psi\in\R\), let \(\Omega\in\BsAR\setminus\set{\emptyset}\), let \(\Theta\defeq\setaa{\psi x+\phi}{x\in\Omega}\), and let \(T\colon\Omega\to\Theta\) be defined by \(T(x)\defeq\psi x+\phi\).
 Let \(\seqska\) be a sequence of complex \tqqa{matrices} with \tBTaa{\phi}{\psi} \(\seq{w_j}{j}{0}{\kappa}\) and let \(\sigma\in\MggqOsg{\kappa}\).
 Then the image measure \(T(\sigma)\) of \(\sigma\) under \(T\) belongs to \(\Mggouaag{q}{\kappa}{\Theta}{\seq{w_j}{j}{0}{\kappa}}\).
\ecor
\bproof
 Use \rprop{142.P1100}.
\eproof

\breml{H.R.BTtr}
 Let \(\phi,\psi\in\C\) and let \(\seqska\) be a sequence of complex \tpqa{matrices} with \tBTaa{\phi}{\psi} \(\seq{w_j}{j}{0}{\kappa}\).
 For each \(k\in\mn{0}{\kappa}\), the matrix \(w_k\) is built from the matrices \(\su{0},\su{1},\dotsc,\su{k}\).
 In particular, for each \(m\in\mn{0}{\kappa}\), the \tBTaaa{\phi}{\psi}{\(\seqs{m}\)} coincides with \(\seq{w_j}{j}{0}{m}\).
\erem

 We are going to study the behavior of \tfc{s}, \tfd{s}, and \tfp{s} under binomial transformation of the underlying sequence \(\seqska\in\Fggqka\), starting with some simple observations concerning the arithmetic of the binomial transformation:

\bexaml{H.E.d^j*s}
 Let \(\phi,\psi\in\C\) and let \(\seqska\) be a sequence of complex \tpqa{matrices}.
 Then \(\seq{\phi^js_0}{j}{0}{\kappa}\) is the \tBTaaa{\phi}{0}{\(\seqska \)} and \(\seq{\psi^js_j}{j}{0}{\kappa}\) is the \tBTaaa{0}{\psi}{\(\seqska \)}.
\eexam

 Obviously, the \tBTion{} is linear in the following sense:

\breml{H.R.BTlin}
 Let \(\phi,\psi\in\C\) and let \(\seqska \) and \(\seqt{\kappa}\) be two sequences of complex \tpqa{matrices}.
 Denote by \(\seq{w_j}{j}{0}{\kappa}\) and by \(\seq{v_j}{j}{0}{\kappa}\) the \tBTaaa{\phi}{\psi}{\(\seqska\)} and of \(\seqt{\kappa}\), \tresp{}
 Then the sequence \(\seq{w_j+v_j}{j}{0}{\kappa}\) is the \tBTaaa{\phi}{\psi}{\(\seq{\su{j}+\tu{j}}{j}{0}{\kappa}\)}.
 Furthermore, for each \(\lambda\in\C\), the sequence \(\seq{\lambda w_j}{j}{0}{\kappa}\) is the \tBTaaa{\phi}{\psi}{\(\seq{\lambda\su{j}}{j}{0}{\kappa}\)}.
\erem

\bpropl{F.P.Fgg^b}
 Let \(\eta\in\R\), let \(\theta\in[0,\infp)\), and let \(\seqska\in\Fggqu{\kappa}\) with \tfpf{} \(\fpseqka\).
 Denote by \(\seq{w_j}{j}{0}{\kappa}\) the \tBTaaa{\eta}{\theta}{\(\seqska\)} and by \(\fgpseqka\) the \tfupfa{\theta\ug+\eta}{\theta\obg+\eta}{\(\seq{w_j}{j}{0}{\kappa}\)}.
 Then \(\seq{w_j}{j}{0}{\kappa}\in\Fggquu{\kappa}{\theta\ug+\eta}{\theta\obg+\eta}\).
 Furthermore,
\begin{align}
  \fgpu{2k}&=\theta^k\fpu{2k}&\text{ for all }k&\in\mn{0}{\kappa}\label{F.P.Fgg^b.A}
\intertext{and}
  \fgpu{2k-1}&=\theta^k\fpu{2k-1}&\text{ for all }k&\in\mn{1}{\kappa}.\label{F.P.Fgg^b.B}
\end{align}
\eprop
\bproof
 Denote by \(\kpseq{\kappa}\) the \tkpfa{\(\seqska\)}.
 From \rpropss{ab.R1532}{ab.R1535} we conclude \(\seqska\in\Kggqka\).
 According to~\zitaa{MR3133464}{\cthm{4.12}{3903}}, then \(\seq{\theta^j\kpu{j}}{j}{0}{\kappa}\) is the \tkupfa{\theta\ug+\eta}{\(\seq{w_j}{j}{0}{\kappa}\)}.
 By virtue of \rrem{F.R.fpkap}, hence \(\fgpu{4\ell}=\theta^{2\ell}\fpu{4\ell}\) for all \(\ell\in\NO\) with \(2\ell\leq\kappa\) and \(\fgpu{4\ell+1}=\theta^{2\ell+1}\fpu{4\ell+1}\) for all \(\ell\in\NO\) with \(2\ell+1\leq\kappa\).
 Now assume \(\kappa\geq1\).
 Denote by \(\lkpseq{\kappa-1}\) the \tkpfa{\(\seqsb{\kappa-1}\)} and by \(\seqr{\kappa-1}\) the \tBTaaa{\eta}{\theta}{\(\seqsb{\kappa-1}\)}.
 Because of \rprop{F.L.sabF}, the sequence \(\seqsb{\kappa-1}\) belongs to \(\Fggqu{\kappa-1}\) and, consequently, to \(\Kggqu{\kappa-1}\) by the same reasoning as above.
 According to~\zitaa{MR3133464}{\cthm{4.12}{3903}}, then \(\seq{\theta^j\lkpu{j}}{j}{0}{\kappa-1}\) is the \tkupfa{\theta\ug+\eta}{\(\seqr{\kappa-1}\)}.
 In view of \rrem{K.R.l*k}, the \tkupfa{\theta\ug+\eta}{\(\seq{\theta r_j}{j}{0}{\kappa-1}\)} hence coincides with \(\seq{\theta^{j+1}\lkpu{j}}{j}{0}{\kappa-1}\).
 Let the sequence \(\seq{t_j}{j}{0}{\kappa-1}\) be given by \(t_j\defeq\su{j+1}\) and denote by \(\seq{v_j}{j}{0}{\kappa-1}\) the \tBTaaa{\eta}{\theta}{\(\seq{t_j}{j}{0}{\kappa-1}\)}.
 For all \(j\in\mn{0}{\kappa-1}\), we have, because of~\zitaa{MR3133464}{\crem{4.7}{3902}}, then \(w_{j+1}-\eta w_{j}=\theta v_{j}\) and, by virtue of \(\sub{j}=\obg\su{j}-t_j\) and \rrem{H.R.BTlin}, therefore
 \(
  \theta r_j
  =\theta(\obg w_{j}-v_{j})
  =\theta\obg w_{j}-\theta v_{j}
  =(\theta\obg+\eta)w_{j}-w_{j+1}
 \).
 Thus, \(\seq{\theta^{j+1}\lkpu{j}}{j}{0}{\kappa-1}\) is the \tkupfa{\theta\ug+\eta}{\(\seq{(\theta\obg+\eta)w_{j}-w_{j+1}}{j}{0}{\kappa-1}\)}.
 \rrem{F.R.fpkap} yields then \(\fgpu{4\ell+2}=\theta^{2\ell+1}\fpu{4\ell+2}\) for all \(k\in\NO\) with \(2\ell\leq\kappa-1\) and \(\fgpu{4\ell+3}=\theta^{2\ell+2}\fpu{4\ell+3}\) for all \(k\in\NO\) with \(2\ell+1\leq\kappa-1\).
 Hence, we have shown \eqref{F.P.Fgg^b.A} and \eqref{F.P.Fgg^b.B}.
 In view of \(\theta\geq0\), then \(\seq{w_j}{j}{0}{\kappa}\in\Fggquu{\kappa}{\theta\ug+\eta}{\theta\obg+\eta}\) follows by virtue of \eqref{F.P.Fgg^b.A}, \eqref{F.P.Fgg^b.B}, and \rprop{F.P.FggFP}.
\eproof

\bcorl{F.C.dia^b}
 Let \(\eta\in\R\), let \(\theta\in[0,\infp)\), and let \(\seqska\in\Fggqu{\kappa}\).
 Denote by \(\seq{w_j}{j}{0}{\kappa}\) the \tBTaaa{\eta}{\theta}{\(\seqska \)}.
 Then \(\seq{\theta^{j+1}\dia{j}}{j}{0}{\kappa}\) coincides with the \tfudfa{\theta\ug+\eta}{\theta\obg+\eta}{\(\seq{w_j}{j}{0}{\kappa}\)}.
\ecor
\bproof
 In the case \(\kappa<\infp\), we first extend the sequence \(\seqska \) by virtue of \rcor{ab.R1011} to a sequence \(\seqs{\kappa+1}\in\Fggqu{\kappa+1}\).
 In view of \rremss{F.R.diatr}{H.R.BTtr}, then the assertion follows from \rrem{F.R.f2n-1} and \rprop{F.P.Fgg^b}.
\eproof

\bpropl{F.P.cia^b}
 Let \(\eta\in\R\), let \(\theta\in(0,\infp)\), and let \(\seqska\in\Fggqu{\kappa}\).
 Denote by \(\seq{w_j}{j}{0}{\kappa}\) the \tBTaaa{\eta}{\theta}{\(\seqska\)}.
 Then \(\seqciaka \) coincides with the \tfucfa{\theta\ug+\eta}{\theta\obg+\eta}{\(\seq{w_j}{j}{0}{\kappa}\)}.
\eprop
\bproof
 Denote by \(\fpseqka\) the \tfpfa{\(\seqska\)} and by \(\seq{\mathfrak{p}_j}{j}{0}{\kappa}\) the \tfucfa{\theta\ug+\eta}{\theta\obg+\eta}{\(\seq{w_j}{j}{0}{\kappa}\)}.
 In view of \rdefn{F.D.cia}, then \rprop{F.P.Fgg^b} yields \(\mathfrak{p}_0=\theta^0\fpu{0}=\cia{0}\).
 Using additionally \rcor{F.C.dia^b} and \rrem{A.R.l*A}, for all \(j\in\mn{1}{\kappa}\), we obtain
\begin{multline*}
  \mathfrak{p}_j
  =\ek*{\rk{\theta^j\dia{j-1}}^\varsqrt}^\mpi\rk{\theta^j\fpu{2j}}\ek*{\rk{\theta^j\dia{j-1}}^\varsqrt}^\mpi
  =\rk{\theta^{j/2}\dia{j-1}^\varsqrt}^\mpi\rk{\theta^j\fpu{2j}}\rk{\theta^{j/2}\dia{j-1}^\varsqrt}^\mpi\\
  =\ek*{\theta^{-j/2}\rk{\dia{j-1}^\varsqrt}^\mpi}\rk{\theta^j\fpu{2j}}\ek*{\theta^{-j/2}\rk{\dia{j-1}^\varsqrt}^\mpi}
  =\rk{\dia{j-1}^\varsqrt}^\mpi\fpu{2j}\rk{\dia{j-1}^\varsqrt}^\mpi
  =\cia{j}.\qedhere
\end{multline*}
\eproof

 To give an analogous result for \(\theta<0\), we first consider the pure reflection already treated in \rlem{F.L.-1^j*f}:

\bleml{F.L.-1^j*dc}
 Let \(\seqska\in\Fggqu{\kappa}\) with \tfdf{} \(\seqdiaka \) and \tfcf{} \(\seqciaka\).
 Let the sequence \(\seqr{\kappa}\) be given by \(r_j\defg(-1)^j\su{j}\).
 Denote by \(\seq{\mathfrak{p}_j}{j}{0}{\kappa}\) the \tfucfa{-\obg}{-\ug}{\(\seqr{\kappa}\)}.
 Then \(\seqr{\kappa}\) belongs to \(\Fggquu{\kappa}{-\obg}{-\ug}\) and \(\seqdiaka \) coincides with the \tfudfa{-\obg}{-\ug}{\(\seqr{\kappa}\)}.
 Furthermore, \(\mathfrak{p}_{2k}=\cia{2k}\) for all \(k\in\NO\) with \(2k\leq\kappa\) and \(\mathfrak{p}_{2k+1}=\Pd{2k}-\cia{2k+1}\) for all \(k\in\NO\) with \(2k+1\leq\kappa\).
\elem
\bproof
 Using \rlem{F.L.-1^j*f}, we conclude \(\seqr{\kappa}\in\Fggquu{\kappa}{-\obg}{-\ug}\) from \rprop{F.P.FggFP}.
 According to \rrem{F.R.d+=f} and \rlem{F.L.-1^j*f}, then \(\seqdiaka \) coincides with the \tfudfa{-\obg}{-\ug}{\(\seqr{\kappa}\)}.
 Denote by \(\fpseqka\) the \tfpfa{\(\seqska\)} and by \(\rfpseq{2\kappa}\) the \tfupfa{-\obg}{-\ug}{\(\seqr{\kappa}\)}.
 In view of \rdefn{F.D.cia} and \rlem{F.L.-1^j*f}, we obtain then \(\mathfrak{p}_0=\rfpu{0}=\fpu{0}=\cia{0}\) and
 \[
  \mathfrak{p}_{2k}
  =\rk{\dia{2k-1}^\varsqrt}^\mpi\rfpu{4k}\rk{\dia{2k-1}^\varsqrt}^\mpi
  =\rk{\dia{2k-1}^\varsqrt}^\mpi\fpu{4k}\rk{\dia{2k-1}^\varsqrt}^\mpi
  =\cia{2k}
 \]
 for all \(k\in\N\) with \(2k\leq\kappa\).
 Now we consider an arbitrary integer \(k\in\NO\) with \(2k+1\leq\kappa\).
 Because of \rlem{F.L.-1^j*f} and \rrem{F.R.f2n-1}, we have
 \(\rfpu{4k+2}
  =\fpu{4k+1}
  =\dia{2k}-\fpu{4k+2}\).
 Observe that the matrix \(\dia{2k}\) is \tnnH{}, according to \rprop{ab.C0929}.
 Hence, we conclude then \(\Pd{2k}=\rk{\dia{2k}^\varsqrt}^\mpi\dia{2k}\rk{\dia{2k}^\varsqrt}^\mpi
 \) in the same way as we derived \eqref{F.L.cs-es.34}.
 Consequently, we get
\[
  \mathfrak{p}_{2k+1}
  =\rk{\dia{2k}^\varsqrt}^\mpi\rfpu{4k+2}\rk{\dia{2k}^\varsqrt}^\mpi
  =\rk{\dia{2k}^\varsqrt}^\mpi\dia{2k}\rk{\dia{2k}^\varsqrt}^\mpi-\rk{\dia{2k}^\varsqrt}^\mpi\fpu{4k+2}\rk{\dia{2k}^\varsqrt}^\mpi
  =\Pd{2k}-\cia{2k+1}.\qedhere
\]
\eproof

\bpropl{F.P.FggBT}
 Let \(\eta\in\R\), let \(\theta\in(-\infty,0)\), and let \(\seqska\in\Fggqu{\kappa}\) with \tfdf{} \(\seqdiaka \) and \tfcf{} \(\seqciaka\).
 Denote by \(\seq{w_j}{j}{0}{\kappa}\) the \tBTaaa{\eta}{\theta}{\(\seqska\)} and by \(\seq{\mathfrak{q}_j}{j}{0}{\kappa}\) the \tfucfa{\theta\obg+\eta}{\theta\ug+\eta}{\(\seq{w_j}{j}{0}{\kappa}\)}.
 Then \(\seq{w_j}{j}{0}{\kappa}\) belongs to \(\Fggquu{\kappa}{\theta\obg+\eta}{\theta\ug+\eta}\) and \(\seq{\rk{-\theta}^{j+1}\dia{j}}{j}{0}{\kappa}\) coincides with the \tfudfa{\theta\obg+\eta}{\theta\ug+\eta}{\(\seq{w_j}{j}{0}{\kappa}\)}.
 Furthermore, \(\mathfrak{q}_{2k}=\cia{2k}\) for all \(k\in\NO\) with \(2k\leq\kappa\) and \(\mathfrak{q}_{2k+1}=\Pd{2k}-\cia{2k+1}\) for all \(k\in\NO\) with \(2k+1\leq\kappa\).
\eprop
\bproof
 Let the sequence \(\seqr{\kappa}\) be given by \(r_j\defg(-1)^j\su{j}\).
 In view of \rexam{H.E.d^j*s}, then \(\seqr{\kappa}\) is the \tBTaaa{0}{-1}{\(\seqska \)}.
 Denote by \(\seq{\mathfrak{p}_j}{j}{0}{\kappa}\) the \tfucfa{-\obg}{-\ug}{\(\seqr{\kappa}\)}.
 According to \rlem{F.L.-1^j*dc}, then \(\seqr{\kappa}\in\Fggquu{\kappa}{-\obg}{-\ug}\), the \tfudfa{-\obg}{-\ug}{\(\seqr{\kappa}\)} is \(\seqdiaka \), and furthermore \(\mathfrak{p}_{2k}=\cia{2k}\) for all \(k\in\NO\) with \(2k\leq\kappa\) and \(\mathfrak{p}_{2k+1}=\Pd{2k}-\cia{2k+1}\) for all \(k\in\NO\) with \(2k+1\leq\kappa\).
 Denote by \(\seq{u_j}{j}{0}{\kappa}\) the \tBTaaa{\eta}{-\theta}{\(\seqr{\kappa}\)}.
 In view of \(-\theta>0\), we can apply \rprop{F.P.Fgg^b} to the sequence \(\seqr{\kappa}\in\Fggquu{\kappa}{-\obg}{-\ug}\) and obtain \(\seq{u_j}{j}{0}{\kappa}\in\Fggquu{\kappa}{\obg\theta+\eta}{\theta\ug+\eta}\).
 Using \rcor{F.C.dia^b} in a similar way, we see that \(\seq{\rk{-\theta}^{j+1}\dia{j}}{j}{0}{\kappa}\) is the \tfudfa{\theta\obg+\eta}{\theta\ug+\eta}{\(\seq{u_j}{j}{0}{\kappa}\)}.
 Furthermore, by the same reasoning, \rprop{F.P.cia^b} yields that \(\seq{\mathfrak{p}_j}{j}{0}{\kappa}\) is the \tfucfa{\theta\obg+\eta}{\theta\ug+\eta}{\(\seq{u_j}{j}{0}{\kappa}\)}.
 Since \(\seq{w_j}{j}{0}{\kappa}\) coincides with \(\seq{u_j}{j}{0}{\kappa}\), by virtue of~\zitaa{MR2805417}{\clem{4.6}{461}}, the proof is complete.
\eproof

\section{Matricial canonical moments}\label{F.CM}
 On the basis of \rprop{I.P.ab8Fgg} we now rewrite our results in terms of matrix measures.
 In this way several results from the scalar theory of canonical moments can be generalized to the matrix case.

\bdefnl{F.D.meacia}
 Let \(\sigma\in\MggqF\) with \tfpmf{} \(\seqmpm{\sigma}\).
 Denote by \(\seqmcm{\sigma}\) the \tfcfa{\(\seqmpm{\sigma}\)} and by \(\seqmdm{\sigma}\) the \tfdfa{\(\seqmpm{\sigma}\)}.
 Then we call \(\seqmcm{\sigma}\) the \emph{\tfmcmfa{\(\sigma\)}} and \(\seqmdm{\sigma}\) the \emph{\tfmdmfa{\(\sigma\)}}.
\edefn

 In view of \rprop{I.P.ab8Fgg}, we can now translate several results obtained for \tFnnd{} sequences into the language of \tnnH{} measures on \(\ab\):
 
\bthml{F.T.Mabcia}
 Let \(\Pi_\ab\colon\MggqF\to\es{q}{\infi}{\bam}\) be defined by \(\sigma\mapsto\seqmcm{\sigma}\).
 Then the mapping \(\Pi_\ab\) is well defined and bijective.
\ethm
\bproof
 Use \rprop{I.P.ab8Fgg} and \rthm{F.T.Fggcia}.
\eproof

\bpropl{F.P.MFed}
 Let \(\sigma\in\MggqF\).
 Then \(\mcm{\sigma}{0}\lgeq\Oqq\) and \(\mcm{\sigma}{j}\in\matint{\Oqq}{\OPu{\ran{\mdm{\sigma}{j-1}}}}\) for all \(j\in\N\).
 Furthermore, \(\mdm{\sigma}{0}=\ba\mcm{\sigma}{0}\) and
\(
 \mdm{\sigma}{j}
 =\ba\rk{\mdm{\sigma}{j-1}}^\varsqrt \mcm{\sigma}{j}\rk{\Iq-\mcm{\sigma}{j}}\rk{\mdm{\sigma}{j-1}}^\varsqrt 
\)
 and, in particular,
\(
 \det\rk{\mdm{\sigma}{j}}
 =\ba^q\det\rk{\mdm{\sigma}{j-1}}\det\rk{\mcm{\sigma}{j}}\det\rk{\Iq-\mcm{\sigma}{j}}
\)
 for all \(j\in\N\).
\eprop
\bproof
 Use \rpropsss{I.P.ab8Fgg}{F.P.ed}{F.L.1-e}.
\eproof

 Consider a \tnnH{} measure \(\sigma\in\MggqF\).
 Because of \rprop{ab.C1101}, the column spaces \(\mathcal{R}_j\defeq\ran{\mdm{\sigma}{j}}\) build a nested decreasing sequence
\(
 \Cq
 \supseteq\mathcal{R}_0
 \supseteq\mathcal{R}_1
 \supseteq\mathcal{R}_2
 \supseteq\dotsb
 \supseteq\mathcal{R}_j
 \supseteq\dotsb
 \supseteq\set{\Ouu{q}{1}}
\),
 which, due to the finite dimensional situation, necessarily stabilizes at some index \(j_0\in\NO\), \tie{},
\(
 \mathcal{R}_{j_0-1}
 \supsetneqq\mathcal{R}_{j_0}
 =\mathcal{R}_{j_0+1}
 =\mathcal{R}_{j_0+2}
 =\dotsb
\)
 In particular, the matrices \(\mathbb{P}_j\defeq\OPu{\ran{\mdm{\sigma}{j}}}\) corresponding to the orthogonal projections onto the subspaces \(\mathcal{R}_j\) form, by virtue of \rrem{A.R.P<P}, a decreasing sequence
\[
 \Iq
 \lgeq\mathbb{P}_0
 \lgeq\mathbb{P}_1
 \lgeq\mathbb{P}_2
 \lgeq\dotsb
 \lgeq\mathbb{P}_{j_0-1}
 \lgeq\mathbb{P}_{j_0}
 =\mathbb{P}_{j_0+1}
 =\mathbb{P}_{j_0+2}
 =\dotsb
 \lgeq\Oqq
\]
 with \(\mathbb{P}_{j_0-1}\neq\mathbb{P}_{j_0}\).
 Consequently, the matricial intervals \(\matint{\Oqq}{\OPu{\ran{\mdm{\sigma}{j-1}}}}\) in \rprop{F.P.MFed} are nested decreasing and remain unchanged for \(j>j_0\).
 Moreover, the modified matricial interval lengths \(d_j\defeq\frac{4}{\bam}\mdm{\sigma}{j}\) build itself a decreasing sequence \(d_0\lgeq d_1\lgeq d_2\lgeq\dotsb\lgeq\Oqq\) of \tnnH{} \tqqa{matrices}, due to \rpropss{ab.C0929}{ab.P1057}, which necessarily converges to a \tnnH{} limit with column space contained in \(\mathcal{R}_{j_0}\).
 The situation is transparent for \tnnH{} measures \(\sigma\) concentrated on a finite subset of \(\ab\):
 
 Let \(\MggqmolO\) be the set of all \(\sigma\in\MggqO\) for which there exists a finite subset \(B\) of \(\Omega\) satisfying \(\sigma\rk{\Omega\setminus B}=\Oqq\).
 The elements of \(\MggqmolO\) are said to be \emph{molecular}.
 Obviously, \(\MggqmolO\) is the set of all \(\sigma\in\MggqO\) for which there exist an \(m\in\N\) and sequences \(\seq{\xi_\ell}{\ell}{1}{m}\) and \(\seq{A_\ell}{\ell}{1}{m}\) from \(\Omega\) and \(\Cggq\), \tresp{}, such that \(\sigma=\sum_{\ell=1}^m\Kronu{\xi_\ell}A_\ell\), where \(\Kronu{\xi_\ell}\) is the Dirac measure on \(\BsAu{\Omega}\) with unit mass at \(\xi_\ell\).
 In particular, we have \(\MggqmolO\subseteq\Mggoua{q}{\infi}{\Omega}\).

\bpropl{F.P.mol}
 Let \(\sigma\in\MggqF\).
 Then the following statements are equivalent:
\baeqi{0}
 \il{F.P.mol.i} \(\sigma\) is molecular.
 \il{F.P.mol.ii} \(\mdm{\sigma}{k}=\Oqq\) for some \(k\in\N\).
 \il{F.P.mol.iii} \(\mcm{\sigma}{\ell}=\Oqq\) for some \(\ell\in\N\).
 \il{F.P.mol.iv} \(\rk{\mcm{\sigma}{m}}^2=\mcm{\sigma}{m}\) for some \(m\in\N\).
\eaeqi
 If \(k\in\N\) fulfills \(\mdm{\sigma}{k}=\Oqq\), then \(\mdm{\sigma}{j}=\Oqq\) and \(\mcm{\sigma}{j}=\Oqq\) for all \(j\in\minf{k+1}\).
\eprop
\bproof
 Denote by \(\seqsinf\) the \tfpmfa{\(\sigma\)} and by \(\fpseq{\infi}\) the \tfpfa{\(\seqsinf\)}.
 In view of \rprop{I.P.ab8Fgg}, we have \(\seqsinf\in\Fggqinf\).
 By virtue of \rdefn{F.D.meacia}, furthermore \(\seqmcm{\sigma}\) is the \tfcfa{\(\seqsinf\)} and \(\seqmdm{\sigma}\) is the \tfdfa{\(\seqsinf\)}.
 Consequently, \rthmp{F.C.ecase}{F.C.ecase.d} yields the equivalence \aequ{F.P.mol.ii}{F.P.mol.iv}, whereas the implication \impl{F.P.mol.iii}{F.P.mol.ii} can be seen from \rpartss{F.C.ecase.a}{F.C.ecase.b} of \rthm{F.C.ecase}.
 
\bimp{F.P.mol.i}{F.P.mol.iii}
 Suppose that \(\sigma\) is molecular, \tie{}, \(\sigma=\sum_{\ell=1}^m\Kronu{\xi_\ell}A_\ell\) for some \(m\in\N\), some numbers \(\xi_1<\xi_2<\dotsc<\xi_m\) from \(\ab\), and some matrices \(A_1,A_2.\dotsc,A_m\) from \(\Cggq\).
 For all \(j\in\NO\), then \(\su{j}=\sum_{\ell=1}^m\xi_\ell^jA_\ell\) holds true.
 Using~\zitaa{MR2570113}{\clem{2.40(d)}{777}}, we can thus conclude \(\rank\Hu{m}=\rank\Hu{m-1}\).
 By virtue of \rprop{F.L.FPdet}, hence \(\rank\fpu{4m}=0\), \tie{}, \(\fpu{4m}=\Oqq\) follows.
 In view of \rdefn{F.D.cia}, we obtain therefore \(\mcm{\sigma}{2m}=\Oqq\). 
\eimp

\bimp{F.P.mol.ii}{F.P.mol.i}
 Let \(k\in\N\) with \(\mdm{\sigma}{k}=\Oqq\).
 By virtue of \rthmp{F.C.ecase}{F.C.ecase.d}, then \(\mdm{\sigma}{j}=\Oqq\) and \(\mcm{\sigma}{j}=\Oqq\) for all \(j\in\minf{k+1}\).
 Hence, using \rprop{F.L.cia-fp}, we obtain \(\fpu{2j}=\Oqq\) for all \(j\in\minf{k+1}\).
 In particular \(\fpu{4\ell}=\Oqq\) for all \(\ell\in\minf{k+1}\), which, by virtue of \rprop{F.L.FPdet}, implies \(\rank\Hu{n}=\rank\Hu{k}\) for all \(n\in\minf{k}\).
 Consequently, \(\rank\Hu{\rk{k+1}q}=\rank\Hu{k}\).
 Since \(\Hu{k}\) is a \taaa{\rk{k+1}q}{\rk{k+1}q}{matrix}, then \(\rank\Hu{\rk{k+1}q}\leq\rk{k+1}q\) follows.
 Taking into account that \(\mu\colon\BsAR\to\Cggq\) defined by \(\mu\rk{B}\defeq\sigma\rk{B\cap\ab}\) belongs to \(\MggquR{\infi}\) and that \(\int_\R x^j\mu\rk{\dif x}=\su{j}\) holds true for all \(j\in\NO\), the application of~\zitaa{MR2570113}{\clem{C.3}{824}} shows that \(\mu\) and hence \(\sigma\) is molecular.
\eimp
\eproof

 We now consider the behavior of the sequences \(\seqmcm{\sigma}\) and  \(\seqmdm{\sigma}\) under some transformations of the underlying \tnnH{} measure \(\sigma\in\MggqF\).
 The following three remarks are readily checked:

\breml{F.R.sm+sm}
 Let \(\sigma,\mu\in\MggqF\).
 Then \(\sigma+\mu\) belongs to \(\MggqF\) and has \tfpmf{} \(\seq{\mpm{\sigma}{j}+\mpm{\mu}{j}}{j}{0}{\infi}\).
\erem
 
\breml{F.R.VsmV}
 Let \(A\in\Cpq\) and let \(\sigma\in\MggqF\).
 Then \(A\sigma A^\ad\) belongs to \(\Mggoa{p}{\ab}\) and has \tfpmf{} \(\seq{A\mpm{\sigma}{j}A^\ad}{j}{0}{\infi}\).
\erem

\breml{F.R.l*sm}
 Let \(\lambda\in[0,\infp)\) and let \(\sigma\in\MggqF\).
 Then \(\lambda\sigma\) belongs to \(\MggqF\) and has \tfpmf{} \(\seq{\lambda\mpm{\sigma}{j}}{j}{0}{\infi}\).
\erem

\breml{F.R.sd+sd}
 Let \(\sigma,\mu\in\MggqF\).
 From \rrem{F.R.sm+sm} and \rprop{F.R.d+d}, we easily see that \(\sigma+\mu\in\MggqF\) and \(\mdm{\sigma+\mu}{j}\lgeq\mdm{\sigma}{j}+\mdm{\mu}{j}\) for all \(j\in\NO\).
\erem

\breml{F.R.A.sd.A^*}
 Let \(A\in\Cpq\) and let \(\sigma\in\MggqF\).
 \rrem{F.R.VsmV} and \rprop{F.R.A.d.A^*} show that \(A\sigma A^\ad\in\Mggoa{p}{\ab}\) and \(\mdm{A\sigma A^\ad}{j}\lgeq A\mdm{\sigma}{j}A^\ad\) for all \(j\in\NO\).
\erem

\bleml{F.L.VsV}
 Let \(V\in\Cqp\) with \(VV^\ad=\Iq\) and let \(\sigma\in\MggqF\).
 Then \(V^\ad\sigma V\) belongs to \(\Mggoa{p}{\ab}\) and has \tfmcmf{} \(\seq{V^\ad\mcm{\sigma}{j}V}{j}{0}{\infi}\) and \tfmdmf{} \(\seq{V^\ad\mdm{\sigma}{j}V}{j}{0}{\infi}\).
\elem
\bproof
 In view of \rprop{I.P.ab8Fgg} and \rdefn{F.D.meacia}, this is an immediate consequence of \rremss{F.R.VsmV}{F.R.U.d.V} and \rlem{F.R.UeU^*}.
\eproof

\bleml{F.L.l*s}
 Let \(\lambda\in(0,\infp)\) and let \(\sigma\in\MggqF\).
 Then \(\lambda\sigma\) belongs to \(\MggqF\) and has \tfmdmf{} \(\seq{\lambda\mdm{\sigma}{j}}{j}{0}{\infi}\).
 Furthermore, \(\mcm{\lambda\sigma}{0}=\lambda\mcm{\sigma}{0}\) and \(\mcm{\lambda\sigma}{j}=\mcm{\sigma}{j}\) for all \(j\in\N\).
\elem
\bproof
 In view of \rprop{I.P.ab8Fgg} and \rdefn{F.D.meacia}, this is an immediate consequence of \rremsss{F.R.l*sm}{F.R.l*d}{F.R.l*e}.
\eproof

\bpropnl{for the scalar case, \tcf{}~\zita{MR0254899} or~\zitaa{MR1468473}{\cthm{1.3.2}{13}}}{F.P.DS132}
 Let \(\eta\in\R\), let \(\theta\in(0,\infp)\), and let \(T\colon\ab\to[\theta\ug+\eta,\theta\obg+\eta]\) be defined by \(T(x)\defeq\theta x+\eta\).
 Furthermore, let \(\sigma\in\MggqF\) with image measure \(T\rk{\sigma}\) under \(T\).
 Then \(T\rk{\sigma}\) belongs to \(\Mggoa{q}{[\theta\ug+\eta,\theta\obg+\eta]}\), the sequence \(\seqmcm{\sigma}\) coincides with the \tfmcmfa{\(T\rk{\sigma}\)}, and \(\seq{\theta^{j+1}\mdm{\sigma}{j}}{j}{0}{\infi}\) is the \tfmdmfa{\(T\rk{\sigma}\)}.
\eprop
\bproof
 Denote by \(\seqsinf\) the \tfpmfa{\(\sigma\)}.
 In view of \rprop{I.P.ab8Fgg}, we have \(\seqsinf\in\Fggqinf\).
 Furthermore, \(\seqmcm{\sigma}\) is the \tfcfa{\(\seqsinf\)} and \(\seqmdm{\sigma}\) is the \tfdfa{\(\seqsinf\)}, by \rdefn{F.D.meacia}.
 Denote by \(\seq{w_j}{j}{0}{\infi}\) the \tBTaaa{\eta}{\theta}{\(\seqsinf\)}.
 From \rprop{F.P.cia^b} and \rcor{F.C.dia^b} we conclude then, that \(\seqmcm{\sigma}\) coincides with the \tfucfa{\theta\ug+\eta}{\theta\obg+\eta}{\(\seq{w_j}{j}{0}{\infi}\)} and that \(\seq{\theta^{j+1}\mdm{\sigma}{j}}{j}{0}{\infi}\) is the \tfudfa{\theta\ug+\eta}{\theta\obg+\eta}{\(\seq{w_j}{j}{0}{\infi}\)}.
 By virtue of \rcor{B.L.BTmom}, we see that \(T\rk{\sigma}\) belongs to \(\Mggoa{q}{[\theta\ug+\eta,\theta\obg+\eta]}\) and that \(\seq{w_j}{j}{0}{\infi}\) is the \tfpmfa{\(T\rk{\sigma}\)}.
 Taking again into account \rdefn{F.D.meacia}, the proof is complete.
\eproof

 Using \rprop{F.P.FggBT} instead of \rprop{F.P.cia^b} and \rcor{F.C.dia^b}, we can prove the following result in a similar way:

\bpropnl{for the scalar case, \tcf{}~\zita{MR0254899} or~\zitaa{MR1468473}{\cthm{1.3.3}{14}}}{F.P.DS133}
 Let \(\eta\in\R\), let \(\theta\in(-\infty,0)\), and let \(T\colon\ab\to[\theta\obg+\eta,\theta\ug+\eta]\) be defined by \(T(x)\defeq\theta x+\eta\).
 Let \(\sigma\in\MggqF\) with image measure \(T\rk{\sigma}\) under \(T\).
 For all \(k\in\NO\), then \(\mcm{T\rk{\sigma}}{2k}=\mcm{\sigma}{2k}\) and \(\mcm{T\rk{\sigma}}{2k+1}=\OPu{\ran{\mdm{\sigma}{2k}}}-\mcm{\sigma}{2k+1}\).
 Furthermore, \(\seq{\rk{-\theta}^{j+1}\mdm{\sigma}{j}}{j}{0}{\infi}\) is the \tfmdmfa{\(T\rk{\sigma}\)}.
\eprop

 In the remaining part of this section we describe symmetry for \tnnH{} measures on \(\rk{\ab,\BsAF}\) in terms of their \tmcm{s}.
 To introduce the corresponding notion of symmetry, we consider the affine bijective function \(R\colon\ab\to\ab\) defined by \(R(x)\defeq\ug+\obg-x\).
 A measure \(\sigma\in\MggqF\) is called \emph{symmetric}, if it coincides with its image measure under \(R\), \tie{}, if \(R(\sigma)=\sigma\) holds true.
 
\bdefnl{F.D.symseq}
 Let \(\eta\in\C\).
 We call a sequence \(\seqska\) of complex \tpqa{matrices} \emph{symmetric with respect to \(\eta\)}, if it coincides with its \tBTaa{\eta}{-1}.
\edefn

\bleml{F.L.sym}
 A measure \(\sigma\in\MggqF\) is symmetric if and only if \(\seqmpm{\sigma}\) is symmetric with respect to \(\ug+\obg\).
\elem
\bproof
 Let \(\eta\defeq\ug+\obg\) and let \(R\colon\ab\to\ab\) be defined by \(R(x)\defeq\eta-x\).
 Furthermore, denote by \(\rho\) the image measure of \(\sigma\) under \(R\) and by \(\seqr{\infi}\) the \tBTaaa{\eta}{-1}{\(\seqmpm{\sigma}\)}.
 According to \rcor{B.L.BTmom}, then \(\rho\in\MggqFag{r}{\infi}\), \tie{}, \(\rho\in\MggqF\) and \(\seqr{\infi}\) is the \tfpmfa{\(\rho\)}.
 From \rprop{I.P.ab8Fgg} we conclude hence, that \(\sigma\) is symmetric if and only if \(\mpm{\sigma}{j}=r_j\) for all \(j\in\NO\).
\eproof

\bpropl{F.P.symseq}
 Let \(\seqska\in\Fggqu{\kappa}\) with \tfpf{} \(\fpseqka\), \tfdf{} \(\seqdiaka\), \tfmf{} \(\seq{\mi{j}}{j}{0}{\kappa}\), and \tfcf{} \(\seqciaka\).
 Then the following statements are equivalent:
\baeqi{0}
 \il{F.P.symseq.i} \(\seqska\) is symmetric with respect to \(\ug+\obg\).
 \il{F.P.symseq.v} \(\fpu{4k+1}=\fpu{4k+2}\) for all \(k\in\NO\) with \(2k+1\leq\kappa\).
 \il{F.P.symseq.ii} \(\dia{2k+1}=\frac{\bam}{4}\dia{2k}\) for all \(k\in\NO\) with \(2k+1\leq\kappa\).
 \il{F.P.symseq.iv} \(\cia{2k+1}=\frac{1}{2}\Pd{2k}\) for all \(k\in\NO\) with \(2k+1\leq\kappa\).
 \il{F.P.symseq.iii} \(\su{2k+1}=\mi{2k}\) for all \(k\in\NO\) with \(2k+1\leq\kappa\).
\eaeqi
\eprop
\bproof
 Let \(\eta\defeq\ug+\obg\).
 Observe, that \(-\obg+\eta=\ug\) and \(-\ug+\eta=\obg\).
 Denote by \(\seqr{\kappa}\) the \tBTaaa{\eta}{-1}{\(\seqska\)} and by \(\seq{\mathfrak{q}_j}{j}{0}{\kappa}\) the \tfucfa{\ug}{\obg}{\(\seqr{\kappa}\)}.
 Applying \rprop{F.P.FggBT} with \(\theta=-1\), we conclude then that \(\seqr{\kappa}\) belongs to \(\Fggquu{\kappa}{\ug}{\obg}\) as well, and, furthermore, that \(\mathfrak{q}_{2k}=\cia{2k}\) for all \(k\in\NO\) with \(2k\leq\kappa\) and that \(\mathfrak{q}_{2k+1}=\Pd{2k}-\cia{2k+1}\) for all \(k\in\NO\) with \(2k+1\leq\kappa\).

\bimp{F.P.symseq.i}{F.P.symseq.iv}
 Suppose that \(\seqska\) is symmetric with respect to \(\eta\).
 According to \rdefn{F.D.symseq}, then \(\su{j}=r_j\) for all \(j\in\mn{0}{\kappa}\).
 In particular, \(\cia{2k+1}=\mathfrak{q}_{2k+1}\) for all \(k\in\NO\) with \(2k+1\leq\kappa\).
 Consequently, \(\cia{2k+1}=\frac{1}{2}\Pd{2k}\) for all \(k\in\NO\) with \(2k+1\leq\kappa\).
\eimp

\bimp{F.P.symseq.iv}{F.P.symseq.i}
 Suppose that \(\cia{2k+1}=\frac{1}{2}\Pd{2k}\) is fulfilled for all \(k\in\NO\) with \(2k+1\leq\kappa\).
 Then \(\mathfrak{q}_{2k+1}=\frac{1}{2}\Pd{2k}\) follows for all \(k\in\NO\) with \(2k+1\leq\kappa\).
 Hence, \(\cia{j}=\mathfrak{q}_{j}\) for all \(j\in\mn{0}{\kappa}\).
 Since \(\seqska\) and \(\seqr{\kappa}\) both belong to \(\Fggqka\), we obtain from \rthm{F.T.Fggcia} therefore \(\su{j}=r_j\) for all \(j\in\mn{0}{\kappa}\).
 Hence, \(\seqska\) is symmetric with respect to \(\eta\).
\eimp

 The equivalence of~\ref{F.P.symseq.iv} and~\ref{F.P.symseq.iii} follows from \rthmp{F.C.ecase}{F.C.ecase.c}, whereas \rprop{ab.P1057} yields the equivalence of~\ref{F.P.symseq.iii} and~\ref{F.P.symseq.ii}.
 The equivalence of~\ref{F.P.symseq.v} and~\ref{F.P.symseq.iii} is a consequence of \rdefnsss{F.D.fpf}{F.D.mi}{F.N.AB}.
\eproof

\bcorl{F.C.s1=0}
 Let \(\rho\in(0,\infp)\) and let \(\seqska\in\Fgguuuu{q}{\kappa}{-\rho}{\rho}\) with \tfudf{-\rho}{\rho} \(\seqdiaka \), \tfumf{-\rho}{\rho} \(\seq{\mi{j}}{j}{0}{\kappa}\), and \tfucf{-\rho}{\rho} \(\seqciaka\).
 Then \(\seqska\) is symmetric with respect to \(0\) if and only if \(\su{2k+1}=\Oqq\) for all \(k\in\NO\) with \(2k+1\leq\kappa\).
 In this case, \(\fpu{4k+1}=\fpu{4k+2}=\omg{2k}=-\umg{2k}\), \(\dia{2k}=2\omg{2k}\), \(\dia{2k+1}=\frac{\bam}{2}\omg{2k}\), \(\mi{2k}=\Oqq\), and \(\cia{2k+1}=\frac{1}{2}\OPu{\ran{\omg{2k}}}\) for all \(k\in\NO\) with \(2k+1\leq\kappa\).
\ecor
\bproof
 According to \rexam{H.E.d^j*s}, the sequence \(\seq{(-1)^j\su{j}}{j}{0}{\kappa}\) is the \tBTaaa{0}{-1}{\(\seqska\)}.
 In view of \rdefn{F.D.symseq}, consequently \(\seqska\) is symmetric with respect to \(0\) if and only if \(\su{j}=(-1)^j\su{j}\) for all \(j\in\mn{0}{\kappa}\).
 The latter is obviously equivalent to the condition that \(\su{2k+1}=\Oqq\) for all \(k\in\NO\) with \(2k+1\leq\kappa\) holds true.
 Now suppose that \(\seqska\) is symmetric with respect to \(0\).
 Consider an arbitrary \(k\in\NO\) with \(2k+1\leq\kappa\).
 Then \(\su{2k+1}=\Oqq\), whereas \rprop{F.P.symseq} yields \(\su{2k+1}=\mi{2k}\).
 Consequently, \(\mi{2k}=\Oqq\).
 In view of \rdefn{F.D.mi}, then \(\omg{2k}=-\umg{2k}\) follows.
 Hence, by virtue of \rdefnsss{F.D.dia}{F.D.fpf}{F.N.AB}, we obtain \(\dia{2k}=2\omg{2k}\) and, furthermore, \(\fpu{4k+1}=\usc{2k+1}=\su{2k+1}-\umg{2k}=\omg{2k}\) and \(\fpu{4k+2}=\osc{2k+1}=\omg{2k}-\su{2k+1}=\omg{2k}\).
 Taking again into account \rprop{F.P.symseq}, we can conclude \(\dia{2k+1}=\frac{\bam}{4}\dia{2k}=\frac{\bam}{2}\omg{2k}\) and \(\cia{2k+1}=\frac{1}{2}\Pd{2k}=\frac{1}{2}\OPu{\ran{\omg{2k}}}\).
\eproof

\bpropnl{for the scalar case, \tcf{}~\zita{MR0254899} or~\zitaa{MR1468473}{\ccor{1.3.4}{14}}}{F.P.symmea}
 A measure \(\sigma\in\MggqF\) is symmetric if and only if \(\mcm{\sigma}{2k+1}=\frac{1}{2}\OPu{\ran{\mdm{\sigma}{2k}}}\) for all \(k\in\NO\).
\eprop
\bproof
 Denote by \(\seqsinf\) the \tfpmfa{\(\sigma\)}.
 In view of \rprop{I.P.ab8Fgg}, we have \(\seqsinf\in\Fggqinf\).
 By virtue of \rdefn{F.D.meacia}, we see that \(\seqmcm{\sigma}\) is the \tfcfa{\(\seqsinf\)} and that \(\seqmdm{\sigma}\) is the \tfdfa{\(\seqsinf\)}.
 Consequently, \rprop{F.P.symseq} yields that \(\seqsinf\) is symmetric with respect to \(\ug+\obg\) if and only if \(\mcm{\sigma}{2k+1}=\frac{1}{2}\OPu{\ran{\mdm{\sigma}{2k}}}\) for all \(k\in\NO\).
 Applying \rlem{F.L.sym} completes the proof.
\eproof

\section{Interrelations with Dette's and Studden's approach}\label{F.s2.cm}
 As already mentioned in the introduction, there exist far reaching connections of our considerations to the concept of \emph{canonical moments}.
 In this section, we are going to derive generalizations of corresponding results from~\zita{MR1883272} for the general case of not necessarily invertible \tbHms{}, \tie{} for sequences \(\seqs{\ell}\in\Fggqu{\ell}\) corresponding to points \(\col\rk{S_0,S_1,\dotsc,S_{\ell}}\) from the moment space \eqref{msp}, not necessarily located in its interior.
 In~\zita{MR1883272}, the authors considered the interval \(\Omega=[0,1]\), \tie{}, \(\ug=0\) and \(\obg=1\).
 Using our notation, we have for the \tbHms{} \(\underline{H}_{n}\) and \(\overline{H}_{n}\), introduced in~\zitaa{MR1883272}{\ceqssptp{2.11}{2.12}{176}{177}}, the relations \(\underline{H}_{2m}=\Hu{m}\), \(\overline{H}_{2m}=\Ku{m-1}-\Gu{m-1}=\Hu{0,m-1,1}\), \(\underline{H}_{2m+1}=\Ku{m}=\Hu{0,m}\), and \(\overline{H}_{2m+1}=\Hu{m}-\Ku{m}=\Hu{m,1}\).
 Furthermore, for the matrices \(S_{n+1}^-\) and \(S_{n+1}^+\) introduced in~\zitaa{MR1883272}{\ceqss{2.13}{2.14}{177}}, we obtain the connections \(S_{2m}^-=\Tripu{m}=\umg{2m-1}\), \(S_{2m+1}^-=\Tripu{0,m}=\umg{2m}\), \(S_{2m}^+=\su{2m-1}-\Tripu{0,m-1,1}=\omg{2m-1}\), and \(S_{2m+1}^+=\su{2m}-\Tripu{m,1}=\omg{2m}\).
 The differences \(D_k\) from~\zitaa{MR1883272}{\ceq{2.17}{178}} can be then written as
\(
 D_k
 =S_k^+-S_k^-
 =\omg{k-1}-\umg{k-1}
 =\dia{k-1}
\).
 Because of \(S_k-S_k^-=\su{k}-\umg{k-1}=\usc{k}\) and \(S_k^+-S_k=\omg{k-1}-\su{k}=\osc{k}\), consequently the below introduced matrices \(\dsu{k}\) and \(\dsv{k}\) coincide with the matrices \(U_k\) and \(V_k\), \tresp{}, given in~\zitaa{MR1883272}{\ceqss{2.16}{2.18}{178}} for an interior point of the moment space only.
 In~\zita{MR1883272} the matrix \(U_k\) is called the \emph{\(k\)th matrix canonical moment}. 

 Let \(\seqska\) be a sequence of complex \tpqa{matrices} with \tfdf{} \(\seqdiaka \) and assume \(\kappa\geq1\).
 For all \(k\in\mn{1}{\kappa}\), then let
\begin{align}\label{F.G.cm}
 \dsu{k}&\defeq\dia{k-1}^\mpi\usc{k}&
 &\text{and}&
 \dsv{k}&\defeq\dia{k-1}^\mpi\osc{k}.
\end{align}
 Between \(\dsu{k}\) and \(\dsv{k}\) we have the connection \(\dsu{k}+\dsv{k}=\dia{k-1}^\mpi\dia{k-1}\), by virtue of \rrem{ab.R1420}.
 Consequently, we get
\begin{align}\label{F.G.u+v}
 \dsv{k}&=\dia{k-1}^\mpi\dia{k-1}-\dsu{k}&
 &\text{and}&
 \dsu{k}&=\dia{k-1}^\mpi\dia{k-1}-\dsv{k}.
\end{align}

\bleml{F.R.duv}
 Let \(\seqska\in\Fggqka\).
 For all \(k\in\mn{1}{\kappa}\), then \(\dia{k-1}\dsu{k}=\usc{k}\), \(\ba\osc{k}\dsu{k}=\dia{k}\), \(\dia{k-1}\dsv{k}=\osc{k}\), and \(  \ba\usc{k}\dsv{k}=\dia{k}\).
\elem
\bproof
 We consider an arbitrary integer \(k\in\mn{1}{\kappa}\).
 According to~\zitaa{arXiv:1701.04246}{\cprop{10.18}{38}}, we have \(\ran{\dia{k-1}}=\ran{\usc{k}}+\ran{\osc{k}}\).
 Hence, \(\ran{\usc{k}}\subseteq\ran{\dia{k-1}}\).
 Because of \rremp{R.AA+B}{R.AA+B.a}, then \(\dia{k-1}\dia{k-1}^\mpi\usc{k}=\usc{k}\) follows.
 By virtue of \eqref{F.G.cm}, consequently \(\dia{k-1}\dsu{k}=\usc{k}\).
 In view of \eqref{F.G.cm}, furthermore \(\ba\osc{k}\dsu{k}=\ba\osc{k}\dia{k-1}^\mpi\usc{k}=\dia{k}\), because of \rprop{ab.C1343}.
 The identities concerning \(\dsv{k}\) can be proved in a similar way.
\eproof

 By virtue of \rdefn{F.N.AB}, we obtain in the setting of \rlem{F.R.duv} consequently
\(\su{k}=\umg{k-1}+\dia{k-1}\dsu{k}\) and \(\su{k}=\omg{k-1}-\dia{k-1}\dsv{k}\).

 For arbitrarily given \(m\in\N\) and complex \tqqa{matrices} \(A_1,A_2,\dotsc,A_m\), we write
\(
 \prodr_{\ell=1}^m A_\ell
 \defeq A_1A_2\dotsm A_m
\)
 for their product with indicated order.

\bthmnl{\tcf{}~\zitaa{MR1883272}{\cthm{2.7}{178}} for the interior point case}{F.T.diax}
 Let \(\seqska\in\Fggqka\) and assume \(\kappa\geq1\).
 For all \(j\in\mn{1}{\kappa}\), then \(\dsu{j}\dsv{j}=\dsv{j}\dsu{j}\) and
 \beql{F.T.diax.A}
  \dia{j}
  =\ba^{j+1}\su{0}\prodr_{\ell=1}^j\dsu{\ell}\dsv{\ell}.
 \eeq
\ethm
\bproof
 According to \rlem{F.R.duv}, we have
 \(
  \ba\dia{j-1}\dsu{j}\dsv{j}
  =\ba\usc{j}\dsv{j}
  =\dia{j}
 \)
 for all \(j\in\mn{1}{\kappa}\).
 Furthermore, \(\ba\su{0}=\dia{0}\), by virtue of \eqref{F.G.d01}.
 Hence, we get \eqref{F.T.diax.A} for all \(j\in\mn{1}{\kappa}\).
 Because of \(\bam>0\), we conclude \(\usc{j}\dia{j-1}^\mpi\osc{j}=\osc{j}\dia{j-1}^\mpi\usc{j}\) for all \(j\in\mn{1}{\kappa}\) from \rprop{ab.C1343}.
 In view of \eqref{F.G.cm}, consequently \(\dsu{j}\dsv{j}=\dsv{j}\dsu{j}\).
\eproof

 Observe that the identities~\zitaa{MR1883272}{\ceqss{2.20}{2.21}{178}}, which were essential for the proof of~\zitaa{MR1883272}{\cthm{2.7}{178}}, correspond to the representations from \rprop{ab.C1343} and \rthm{ab.P1422}.
 For the characterization~\zitaa{MR1883272}{\cthm{4.1}{189}} of systems of monic matrix polynomials, orthogonal with respect to a \tnnH{} measure on \(([0,1],\BsAu{[0,1]})\), the authors introduced in~\zitaa{MR1883272}{\cfo{2.33}{183}} the following matrices:

 Let \(\seqska\) be a sequence of complex \tpqa{matrices} and assume \(\kappa\geq1\).
 Then let \(\dsz{1}\defeq\dsu{1}\), \(\dsc{1}\defeq\dsv{1}\), and, for all \(k\in\mn{2}{\kappa}\), let \(\dsz{k}\defeq\dsv{k-1}\dsu{k}\) and \(\dsc{k}\defeq\dsu{k-1}\dsv{k}\).
 In view of \eqref{F.G.u+v} and \rrem{ab.R1052}, the sequences \(\seq{\dsz{k}}{k}{1}{\kappa}\) and \(\seq{\dsc{k}}{k}{1}{\kappa}\) can be interpreted as matricial \emph{chain sequences} (\tcf{}~\zitaa{MR1883272}{\cthm{4.1(b)}{189}}).

\bthmnl{\tcf{}~\zitaa{MR1883272}{\cthm{2.8}{183}} for the interior point case}{F.T.chain}
 Let \(\seqska\in\Fggqka\) and assume \(\kappa\geq1\).
 Then
 \(
  \ba\usc{k-1}\dsz{k}
  =\usc{k}
 \)
 for all \(k\in\mn{1}{\kappa}\) and
 \(
  \ba\osc{k-1}\dsc{k}
  =\osc{k}
 \)
 for all \(k\in\mn{2}{\kappa}\).
\ethm
\bproof
 From \rthm{ab.P1422} and \rlem{F.R.duv}, we conclude
 \(\ba\usc{0}\dsz{1}
  =\dia{0}\dsz{1}
  =\dia{0}\dsu{1}
  =\usc{1}\).
 Now assume \(\kappa\geq2\).
 We consider an arbitrary integer \(k\in\mn{2}{\kappa}\).
 Using \rlem{F.R.duv}, we obtain then
 \(\ba\usc{k-1}\dsz{k}
  =\ba\usc{k-1}\dsv{k-1}\dsu{k}
  =\dia{k-1}\dsu{k}
  =\usc{k}\).
 In a similar way, we get \(\ba\osc{k-1}\dsc{k}=\osc{k}\).
\eproof

\bcornl{\tcf{}~\zitaa{MR1883272}{\ceqss{2.34}{2.35}{183}} for the interior point case}{F.C.DS23435}
 If \(\seqska\in\Fggqka\), then \(\usc{k}=\ba^k\su{0}\prodr_{\ell=1}^k\dsz{\ell}\) and \(\osc{k}=\ba^k\su{0}\prodr_{\ell=1}^k\dsc{\ell}\) for all \(k\in\mn{1}{\kappa}\).
\ecor
\bproof
 Note that \(\ba\su{0}\dsz{1}=\dia{0}\dsz{1}=\dia{0}\dsu{1}=\usc{1}\) and likewise \(\ba\su{0}\dsc{1}=\osc{1}\), because of \eqref{F.G.d01} and \rlem{F.R.duv}.
 Applying \rthm{F.T.chain} completes the proof.
\eproof

 Observe that the identities stated in \rcor{F.C.DS23435} were essentially used in the proof of~\zitaa{MR1883272}{\cthm{2.8}{183}}.

\appendix
\section{Some facts from matrix theory}\label{M.S}
 This appendix contains a summary of results from matrix theory which are used in this paper.
 What concerns results on Moore--Penrose inverses of matrices we refer to~\zitaa{MR1152328}{\csec{1}}.

\breml{A.R.A++*}%
 If \(A\in\Cpq\), then \(\rk{A^\mpi}^\mpi=A\) and \(\rk{A^\ad}^\mpi=\rk{A^\mpi}^\ad\).
\erem

\breml{A.R.l*A}
 Let \(\eta\in\C\) and let \(A\in\Cpq\).
 Then \(\rk{\eta A}^\mpi=\eta^\mpi A^\mpi\).
\erem

\breml{A.R.UA+V}%
 Let \(U\in\Coo{u}{p}\) with \(U^\ad U=\Ip\), let \(V\in\Coo{q}{v}\) with \(VV^\ad=\Iq\), and let \(A\in\Cpq\).
 Then \(\rk{UAV}^\mpi=V^\ad A^\mpi U^\ad\).
\erem

\breml{ab.R1842}%
 Let \(A\in\Cpq\).
 Then \(\ran{A}=\ran{AA^\ad}\) and \(\nul{A}=\nul{A^\ad A}\).
 Furthermore, \(\ran{A^\ad}=\ek{\nul{A}}^\bot=\ran{A^\mpi}\) and \(\nul{A^\ad}=\ek{\ran{A}}^\bot=\nul{A^\mpi}\).
\erem

\breml{R.AA+B}
 Let \(A\in\Cpq\).
 Then:
 \benui
  \il{R.AA+B.a} Let \(B\in\Coo{p}{m}\).
  Then \(\ran{B}\subseteq\ran{A}\) if and only if \(AA^\mpi B=B\).
  \item Let \(C\in\Coo{n}{q}\).
  Then \(\nul{A}\subseteq\nul{C}\) if and only if \(CA^\mpi A=C\).
 \eenui
\erem

\breml{A.R.XAX<=XBX}
 Let \(A,B\in\CHq\) with \(A\lleq B\) and let \(X\in\Cqp\).
 Then the matrices \(X^\ad AX\) and \(X^\ad BX\) are both \tH{} with \(X^\ad AX\lleq X^\ad BX\).
\erem

\breml{A.R.Vsqrt}
 Let \(A\in\Cggq\) and let \(V\in\Cqp\) with \(VV^\ad=\Iq\).
 Then \(V^\ad AV\in\Cggp\) and \(\rk{V^\ad AV}^\varsqrt=V^\ad A^\varsqrt V\).
\erem

\breml{A.R.A+>}%
 If \(A\in\Cggq\), then \(A^\mpi\in\Cggq\) and \(\rk{A^\mpi}^\varsqrt=\rk{A^\varsqrt}^\mpi\).
\erem

\breml{A.R.r-sqrt}
 If \(A\in\Cggq\), then \(\ran{A^\varsqrt}=\ran{A}\) and \(\nul{A^\varsqrt}=\nul{A}\).
\erem

\breml{A.R.0<P<1}
 If \(\mathcal{U}\) is a linear subspace of the unitary space \(\Cp\), then \(\Opp\lleq\OPu{\mathcal{U}}\lleq\Ip\) and \(\OPu{\mathcal{U}}+\OPu{\mathcal{U}^\orth}=\Ip\).
\erem

\breml{A.R.P<P}
 Let \(\mathcal{U}\) and \(\mathcal{V}\) be linear subspaces of the unitary space \(\Cp\).
 Then \(\mathcal{U}\subseteq\mathcal{V}\) if and only if \(\OPu{\mathcal{U}}\lleq\OPu{\mathcal{V}}\).
\erem

\breml{ab.R1052}
 If \(A\in\Cpq\), then \(\OPu{\ran{A}}=AA^\mpi\) and \(\OPu{\ran{A^\ad}}=A^\mpi A\).
\erem

\bleml{A.R.rA<rB}
 Let \(A,B\in\CHq\) with \(\Oqq\lleq A\lleq B\).
 Then:
\benui
 \il{A.R.rA<rB.a} \(\ran{A}\subseteq\ran{B}\) and \(\nul{B}\subseteq\nul{A}\).
 \il{A.R.rA<rB.b} \(\Oqq\lleq\OPu{\ran{A}}B^\mpi\OPu{\ran{A}}\lleq A^\mpi\).
\eenui
\elem
\bproof
 \eqref{A.R.rA<rB.a} This is a well-known fact.%

 \eqref{A.R.rA<rB.b} Let \(Q\defeq A^\varsqrt\) and let \(R\defeq B^\varsqrt\).
 By virtue of \rrem{A.R.XAX<=XBX} and \eqref{mpi}, then \(K\defeq QR^\mpi\) fulfills
 \[
  K^\ad K
  =\rk{R^\mpi}^\ad AR^\mpi
  \lleq\rk{R^\mpi}^\ad BR^\mpi
  =\rk{R^\mpi}^\ad R^\ad RR^\mpi
  =\rk{RR^\mpi}^\ad RR^\mpi
  =RR^\mpi.
 \]
 In view of \rremss{ab.R1052}{A.R.0<P<1}, we have \(RR^\mpi=\OPu{\ran{R}}\lleq\Iq\).
 Thus, \(K^\ad K\lleq\Iq\) follows.
 In particular, \(KK^\ad\lleq\Iq\) holds true (see, \teg{}~\zitaa{MR1152328}{\clem{1.1.12}{19}}).
 Because of \rrem{A.R.A+>}, we have \(B^\mpi\lgeq\Oqq\) and furthermore \(Q^\mpi\rk{Q^\mpi}^\ad=A^\mpi\) and \(R^\mpi\rk{R^\mpi}^\ad=B^\mpi\).
 According to \rrem{A.R.XAX<=XBX}, we obtain then
 \[
  A^\mpi
  \lgeq  Q^\mpi KK^\ad\rk{Q^\mpi}^\ad
  =Q^\mpi QR^\mpi\rk{R^\mpi}^\ad Q^\ad\rk{Q^\mpi}^\ad
  =\rk{Q^\mpi Q} B^\mpi \rk{Q^\mpi Q}^\ad
  \lgeq\Oqq.
 \]
 From \eqref{mpi} and \rremss{ab.R1052}{A.R.r-sqrt}, we get moreover
 \(
  \rk{Q^\mpi Q}^\ad
  =Q^\mpi Q
  =\OPu{\ran{Q^\ad}}
  =\OPu{\ran{Q}}
  =\OPu{\ran{A}}
 \).
\eproof

\bibliography{176arxiv}
\bibliographystyle{bababbrv}

\vfill\noindent
\begin{minipage}{0.5\textwidth}
 Universit\"at Leipzig\\
 Fakult\"at f\"ur Mathematik und Informatik\\
 PF~10~09~20\\
 D-04009~Leipzig
\end{minipage}
\begin{minipage}{0.49\textwidth}
 \begin{flushright}
  \texttt{
   fritzsche@math.uni-leipzig.de\\
   kirstein@math.uni-leipzig.de\\
   maedler@math.uni-leipzig.de
  } 
 \end{flushright}
\end{minipage}

\end{document}